\documentclass[final,3p]{elsarticle}
 \usepackage{graphics}
 \usepackage{graphicx}
 \usepackage{epsfig}
\usepackage{amssymb}
 \usepackage{amsthm}
 \usepackage{lineno}
 \usepackage{amsmath}
   \numberwithin{equation}{section}
\newtheorem{thm}{Theorem}[section]
\newtheorem{cor}[thm]{Corollary}
\newtheorem{lem}[thm]{Lemma}

\newtheorem{defn}[thm]{Definition}
\newtheorem{rem}[thm]{Remark}

\biboptions{sort&compress,square}
\allowdisplaybreaks
\begin{document}
\begin{frontmatter}
\author[rvt1,rvt2]{Jin Hong}
\ead{jhong@nenu.edu.cn}
\author[rvt2]{Yong Wang\corref{cor2}}
\ead{wangy581@nenu.edu.cn}
\cortext[cor2]{Corresponding author.}
\address[rvt1]{School of Mathematics and Statistics, Yili Normal University, Yining 835000, China}
\address[rvt2]{School of Mathematics and Statistics, Northeast Normal University, Changchun, 130024, China}

\title{The $J$-twist of the Dirac operator, spectral Einstein functionals and
	the noncommutative residue}
\begin{abstract}
	We study spectral Einstein functionals associated with the
	$J$-twist of the Dirac operator on almost product Riemannian
	spin manifolds. On closed even-dimensional manifolds, we
	compute the corresponding noncommutative residue in terms
	of the Einstein tensor, the curvature of the connection
	determined by $D_J^2$, and its zero-order endomorphism.
	For four-dimensional compact manifolds with boundary,
	we consider two boundary residue functionals obtained by
	inserting two covariant derivatives in the first factor
	of the compositions
	$\pi^+(D_J^{-2})\circ\pi^+(D_J^{-2})$ and
	$\pi^+(D_J^{-1})\circ\pi^+(D_J^{-3})$.
	Under the stated collar and almost product structure
	assumptions, we calculate their interior and boundary
	contributions and obtain
	D\k{a}browski--Sitarz--Zalecki type formulas.
	We also present applications of these results in physics
	and noncommutative geometry, focusing on spectral actions,
	spectral distance, and mean-curvature boundary contributions.
\end{abstract}
\end{frontmatter}
\section{Introduction}
\label{1}
The noncommutative residue, introduced by Wodzicki and Guillemin
\cite{Wo1,Gu}, is a trace on the algebra of classical
pseudodifferential operators on a closed manifold. Connes' trace theorem relates the noncommutative residue
to the Dixmier trace for classical pseudodifferential
operators of order $-n$ on a closed $n$-dimensional
manifold \cite{Co2}. Its relation
to Riemannian curvature is expressed by the
Kastler--Kalau--Walze theorem: for the Dirac operator $D$ on a
closed even-dimensional Riemannian spin manifold of dimension
$n\geq 4$, the residue $\operatorname{Wres}(D^{-n+2})$ is
proportional to the integral of the scalar curvature
\cite{Ka,KW}. The proof relates the noncommutative residue to
the heat kernel expansion of $D^2$ and the Lichnerowicz formula
\cite{KW,PBG}. For compact manifolds with boundary, Fedosov,
Golse, Leichtnam and Schrohe \cite{FGLS} constructed a
noncommutative residue in Boutet de Monvel's calculus. Their
residue is a trace whose formula consists of an interior integral
and a boundary integral. Grubb and Schrohe \cite{GS} studied the relation between
this boundary residue and resolvent and heat trace
expansions in Boutet de Monvel's calculus. Using this residue, Wang \cite{Wa3,Wa4}
established Kastler--Kalau--Walze type theorems for manifolds
with boundary, including results for the Dirac operator and
the signature operator in dimensions three and four.

D\k{a}browski, Sitarz and Zalecki \cite{DL} introduced spectral
metric and Einstein functionals using the noncommutative residue.
For the scalar Laplacian on a closed even-dimensional Riemannian
manifold, they defined bilinear functionals of vector fields
that are proportional to the integrals of the metric tensor
and the Einstein tensor evaluated on these vector fields.
They also obtained the corresponding formulas for the spin
Laplacian and considered Laplace-type operators on vector
bundles. In the latter case, the spectral Einstein functional
contains contributions from the curvature of the connection
and the zero-order endomorphism in the Laplace-type decomposition.
These constructions were further formulated in the setting
of spectral triples \cite{DL}. Wang, Wang and Wu \cite{WWw}
studied spectral Einstein functionals associated with Dirac
operators with torsion and computed the corresponding
noncommutative residues on four-dimensional manifolds with
boundary. Hong and Yang \cite{HYSE} studied spectral Einstein
functionals associated with sub-Dirac operators and proved
D\k{a}browski--Sitarz--Zalecki type theorems on four-dimensional
foliated manifolds with boundary.

For an almost product Riemannian spin manifold $(M,g,J)$,
the tensor field $J$ satisfies $J^2=\mathrm{Id}$ and
$g(JX,JY)=g(X,Y)$. The $J$-twist of the Dirac operator is
defined by
$D_J=\sum_{j=1}^{n}c(Je_j)\nabla^{S(TM)}_{e_j}$,
where $\{e_j\}_{j=1}^{n}$ is a local orthonormal frame,
$c$ denotes Clifford multiplication, and $\nabla^{S(TM)}$
is the spin connection. Liu and Wang \cite{LWJ} derived a
Lichnerowicz-type formula for $D_J^2$ and proved
Kastler--Kalau--Walze type theorems on closed almost product
Riemannian spin manifolds and on three- and four-dimensional
almost product Riemannian spin manifolds with boundary.
Their formulas involve the Riemann curvature tensor, $J$,
and covariant derivatives of $J$. The Lichnerowicz-type formula
for $D_J^2$, together with the residue formulas for Laplace-type
operators in \cite{DL}, provides the basis for computing the
spectral Einstein functional associated with $D_J^2$.
 
 The aim of this paper is to study spectral Einstein functionals associated with the $J$-twist $D_J$ of the Dirac operator on compact almost product Riemannian spin manifolds, both with and without boundary. In Section 2, we compute the spectral Einstein functional associated with $D_J^2$ on closed even-dimensional manifolds, expressing it in terms of the Einstein tensor and additional geometric terms involving $J$ and its covariant derivatives. In Sections 3 and 4, we compute the noncommutative residues
 $\widetilde{\operatorname{Wres}}
 \left[
 \pi^{+}\left(\widetilde{\nabla}^{J}_{X}\widetilde{\nabla}^{J}_{Y}D_J^{-2}\right)
 \circ\pi^{+}\left(D_J^{-2}\right)
 \right]$
 and
 $\widetilde{\operatorname{Wres}}
 \left[
 \pi^{+}\left(\widetilde{\nabla}^{J}_{X}\widetilde{\nabla}^{J}_{Y}D_J^{-1}\right)
 \circ\pi^{+}\left(D_J^{-3}\right)
 \right]$
 on four-dimensional compact almost product Riemannian spin manifolds with boundary, respectively. 
 Under the collar metric assumptions stated below, we obtain
 D\k{a}browski--Sitarz--Zalecki type formulas for these two
 functionals. Their interior contributions coincide, since
 the compositions before boundary truncation agree, modulo
 smoothing operators, with
 $\widetilde{\nabla}^{J}_{X}
 \widetilde{\nabla}^{J}_{Y}D_J^{-4}$.
 Their boundary contributions generally differ, reflecting
 the dependence of the boundary residue on the factorization
 into the two truncated operators.
 For the second functional, the boundary formula is derived
 under the additional assumption that
 \[
 J(\partial_{x_n})=\partial_{x_n}
 \]
 throughout a collar neighborhood of $\partial M$, where
 $\partial_{x_n}$ is the inward unit normal vector field
 to the collar slices.

\section{Spectral functionals for the $J$-twist of the Dirac operator }

We give some definitions and basic notions which we will use in this paper.

Let $M$ be a smooth compact oriented Riemannian manifold
of dimension $n$ without boundary, and let $N$ be a vector
bundle over $M$.
We say that $P$ is a differential operator of Laplace type, if it has locally the form
\begin{equation}
	P=-(g^{ij}\partial_i\partial_j+A^i\partial_i+B),\nonumber
\end{equation}
where $\partial_{i}$  is a natural local frame on $TM,$ $(g^{ij})_{1\leq i,j\leq n}$ is the inverse matrix associated to the metric
matrix  $(g_{ij})_{1\leq i,j\leq n}$ on $M,$ $A^{i}$ and $B$ are smooth sections of $\textrm{End}(N)$ on $M$ (endomorphism).
If $P$ has the above local form, then there is a unique
connection $\nabla$ on $N$ and a unique endomorphism $E_0$ such that
\begin{equation}
	P=-[g^{ij}(\nabla_{\partial_{i}}\nabla_{\partial_{j}}- \nabla_{\nabla^{L}_{\partial_{i}}\partial_{j}})+E_0],\nonumber
\end{equation}
where $\nabla^{L}$ is the Levi-Civita connection on $M$. Moreover
(with local frames of $T^{*}M$ and $N$), $\nabla_{\partial_{i}}=\partial_{i}+\omega_{i} $
and $E_0$ are related to $g^{ij}$, $A^{i}$ and $B$ through
\begin{eqnarray}
	&&\omega_{i}=\frac{1}{2}g_{ij}\big(A^{j}+g^{kl}\Gamma_{ kl}^{j} \texttt{id}\big),\nonumber\\
	&&E_0=B-g^{ij}\big(\partial_{i}(\omega_{j})+\omega_{i}\omega_{j}-\omega_{k}\Gamma_{ ij}^{k} \big),\nonumber
\end{eqnarray}
where $\Gamma_{ kl}^{j}$ is the  Christoffel coefficient of $\nabla^{L}$.

Let $(M,g,J)$ be a closed oriented almost product
Riemannian spin manifold of even dimension $n=2m\geq4$,
with a fixed spin structure. Let $S(TM)$ denote the
complex spinor bundle and $\nabla^{S(TM)}$ its spin
connection. The tensor field $J$ satisfies
\[
J^2=\mathrm{id},
\qquad
g(JX,JY)=g(X,Y).
\]
We use the Clifford convention
\[
c(X)c(Y)+c(Y)c(X)
=-2g(X,Y)\,\mathrm{id},
\] here $\texttt{id}$ stands for the identity map.

We use the curvature convention
\[
R^L(X,Y)
=
\nabla^L_X\nabla^L_Y
-\nabla^L_Y\nabla^L_X
-\nabla^L_{[X,Y]},
\qquad
R(X,Y,Z,W)=g(R^L(X,Y)W,Z).
\]
The Ricci tensor and scalar curvature are defined by
\[
\operatorname{Ric}(X,Y)
=
\sum_{a=1}^{n}g(R^L(e_a,X)Y,e_a),
\qquad
s=\sum_{a=1}^{n}\operatorname{Ric}(e_a,e_a),
\]
where $\{e_a\}_{a=1}^{n}$ is a local orthonormal frame.
In particular,
\[
\sum_{i,j=1}^{n}R(e_i,e_j,e_j,e_i)=-s.
\]

 We can also define on almost product Riemannian spin manifold the following the $J$-twist of the Dirac operator $D_{J}$ by
\begin{align}
	&D_{J}:=\sum_{j=1}^{n}c(e_{j})\nabla^{S(TM)}_{J(e_{j})}=\sum_{j=1}^{n}c[J(e_{j})]\nabla^{S(TM)}_{e_{j}} .
\end{align}
The connection associated with the Laplace-type operator
$D_J^2$ is
\begin{equation}\label{eq:connection-J}
	\widetilde{\nabla}^{J}_X
	=
	\nabla^{S(TM)}_X
	-\frac12\sum_{\alpha=1}^{n}
	c(Je_\alpha)c((\nabla^L_{e_\alpha}J)X).
\end{equation}
We write
\[
\Delta^{\widetilde{\nabla}^{J}}
=
-\sum_{j=1}^{n}
\left(
\widetilde{\nabla}^{J}_{e_j}
\widetilde{\nabla}^{J}_{e_j}
-
\widetilde{\nabla}^{J}_{\nabla^L_{e_j}e_j}
\right),
\qquad
D_J^2=\Delta^{\widetilde{\nabla}^{J}}+E_J.
\]
For $P=D_J^2$, the two conventions are related by
$E_J=-E_0$.
Negative integer powers of $D_J$ are understood as powers
of a classical parametrix. Choices differing by smoothing
operators yield the same noncommutative residue.

Let $\partial_i=\partial_{x_i}$ and write
\[
\xi=\sum_{j=1}^{n}\xi_j\,dx_j,
\qquad
\nabla^L_{\partial_i}\partial_j
=\sum_{k=1}^{n}\Gamma_{ij}^{k}\partial_k.
\]
For a local orthonormal frame $\{e_k\}_{k=1}^{n}$,
we use the connection one-forms
\[
\omega_{k,l}(Z)=g(\nabla^L_Ze_l,e_k).
\]
In the induced local spinor trivialization,
\[
\nabla^{S(TM)}_{\partial_i}
=\partial_i+\sigma_i,
\qquad
\sigma_i
=-\frac14\sum_{k,l=1}^{n}
\omega_{k,l}(\partial_i)c(e_k)c(e_l).
\]
We use the contracted quantities
\[
\Gamma^k
=\sum_{i,j=1}^{n}g^{ij}\Gamma_{ij}^{k},
\qquad
\sigma^j=\sum_{i=1}^{n}g^{ij}\sigma_i,
\]
and
\[
\xi^j=\sum_{i=1}^{n}g^{ij}\xi_i,
\qquad
\partial^j=\sum_{i=1}^{n}g^{ij}\partial_i.
\]
For a multi-index $\alpha$, we set
\[
D_x^\alpha=(-i)^{|\alpha|}\partial_x^\alpha.
\]

 We recall the following Lichnerowicz-type formula;
see \cite{LWJ}. We use the curvature and Laplace-type
conventions fixed above.

\begin{lem}
	The following identity holds:
\begin{align}\label{3d6}
	D_{J}^{2}
	&=-g^{ij}\left(
	\widetilde{\nabla}^{J}_{\partial_i}
	\widetilde{\nabla}^{J}_{\partial_j}
	-\widetilde{\nabla}^{J}_{\nabla^L_{\partial_i}\partial_j}
	\right)-\frac{1}{8}\sum_{i,j,k,l=1}^{n}R(J(e_{i}), J(e_{j}), e_{k}, e_{l})c(e_{i})c(e_{j})c(e_{k})c(e_{l})\\
	&-\frac{1}{2}\sum_{\nu,j=1}^{n}c[\nabla_{e_{j}}^{L}(J)e_{\nu}]c[(\nabla^{L}_{e_{\nu}}J)e_{j}]-\frac{1}{2}\sum_{\nu,j=1}^{n}c[J(e_{\nu})]c[(\nabla^{L}_{e_{j}}(\nabla^{L}_{e_{\nu}}(J)))e_{j}-(\nabla^{L}_{\nabla^{L}_{e_{j}}e_{\nu}}(J))e_{j}]\nonumber\\
	&+\frac{1}{4}\sum_{\alpha,\nu,j=1}^{n}c[J(e_{\alpha})]c[(\nabla^{L}_{e_{\alpha}}J)e_{j}]c[J(e_{\nu})]c[(\nabla^{L}_{e_{\nu}}J)e_{j}],\nonumber
\end{align}
\end{lem}

We next compute the local coordinate expression of $D_J^3$.
\begin{align}\label{3f1}
	D_{J}^{3}=&\sum^{n}_{r,a=1}c[J(e_{r})]\langle e_{r},dx_{a}\rangle(-g^{ij}\partial_{a}\partial_{i}\partial_{j})+\sum^{n}_{r,a=1}c[J(e_{r})]\langle e_{r},dx_{a}\rangle\Big(-(\partial_{a}g^{ij})\partial_{i}\partial_{j}-g^{ij}(2\sigma_{i}\partial_{j}\nonumber\\
	&\qquad-\Gamma^{k}_{ij}\partial_{k})\partial_{a}+\sum^{n}_{\alpha,\beta,k=1}c[J(e_{\alpha})]c[(\nabla^{L}_{e_{\alpha}}J)e_{\beta}]\langle e_{\beta}, dx_{k}\rangle\partial_{a}\partial_{k}\Big)+\sum_{r,a=1}^{n}
	c[J(e_r)]\langle e_r,dx_a\rangle\sigma_a(-g^{ij}\partial_{i}\partial_{j})\nonumber\\
	&+\sum^{n}_{r,a=1}c[J(e_{r})]\langle e_{r},dx_{a}\rangle\Big[-2(\partial_{a}g^{ij})\sigma_{i}\partial_{j}-2g^{ij}(\partial_{a}\sigma_{i})\partial_{j}+(\partial_{a}g^{ij})\Gamma^{k}_{ij}\partial_{k}+g^{ij}(\partial_{a}\Gamma^{k}_{ij})\partial_{k}\nonumber\\
	&\qquad+\sum^{n}_{\alpha,\beta,k=1}\partial_{a}\Big(c[J(e_{\alpha})]c[(\nabla^{L}_{e_{\alpha}}J)e_{\beta}]\Big)\langle e_{\beta}, dx_{k}\rangle\partial_{k}+\sum^{n}_{\alpha,\beta,k=1}c[J(e_{\alpha})]c[(\nabla^{L}_{e_{\alpha}}J)e_{\beta}]\Big(\partial_{a}\langle e_{\beta}, dx_{k}\rangle\Big)\partial_{k}\Big]\nonumber\\
	&+\sum_{r,a=1}^{n}
	c[J(e_r)]\langle e_r,dx_a\rangle
	\Biggl[
	-\frac18\sum_{i,j,k,l=1}^{n}
	R(J(e_i),J(e_j),e_k,e_l)
	c(e_i)c(e_j)c(e_k)c(e_l)
	\nonumber\\
	&\qquad
	-g^{ij}\Bigl(
	\partial_i\sigma_j+\sigma_i\sigma_j
	-\Gamma_{ij}^{k}\sigma_k
	\Bigr)+\sum_{\alpha,\beta,k=1}^{n}
	c[J(e_\alpha)]
	c[(\nabla_{e_\alpha}^{L}J)e_\beta]
	\langle e_\beta,dx_k\rangle\sigma_k
	\Biggr]\partial_a\nonumber\\
	&+\sum^{n}_{r,a=1}c[J(e_{r})]\langle e_{r},dx_{a}\rangle\partial_{a}\Biggl[
	-\frac18\sum_{i,j,k,l=1}^{n}
	R(J(e_i),J(e_j),e_k,e_l)
	c(e_i)c(e_j)c(e_k)c(e_l)
	\nonumber\\
	&\qquad
	-g^{ij}\Bigl(
	\partial_i\sigma_j+\sigma_i\sigma_j
	-\Gamma_{ij}^{k}\sigma_k
	\Bigr)+\sum_{\alpha,\beta,k=1}^{n}
	c[J(e_\alpha)]
	c[(\nabla_{e_\alpha}^{L}J)e_\beta]
	\langle e_\beta,dx_k\rangle\sigma_k
	\Biggr]\nonumber\\
	&+\sum_{r,a=1}^{n}
	c[J(e_r)]\langle e_r,dx_a\rangle\sigma_a\Big(-2\sigma^{j}\partial_{j}+\Gamma^{k}\partial_{k}-\frac{1}{8}\sum_{i,j,k,l=1}^{n}R(J(e_{i}), J(e_{j}), e_{k}, e_{l})c(e_{i})c(e_{j})c(e_{k})c(e_{l})\nonumber\\
	&\qquad-g^{ij}\big((\partial_{i}\sigma_{j})+\sigma_{i}\sigma_{j}-\Gamma^{k}_{ij}\sigma_{k}\big)+\sum^{n}_{\alpha,\beta,k=1}c[J(e_{\alpha})]c[(\nabla^{L}_{e_{\alpha}}J)e_{\beta}]\langle e_{\beta}, dx_{k}\rangle(\partial_{k}+\sigma_{k})\Big).
\end{align}

To compute the spectral Einstein functional associated
with $D_J^2$, we use the following residue formula for
Laplace-type operators due to D\k{a}browski, Sitarz
and Zalecki \cite{DL}.

Let $V$, $W$ be a pair of vector fields on a compact Riemannian manifold $M$, of dimension $n = 2m$. Consider the Laplace-type operator $\Delta_{T,E}=-\sum_{j=1}^{n}(\nabla^{T}_{e_j}\nabla^{T}_{e_j}-\nabla^{T}_{\nabla_{e_j}^L e_j})+E $
acting on sections of a vector bundle $\overline{E}$, where $\nabla^{T}$ is a connection on $\overline{E}$, and $E$ is an endomorphism of $\overline{E}$. We use the normalization 
$\operatorname{Wres}(P):= \int_{M} \int_{{|\xi|=1}} \operatorname{tr}\sigma_{-n}(P)(x, \xi) \sigma(\xi) \mathrm{d}x$,
where $P$ is a classical pseudodifferential operator.
We define the spectral Einstein functional by
\[
G^{\Delta_{T,E}}(V,W)
:=
\operatorname{Wres}
\left(
\nabla_V^T\nabla_W^T\Delta_{T,E}^{-m}
\right).
\]
The following result holds.
\begin{lem}\cite{DL}
The Einstein functional is equal to
 \begin{equation}
Wres\big(\nabla^{T}_{V}\nabla^{T}_{W}\Delta_{T,E}^{-m}\big)=\frac{\upsilon_{n-1}}{6}
\operatorname{rank}(\overline E)
\int_M G(V,W)\,d\mathrm{vol}_g
 +\frac{\upsilon_{n-1}}{2}\int_{M}F(V,W)vol_{g}+\frac{\upsilon_{n-1}}{2}\int_{M}(\mathrm{Tr}E)g(V,W)vol_{g},
\end{equation}
where $G(V,W)$ denotes the Einstein tensor evaluated on the two vector fields, $F(V,W)=\mathrm{Tr}(V_{a}W_{b}F_{ab})$ and
$F_{ab}$ is the curvature tensor of the connection $\nabla^{T}$, $\mathrm{Tr}E$ denotes the trace of $E$ and $\upsilon_{n-1}=\frac{2\pi^{m}}{\Gamma(m)}$. For $\overline E=S(TM)$, one has
$\operatorname{rank}(\overline E)=2^m$.
\end{lem}
The purpose of this section is to prove the following.

\begin{thm}\label{thm1}
For the Laplace (type) operator $\Delta_{J}=D_{J}^2$, the Einstein functional is equal to

	\begin{align}
		&\operatorname{Wres}\!\left(
		\widetilde{\nabla}^{J}_{V}
		\widetilde{\nabla}^{J}_{W}
		(D_J)^{-2m}
		\right)\nonumber\\[2pt]
		={}&\frac{2^m\upsilon_{n-1}}{6}
		\int_M
		\left[
		\operatorname{Ric}(V,W)-\frac{s}{2}g(V,W)
		\right]\,d\mathrm{vol}_g\nonumber\\[2pt]
		&+2^{m-2}\upsilon_{n-1}
		\int_M\sum_{l=1}^{n}
		\Big\{
		V\!\left[g\!\left(Je_l,(\nabla^L_{e_l}J)W\right)\right]
		-W\!\left[g\!\left(Je_l,(\nabla^L_{e_l}J)V\right)\right]
		-g\!\left(Je_l,(\nabla^L_{e_l}J)[V,W]\right)
		\Big\}\,d\mathrm{vol}_g\nonumber\\[2pt]
		&+2^{m-3}\upsilon_{n-1}\int_M g(V,W)
		\Bigg[
		-\sum_{i,j=1}^{n}R(Je_i,Je_j,e_j,e_i)
		+2\sum_{\nu,j=1}^{n}
		g\!\left(
		(\nabla^L_{e_j}J)e_\nu,
		(\nabla^L_{e_\nu}J)e_j
		\right)\nonumber\\
		&\qquad
		+2\sum_{\nu,j=1}^{n}
		g\!\left(
		Je_\nu,
		\left[
		\nabla^L_{e_j}(\nabla^L_{e_\nu}J)
		-\nabla^L_{\nabla^L_{e_j}e_\nu}J
		\right]e_j
		\right)
		+\sum_{\alpha,\nu,j=1}^{n}
		g\!\left(Je_\alpha,(\nabla^L_{e_\alpha}J)e_j\right)
		g\!\left(Je_\nu,(\nabla^L_{e_\nu}J)e_j\right)\nonumber\\
		&\qquad
		-\sum_{\nu,j=1}^{n}
		g\!\left(
		(\nabla^L_{e_\nu}J)e_j,
		(\nabla^L_{e_\nu}J)e_j
		\right)
		+\sum_{\alpha,\nu,j=1}^{n}
		g\!\left(Je_\alpha,(\nabla^L_{e_\nu}J)e_j\right)
		g\!\left((\nabla^L_{e_\alpha}J)e_j,Je_\nu\right)
		\Bigg]\,d\mathrm{vol}_g .
\end{align}

where $s$ is the scalar curvature, \(\upsilon_{n-1}=2\pi^m/\Gamma(m)\).
\end{thm}

\begin{proof}
According to the definition of connection $\widetilde{\nabla}^{J}$, we have
 \begin{align}
\widetilde{\nabla}_{X}^{J}=&{\nabla}_{X}^{S(TM)}-\frac{1}{2}\sum_{\alpha=1}^{n}c[J(e_{\alpha})]c[(\nabla^{L}_{e_{\alpha}}J)X]\nonumber\\
=&X+\frac{1}{4}\sum_{j,l=1}^{n}\langle {\nabla}^{L}_{X}e_j,e_l\rangle c(e_j)c(e_l)-\frac{1}{2}\sum_{\alpha=1}^{n}c[J(e_{\alpha})]c[(\nabla^{L}_{e_{\alpha}}J)X]\nonumber\\
:=&X+\overline{A}(X).
\end{align}
Fix $x_0\in M$ and choose a local orthonormal frame
$\{e_a\}_{a=1}^{n}$ that is normal at $x_0$, so that
\[
(\nabla^L_{e_a}e_b)(x_0)=0,
\qquad 1\leq a,b\leq n.
\]
Let $V=\sum_{a=1}^{n}V^{a}e_{a}$, $W=\sum_{b=1}^{n}W^{b}e_{b}$,
seeing that
 \begin{equation}
F(V,W)=\mathrm{Tr}(V_{a}W_{b}F_{ab})=\sum_{a,b=1}^{n}V_{a}W_{b}\operatorname{Tr}^{S(TM)}(F_{e_{a},e_{b}}),
\end{equation}
we obtain
 \begin{align}
F_{e_{a},e_{b}}=e_{a}(\overline{A}(e_{b}))-e_{b}(\overline{A}(e_{a}))+\overline{A}(e_{a})\overline{A}(e_{b})-\overline{A}(e_{b})\overline{A}(e_{a})
-\overline{A}([e_{a},e_{b}]).
\end{align}

In addition, a direct calculation yield
 \begin{align}
\mathrm{Tr}\big(e_{a}(\overline{A}(e_{b}))\big)(x_0)
=&\mathrm{Tr}\Big[e_{a}(\frac{1}{4}\sum_{s,t=1}^{n}\langle {\nabla}^{L}_{e_b}e_s,e_t\rangle c(e_s)c(e_t)-\frac{1}{2}\sum_{l=1}^{n}c[J(e_{l})]c[(\nabla^{L}_{e_{l}}J)e_b])](x_0)\nonumber\\
=&\frac{1}{2}e_{a}\sum_{l=1}^{n}g\big(J(e_{l}), (\nabla^{L}_{e_{l}}J)e_b \big)(x_0)\mathrm{Tr}[Id],
\end{align}
Moreover, since the Levi-Civita connection is torsion-free,
we have
\[
[e_a,e_b](x_0)
=
\bigl(\nabla^L_{e_a}e_b-\nabla^L_{e_b}e_a\bigr)(x_0)
=0.
\]
Consequently,
\begin{align}
	\mathrm{Tr}\bigl(\overline{A}([e_a,e_b])\bigr)(x_0)=0.
\end{align}
so
\begin{align}
\sum_{a,b=1}^{n}V_{a}W_{b}\operatorname{Tr}^{S(TM)}(e_{a}(\overline{A}(e_{b})))=&\frac{1}{2}\sum_{a,b,l=1}^{n}V_{a}W_{b}e_{a}g\big(J(e_{l}), (\nabla^{L}_{e_{l}}J)e_b \big)(x_0)\mathrm{Tr}[Id]\nonumber\\
=&\frac{1}{2}\sum_{l=1}^{n}\bigg(V\langle J(e_{l}), (\nabla^{L}_{e_{l}}J)W \rangle - \langle J(e_{l}), (\nabla^{L}_{e_{l}}J)(\nabla^{L}_V W) \rangle \bigg)\mathrm{Tr}[Id],
\end{align}
for the same reason
\begin{align}
	\sum_{a,b=1}^{n}V_{a}W_{b}\operatorname{Tr}^{S(TM)}(e_{b}(\overline{A}(e_{a})))
	=\frac{1}{2}\sum_{l=1}^{n}\bigg(W\langle J(e_{l}), (\nabla^{L}_{e_{l}}J)V \rangle - \langle J(e_{l}), (\nabla^{L}_{e_{l}}J)(\nabla^{L}_W V) \rangle \bigg)\mathrm{Tr}[Id],
\end{align}

so
\begin{align}
	F(V,W)=&\frac{1}{2}\sum_{l=1}^{n}\bigg(V\langle J(e_{l}), (\nabla^{L}_{e_{l}}J)W \rangle - W \langle J(e_{l}), (\nabla^{L}_{e_{l}}J)V \rangle - \langle J(e_{l}), (\nabla^{L}_{e_{l}}J)([V,W]) \rangle \bigg)\mathrm{Tr}[Id].
\end{align}
In the preceding residue formula, we take
$\nabla^T=\widetilde{\nabla}^{J}$ and $E=E_J$.
Taking the fiberwise trace of the zero-order term
in \eqref{3d6}, we obtain

\begin{align}
	\operatorname{Tr}^{S(TM)}E_J=&-\frac{1}{4}\sum_{i,j=1}^{n}R(Je_i,Je_j,e_j,e_i)\mathrm{Tr}[Id]+\frac{1}{2}\sum_{\nu,j=1}^{n}\langle (\nabla^L_{e_j}J)e_\nu,(\nabla^L_{e_\nu}J) e_j\rangle\mathrm{Tr}[Id]\nonumber\\
	&+\frac{1}{2}\sum_{\nu,j=1}^{n}
	\langle  Je_\nu,\left[ \nabla^L_{e_j}(\nabla^L_{e_\nu}J) -\nabla^L_{\nabla^L_{e_j}e_\nu}J \right]e_j  \rangle\mathrm{Tr}[Id]\nonumber\\
	&+\frac {1}{4}\sum_{\alpha,\nu,j=1}^{n}
	\Big[
	\langle Je_\alpha,(\nabla^L_{e_\alpha}J) e_j\rangle
	\langle Je_\nu,(\nabla^L_{e_\nu}J) e_j\rangle -\langle Je_\alpha,Je_\nu\rangle
	\langle (\nabla^L_{e_\alpha}J) e_j,(\nabla^L_{e_\nu}J) e_j\rangle\nonumber\\
	&+\langle Je_\alpha,(\nabla^L_{e_\nu}J) e_j\rangle
	\langle (\nabla^L_{e_\alpha}J) e_j,Je_\nu\rangle
	\Big]\mathrm{Tr}[Id].
\end{align}
which finishes the proof of Theorem \ref{thm1}.
\end{proof}
\begin{rem}
	When $J=\mathrm{id}$, all covariant derivatives of $J$
	vanish, and hence
	\[
	D_J=D,
	\qquad
	\widetilde{\nabla}^{J}=\nabla^{S(TM)}.
	\]
	The curvature symmetries, the first Bianchi identity,
	and the Clifford relations give
	\[
	-\frac18\sum_{i,j,k,l=1}^{n}
	R(e_i,e_j,e_k,e_l)c(e_i)c(e_j)c(e_k)c(e_l)
	=
	\frac{s}{4}\,\mathrm{id}.
	\]
	Thus
	\[
	E_J=\frac{s}{4}\,\mathrm{id}.
	\]
	Since the curvature of the spin connection has
	vanishing fiberwise trace, Theorem~\ref{thm1}
	reduces to
	\begin{align*}
		&\operatorname{Wres}\left(
		\nabla^{S(TM)}_V\nabla^{S(TM)}_W D^{-2m}
		\right)
		\\
		&\qquad=
		2^m\upsilon_{n-1}
		\int_M
		\left(
		\frac16\operatorname{Ric}(V,W)
		+\frac1{24}s\,g(V,W)
		\right)d\mathrm{vol}_g.
	\end{align*}
\end{rem}

 \section{The residue for the $J$-twist of the Dirac operator $\widetilde{\nabla}_{X}^{J}\widetilde{\nabla}_{Y}^{J}D^{-2}_{J}$ and $D^{-2}_{J}$ } 
In this section, we compute the lower dimensional volume for 4-dimension compact manifolds with boundary and get a Dabrowski-Sitarz-Zalecki type formula in this case.
We use Boutet de Monvel's calculus as in \cite{Wa1},
with the geometric and local trivialization conventions
specified below.
 We will consider $D^{2}_{J}$.
Since $[\sigma_{-4}(\widetilde{\nabla}_{X}^{J}\widetilde{\nabla}_{Y}^{J}D^{-2}_{J}\circ D^{-2}_{J})]|_{M}$ has the same expression as $[\sigma_{-4}(\widetilde{\nabla}_{X}^{J}\widetilde{\nabla}_{Y}^{J}D^{-1}_{J}\circ D^{-3}_{J})]|_{M}$ in the case of manifolds without boundary,
we can therefore use Theorem \ref{thm1} to calculate the first term.

Let $(M,g,J)$ be a compact oriented four-dimensional
almost product Riemannian spin manifold with smooth boundary.
Assume that a collar neighborhood of $\partial M$
is identified with $\partial M\times[0,\varepsilon)$
and that
\[
g=h(x_n)^{-1}g^{\partial M}+dx_n^2,
\qquad h(0)=1,
\]
where $h$ is a smooth positive function and $x_n$
increases in the inward normal direction.
Fix a boundary point $x_0$ and choose geodesic coordinates
on $\partial M$ centered at $x_0$.
Let $\{\bar e_a\}_{a=1}^{n-1}$ be a local orthonormal
frame for $g^{\partial M}$, normal at $x_0$, with
$\bar e_a(x_0)=\partial_{x_a}(x_0)$.
Extend this frame independently of $x_n$ and set
\[
e_a=\sqrt{h(x_n)}\,\bar e_a \quad (a<n),
\qquad e_n=\partial_{x_n}.
\]
We use the local spinor trivialization induced by this
orthonormal frame. In this trivialization, the matrices
$c(e_j)$ are constant.

Along the normal coordinate line through $x_0$, we have
\[
c(dx_a)=\sqrt{h(x_n)}\,c(e_a)\quad(a<n),
\qquad c(dx_n)=c(e_n).
\]
Consequently,
\[
\left.\partial_{x_n}c(dx_a)\right|_{x_0}
=\frac{h'(0)}2c(e_a)(x_0)\quad(a<n),
\qquad
\partial_{x_n}c(dx_n)=0.
\]
All partial derivatives of Clifford symbols are taken
in this trivialization, with the coordinate components
of the cotangent variable held fixed.
The vector fields $X$ and $Y$ are smooth up to the boundary.
\begin{cor}
Let $(M,g,J)$ be a closed oriented four-dimensional
almost product Riemannian spin manifold, and let
$\widetilde{\nabla}^{J}$ be the connection defined
in \eqref{eq:connection-J}.
Then, for $X,Y\in\Gamma(TM)$, the following identity holds:
	\begin{align}
		&\operatorname{Wres}(\widetilde{\nabla}_{X}^{J}\widetilde{\nabla}_{Y}^{J}D_{J}^{-2}
		\circ D_{J}^{-2})\nonumber\\
		=&\frac{4\pi^2}{3}
		\int_M
		\left[
		\operatorname{Ric}(X,Y)-\frac{s}{2}g(X,Y)
		\right]\,d\mathrm{vol}_g
		\nonumber\\[2pt]
		&+2\pi^2
		\int_M\sum_{l=1}^{4}
		\Big\{
		X\!\left[g\!\left(Je_l,(\nabla^L_{e_l}J)Y\right)\right]
		-Y\!\left[g\!\left(Je_l,(\nabla^L_{e_l}J)X\right)\right]
		-g\!\left(Je_l,(\nabla^L_{e_l}J)[X,Y]\right)
		\Big\}\,d\mathrm{vol}_g
		\nonumber\\[2pt]
		&+\pi^2\int_M g(X,Y)
		\Bigg[
		-\sum_{i,j=1}^{4}R(Je_i,Je_j,e_j,e_i)
		+2\sum_{\nu,j=1}^{4}
		g\!\left(
		(\nabla^L_{e_j}J)e_\nu,
		(\nabla^L_{e_\nu}J)e_j
		\right)
		\nonumber\\
		&\qquad
		+2\sum_{\nu,j=1}^{4}
		g\!\left(
		Je_\nu,
		\left[
		\nabla^L_{e_j}(\nabla^L_{e_\nu}J)
		-\nabla^L_{\nabla^L_{e_j}e_\nu}J
		\right]e_j
		\right)
		+\sum_{\alpha,\nu,j=1}^{4}
		g\!\left(Je_\alpha,(\nabla^L_{e_\alpha}J)e_j\right)
		g\!\left(Je_\nu,(\nabla^L_{e_\nu}J)e_j\right)
		\nonumber\\
		&\qquad
		-\sum_{\nu,j=1}^{4}
		g\!\left(
		(\nabla^L_{e_\nu}J)e_j,
		(\nabla^L_{e_\nu}J)e_j
		\right)
		+\sum_{\alpha,\nu,j=1}^{4}
		g\!\left(Je_\alpha,(\nabla^L_{e_\nu}J)e_j\right)
		g\!\left((\nabla^L_{e_\alpha}J)e_j,Je_\nu\right)
		\Bigg]\,d\mathrm{vol}_g,
	\end{align}
	where $s$ is the scalar curvature.
\end{cor}
Let  $\mathbf{S}\left(\mathbf{S}^{\prime}\right)$  be the unit sphere about  $\xi\left(\xi^{\prime}\right)$  and  $\sigma(\xi)\left(\sigma\left(\xi^{\prime}\right)\right)$  be the corresponding canonical $n-1(n-2)$ volume form. Now we recall the main theorem in \cite{FGLS},
\begin{thm}\label{thmb0}\cite{FGLS}
	Let $M$ and  $\partial M$  be connected, 
$\operatorname{dim} M=n \geq 3, A=\left(\begin{array}{lr}\pi^{+} P+G & K \\ T & S\end{array}\right) \in \mathcal{B}$ , 
and denote by  $p$, $b$  and  $s$  the local symbols of  $P$, $G$  and  $S$ , respectively. Define

\begin{align}
   \widetilde{\operatorname{Wres}}(A)= 
	&\int_{M} \int_{\mathbf{S}} \operatorname{Tr}_{E}\left[p_{-n}(x, \xi)\right] \sigma(\xi) \mathrm{d} x\nonumber\\
	& +2 \pi \int_{\partial M} \int_{\mathbf{S}^{\prime}}\left\{\operatorname{Tr}_{E}\left[\left(\operatorname{Tr} b_{-n}\right)\left(x^{\prime}, \xi^{\prime}\right)\right]
	+\operatorname{Tr}_{F}\left[s_{1-n}\left(x^{\prime}, \xi^{\prime}\right)\right]\right\} \sigma\left(\xi^{\prime}\right) \mathrm{d} x^{\prime}\nonumber
\end{align}

Then (a)  $\widetilde{\operatorname{Wres}}([A, B])=0$ , for any  $A, B \in \mathcal{B} $; (b) Up to multiplication by a constant,
$\widetilde{\operatorname{Wres}}$ is the unique
continuous trace on
$\mathcal{B}/\mathcal{B}^{-\infty}$.
\end{thm}

Let $p_1,p_2$ be nonnegative integers with
$p_1+p_2\leq n$.
For a classical pseudodifferential operator $\widetilde{A}$,
let $\sigma_l(\widetilde{A})$ denote the homogeneous
component of degree $l$ of its symbol.
Following (3.5) and (3.6) in \cite{Wa1},
we make the following definition.
\begin{defn}
	The spectral Einstein functional associated with
	the $J$-twist of the Dirac operator on a manifold
	with boundary is defined by
	\begin{equation}
		vol_{n}^{\{p_{1},p_{2}\}}M:=\widetilde{\operatorname{Wres}}[\pi^{+}(\widetilde{\nabla}^J_{X}\widetilde{\nabla}^J_{Y}(D_{J}^2)^{-p_{1}})
		\circ\pi^{+}\bigl((D_J^2)^{-p_2}\bigr)],
	\end{equation}
	where $\pi^{+}(\widetilde{\nabla}^J_{X}\widetilde{\nabla}^J_{Y}(D_{J}^2)^{-p_{1}})$, $\pi^{+}\bigl((D_J^2)^{-p_2}\bigr)$ are
	elements in Boutet de Monvel's algebra \cite{Wa3}.
\end{defn}
Applying formula (2.1.4) in \cite{Wa1}
to the above composition, we obtain
\begin{align}
	&\widetilde{\operatorname{Wres}}[\pi^{+}(\widetilde{\nabla}^J_{X}\widetilde{\nabla}^J_{Y}(D_{J}^{2})^{-p_{1}})
	\circ\pi^{+}\bigl((D_J^2)^{-p_2}\bigr)]\nonumber\\
	&=\int_{M}\int_{|\xi|=1}\mathrm{Tr}_{S(TM)}
	[\sigma_{-n}(\widetilde{\nabla}^J_{X}\widetilde{\nabla}^J_{Y}(D_{J}^2)^{-p_{1}}
	\circ (D_{J}^{2})^{-p_{2}})]\sigma(\xi)\texttt{d}x+\int_{\partial M}\Phi,
\end{align}
where
\begin{align}\label{c4}
	\Phi=&\int_{|\xi'|=1}\int_{-\infty}^{+\infty}\sum_{j,k=0}^{\infty}\sum \frac{(-i)^{|\alpha|+j+k+1}}{\alpha!(j+k+1)!}
	\mathrm{Tr}_{S(TM)}[\partial_{x_{n}}^{j}\partial_{\xi'}^{\alpha}\partial_{\xi_{n}}^{k}\sigma_{r}^{+}
	(\widetilde{\nabla}^J_{X}\widetilde{\nabla}^J_{Y}(D_{J}^{2})^{-p_{1}})(x',0,\xi',\xi_{n})\nonumber\\
	&\times\partial_{x'}^{\alpha}\partial_{\xi_{n}}^{j+1}\partial_{x_{n}}^{k}\sigma_{l}((D_{J}^2)^{-p_{2}})(x',0,\xi',\xi_{n})]
	\texttt{d}\xi_{n}\sigma(\xi')\texttt{d}x' ,
\end{align}
and the sum is taken over $r-k-|\alpha|+\ell-j-1=-n,r\leq2-2p_{1},\ell\leq-2p_{2}$.\\

So we only need to compute $\int_{\partial M}\Phi$.
 
 In the spinor trivialization fixed above, the connection
 \eqref{eq:connection-J} has the local expression
 \[
 \widetilde{\nabla}^J_{\mathbb X}
 =
 \mathbb X+M(\mathbb X)+N(\mathbb X),
 \]
 where
 \begin{align*}
 	M(\mathbb X)
 	&=
 	\frac14\sum_{s,t=1}^{n}
 	\langle\nabla^L_{\mathbb X}e_s,e_t\rangle
 	c(e_s)c(e_t),\\
 	N(\mathbb X)
 	&=
 	-\frac12\sum_{\alpha=1}^{n}
 	c(Je_\alpha)c((\nabla^L_{e_\alpha}J)\mathbb X).
 \end{align*}
 Expanding the composition, we obtain
 
\begin{align}
\widetilde{\nabla}^J_{\mathbb{X}}\widetilde{\nabla}^J_{\mathbb{Y}}&=(\mathbb{X}+M(\mathbb{X})+N(\mathbb{X}))(\mathbb{Y}+M(\mathbb{Y})+N(\mathbb{Y}))\nonumber\\
 &=\mathbb{X}\mathbb{Y}+M(\mathbb{Y})\mathbb{X}+N(\mathbb{Y})\mathbb{X}+M(\mathbb{X})\mathbb{Y}+N(\mathbb{X})\mathbb{Y}+\mathbb{X}[M(\mathbb{Y})]+\mathbb{X}[N(\mathbb{Y})]\nonumber\\
 &~~~~+M(\mathbb{X})M(\mathbb{Y})+M(\mathbb{X})N(\mathbb{Y})+N(\mathbb{X})M(\mathbb{Y})+N(\mathbb{X})N(\mathbb{Y}),
\end{align}
where
$\mathbb{X}=\sum_{j=1}^nX_j\partial_{x_j}, \mathbb{Y}=\sum_{l=1}^nY_l\partial_{x_l}$.

Then we have the following lemmas.
\begin{lem}\label{lem3} The following identities hold:
\begin{align}
 \sigma_{0}(\widetilde{\nabla}^J_{X}\widetilde{\nabla}^J_{Y})=&X[M(Y)]+X[N(Y)]+M(X)M(Y)+M(X)N(Y)+N(X)M(Y)+N(X)N(Y);\\
\sigma_{1}(\widetilde{\nabla}^J_{X}\widetilde{\nabla}^J_{Y})
=&\sqrt{-1}\sum_{j,l=1}^nX_j\frac{\partial{Y_l}}{\partial{x_j}}\xi_l
+\sqrt{-1}\sum_{j=1}^{n}[M(Y)+N(Y)]X_j\xi_j+\sqrt{-1}\sum_{l=1}^{n}[M(X)+N(X)]Y_l\xi_l;\\
\sigma_{2}(\widetilde{\nabla}^J_{X}\widetilde{\nabla}^J_{Y})=&-\sum_{j,l=1}^nX_jY_l\xi_j\xi_l.
\end{align}
\end{lem}

\begin{lem} \label{lem4}\cite{WW} The following identities hold:
\begin{align}
\sigma_{-1}({D}_{J}^{-1})=&\frac{ic[J(\xi)]}{|\xi|^2};\\
\sigma_{-2}({D}_{J}^{-1})=&\frac{c[J(\xi)]\sigma_{0}({D}_{J})c[J(\xi)]}{|\xi|^4}+\frac{c[J(\xi)]}{|\xi|^6}\sum_ {j=1}^{n} c[J(dx_j)]
\Big[\partial_{x_j}(c[J(\xi)])|\xi|^2-c[J(\xi)]\partial_{x_j}(|\xi|^2)\Big];\\
\sigma_{-2}\left(D_{J}^{-2}\right)= & |\xi|^{-2};\\
\sigma_{-3}\left(D_{J}^{-2}\right)= & -\sqrt{-1}|\xi|^{-4}\sum_{k=1}^{n} \xi_{k}\bigg(\Gamma^{k}-2\sigma^{k}+\sum_{\alpha=1}^{n}c[J(e_{\alpha})]c[(\nabla^{L}_{e_{\alpha}}J)(dx_k)^{*}]\bigg)\nonumber\\
&-2\sqrt{-1}|\xi|^{-6} \sum_{j,\alpha,\beta=1}^{n}\xi^{j} \xi_{\alpha} \xi_{\beta} \partial_{j} g^{\alpha \beta}.
\end{align}
where $\sigma_0({D}_{J})=-\frac{1}{4}\sum_{i,j,k=1}^{n}\omega_{j,k}(e_i)c[J(e_i)]c(e_j)c(e_k)$.
\end{lem}

\begin{lem} The following identities hold:
\begin{align}
\sigma_{0}(\widetilde{\nabla}^J_{X} \widetilde{\nabla}^J_{Y} D^{-2}_{J})=&
-\sum_{j,l=1}^nX_jY_l\xi_j\xi_l|\xi|^{-2};\\
\sigma_{-1}(\widetilde{\nabla}^J_{X} \widetilde{\nabla}^J_{Y} D^{-2}_{J})=&
\sigma_{2}(\widetilde{\nabla}^J_{X}\widetilde{\nabla}^J_{Y})\sigma_{-3}(D^{-2}_{J})
+\sigma_{1}(\widetilde{\nabla}^J_{X}\widetilde{\nabla}^J_{Y})\sigma_{-2}(D^{-2}_{J})\nonumber\\
&+\sum_{j=1}^{n}\partial_{\xi_{j}}\big[\sigma_{2}(\widetilde{\nabla}^J_{X}\widetilde{\nabla}^J_{Y})\big]
D_{x_{j}}\big[\sigma_{-2}(D^{-2}_{J})\big].
\end{align}
\end{lem}

Now we  need to compute $\int_{\partial M} \Phi$. When $n=4$, then ${\rm tr}[{\rm \texttt{id}}]=4$, the sum is taken over $
r+l-k-j-|\alpha|=-3,~~r\leq 0,~~l\leq-2,$ then we have the following five cases:

 {\bf case a)~I)}~$r=0,~l=-2,~k=j=0,~|\alpha|=1$.

 By (3.4), we get
\begin{equation}
\label{b24}
\Phi_1=-\int_{|\xi'|=1}\int^{+\infty}_{-\infty}\sum_{|\alpha|=1}
\mathrm{Tr}[\partial^\alpha_{\xi'}\pi^+_{\xi_n}\sigma_{0}(\widetilde{\nabla}^J_{X}\widetilde{\nabla}^J_{Y} D^{-2}_{J})\times
 \partial^\alpha_{x'}\partial_{\xi_n}\sigma_{-2}(D^{-2}_{J})](x_0)d\xi_n\sigma(\xi')dx'.
\end{equation}
By Lemma 2.2 in \cite{Wa3}, for $\beta<n$, then
\begin{equation}
\label{b25}
\partial_{x_\beta}\sigma_{-2}(D^{-2}_{J})(x_0)=
\partial_{x_\beta}(|\xi|^{-2})(x_0)=
-\frac{\partial_{x_\beta}(|\xi|^{2})(x_0)}{|\xi|^4}=0,
\end{equation}
 so $\Phi_1=0$.

 {\bf case a)~II)}~$r=0,~l=-2,~k=|\alpha|=0,~j=1$.

 By (3.4), we get
\begin{equation}
\label{b26}
\Phi_2=-\frac{1}{2}\int_{|\xi'|=1}\int^{+\infty}_{-\infty}
\mathrm{Tr}[\partial_{x_n}\pi^+_{\xi_n}\sigma_{0}(\widetilde{\nabla}^J_{X}\widetilde{\nabla}^J_{Y} D_{J}^{-2})\times
\partial_{\xi_n}^2\sigma_{-2}(D^{-2}_{J})](x_0)d\xi_n\sigma(\xi')dx'.
\end{equation}
By Lemma 3.5, we have
\begin{eqnarray}\label{b237}
\partial_{\xi_n}^2\sigma_{-2}(D^{-2}_{J})(x_0)=\partial_{\xi_n}^2(|\xi|^{-2})(x_0)=\frac{6\xi_n^2-2}{(1+\xi_n^2)^3}.
\end{eqnarray}
It follows that
\begin{align}
\partial_{x_n}\sigma_0
	\left(
	\widetilde{\nabla}^J_X
	\widetilde{\nabla}^J_Y
	D_J^{-2}
	\right)(x_0)
	|_{|\xi'|=1}
	={}&
	\partial_{x_n}\left(
	-\sum_{j,l=1}^{n}
	X_jY_l\xi_j\xi_l|\xi|^{-2}
	\right)(x_0)
	|_{|\xi'|=1}
	\notag\\
	={}&
	\frac{h'(0)}{(1+\xi_n^2)^2}
	\sum_{j,l=1}^{n}
	X_jY_l\xi_j\xi_l-
	\frac{1}{1+\xi_n^2}
	\sum_{j,l=1}^{n}
	\partial_{x_n}(X_jY_l)\xi_j\xi_l .
\end{align}
By integrating formula, we obtain
\begin{align}
	\pi_{\xi_n}^{+}\partial_{x_n}
	\sigma_0\!\left(
	\widetilde{\nabla}^J_X
	\widetilde{\nabla}^J_Y
	D_J^{-2}
	\right)(x_0)
	|_{|\xi'|=1}
	={}&
	\partial_{x_n}\pi_{\xi_n}^{+}
	\sigma_0\!\left(
	\widetilde{\nabla}^J_X
	\widetilde{\nabla}^J_Y
	D_J^{-2}
	\right)(x_0)
	|_{|\xi'|=1}
	\notag\\
	={}&
	-\frac{(i\xi_n+2)h'(0)}{4(\xi_n-i)^2}
	\sum_{j,l=1}^{n-1}X_jY_l\xi_j\xi_l
	-\frac{i\xi_nh'(0)}{4(\xi_n-i)^2}X_nY_n
	\notag\\
	&-
	\frac{ih'(0)}{4(\xi_n-i)^2}
	\sum_{j=1}^{n-1}X_jY_n\xi_j
	-
	\frac{ih'(0)}{4(\xi_n-i)^2}
	\sum_{l=1}^{n-1}X_nY_l\xi_l
	\notag\\
	&+
	\frac{i}{2(\xi_n-i)}
	\sum_{j,l=1}^{n-1}
	\partial_{x_n}(X_jY_l)\xi_j\xi_l
	-
	\frac{i}{2(\xi_n-i)}
	\partial_{x_n}(X_nY_n)
	\notag\\
	&-
	\frac{1}{2(\xi_n-i)}
	\sum_{j=1}^{n-1}
	\partial_{x_n}(X_jY_n)\xi_j
	-
	\frac{1}{2(\xi_n-i)}
	\sum_{l=1}^{n-1}
	\partial_{x_n}(X_nY_l)\xi_l .
\end{align}
We note that $\beta<n,~\int_{|\xi'|=1}\xi_{\beta_{1}}\xi_{\beta_{2}}\cdots\xi_{\beta_{2d+1}}\sigma(\xi')=0$,
so we omit some items that have no contribution for computing {\bf case a)~II)}.
From (3.18) and (3.20), we obtain
\begin{align}
	&\left.
	\operatorname{Tr}\!\left[
	\partial_{x_n}\pi_{\xi_n}^{+}
	\sigma_0\!\left(
	\widetilde{\nabla}^J_X
	\widetilde{\nabla}^J_Y
	D_J^{-2}
	\right)
	\partial_{\xi_n}^{2}\sigma_{-2}(D_J^{-2})
	\right](x_0)
	\right|_{|\xi'|=1}
	\notag\\
	\sim{}&
	2\left[
	\frac{i-3i\xi_n^2}{(\xi_n-i)^4(\xi_n+i)^3}
	+
	\frac{1-3\xi_n^2}{(\xi_n-i)^5(\xi_n+i)^3}
	\right]
	\sum_{j,l=1}^{n-1}X_jY_l\xi_j\xi_l\,h'(0)
	\notag\\
	&+
	2\left[
	\frac{i-3i\xi_n^2}{(\xi_n-i)^4(\xi_n+i)^3}
	-
	\frac{1-3\xi_n^2}{(\xi_n-i)^5(\xi_n+i)^3}
	\right]
	X_nY_n\,h'(0)
	\notag\\
	&+
	\frac{4i(3\xi_n^2-1)}
	{(\xi_n-i)^4(\xi_n+i)^3}
	\sum_{j,l=1}^{n-1}
	\partial_{x_n}(X_jY_l)\xi_j\xi_l
	-
	\frac{4i(3\xi_n^2-1)}
	{(\xi_n-i)^4(\xi_n+i)^3}
	\partial_{x_n}(X_nY_n).
\end{align}
Therefore, we get
\begin{align}
	\Phi_2
	={}&
	\Biggl[
	\left(
	\frac{5}{48}\sum_{j=1}^{n-1}X_jY_j
	-\frac{1}{16}X_nY_n
	\right)h'(0)
	\notag\\
	&\qquad
	-\frac{1}{12}\sum_{j=1}^{n-1}
	\partial_{x_n}(X_jY_j)
	+\frac14\partial_{x_n}(X_nY_n)
	\Biggr]\pi\Omega_3\,dx'.
	\notag
\end{align}
where $\Omega_3=\operatorname{Vol}(S^2)=4\pi$.

  {\bf case a)~III)}~$r=0,~l=-2,~j=|\alpha|=0,~k=1$.

By (3.4), we get
\begin{align}\label{36}
\Phi_3&=-\frac{1}{2}\int_{|\xi'|=1}\int^{+\infty}_{-\infty}
\mathrm{Tr} [\partial_{\xi_n}\pi^+_{\xi_n}\sigma_{0}(\widetilde{\nabla}^J_{X}\widetilde{\nabla}^J_{Y} D_{J}^{-2})\times
\partial_{\xi_n}\partial_{x_n}\sigma_{-2}(D_{J}^{-2})](x_0)d\xi_n\sigma(\xi')dx'\nonumber\\
&=\frac{1}{2}\int_{|\xi'|=1}\int^{+\infty}_{-\infty}
\mathrm{Tr}[\partial_{\xi_n}^2\pi^+_{\xi_n}\sigma_{0}(\widetilde{\nabla}^J_{X}\widetilde{\nabla}^J_{Y} D_{J}^{-2})\times
\partial_{x_n}\sigma_{-2}(D_{J}^{-2})](x_0)d\xi_n\sigma(\xi')dx'.
\end{align}
 By Lemma 3.5, we have
\begin{eqnarray}\label{37}
\partial_{x_n}\sigma_{-2}(D_{J}^{-2})(x_0)|_{|\xi'|=1}
=-\frac{h'(0)}{(1+\xi_n^2)^2}.
\end{eqnarray}
An easy calculation gives
\begin{align}\label{38}
\pi^+_{\xi_n}\sigma_{0}(\widetilde{\nabla}^J_{X}\widetilde{\nabla}^J_{Y} D^{-2}_{J})(x_0)|_{|\xi'|=1}
=&\frac{i}{2(\xi_n-i)}\sum_{j,l=1}^{n-1}X_jY_l\xi_j\xi_l-\frac{i}{2(\xi_n-i)}X_nY_n\nonumber\\
&-\frac{1}{2(\xi_n-i)}\sum_{j=1}^{n-1}X_jY_n\xi_j-\frac{1}{2(\xi_n-i)}\sum_{l=1}^{n-1}X_nY_l\xi_l.
\end{align}
Also, straightforward computations yield
\begin{align}\label{mmmmm}
\partial_{\xi_n}^2\pi^+_{\xi_n}\sigma_{0}(\widetilde{\nabla}^J_{X}\widetilde{\nabla}^J_{Y} D^{-2}_{J})(x_0)|_{|\xi'|=1}
=&\frac{i}{(\xi_n-i)^3}\sum_{j,l=1}^{n-1}X_jY_l\xi_j\xi_l-\frac{i}{(\xi_n-i)^3}X_nY_n\nonumber\\
&-\frac{1}{(\xi_n-i)^3}\sum_{j=1}^{n-1}X_jY_n\xi_j-\frac{1}{(\xi_n-i)^3}\sum_{l=1}^{n-1}X_nY_l\xi_l.
\end{align}

Therefore, we get
\begin{align}\label{41}
\Phi_3&=\frac{1}{2}\int_{|\xi'|=1}\int^{+\infty}_{-\infty}
\mathrm{Tr}[\partial^2_{\xi_n}\pi^+_{\xi_n}\sigma_{0}(\widetilde{\nabla}^J_{X}\widetilde{\nabla}^J_{Y} D^{-2}_{J})\times
\partial_{x_n}\sigma_{-2}(D^{-2}_{J})](x_0)d\xi_n\sigma(\xi')dx'\nonumber\\
&=\left(-\frac{5}{48}\sum_{j=1}^{n-1}X_jY_j+\frac{5}{16}X_nY_n\right)h'(0)\pi\Omega_3dx'.
\end{align}

 {\bf case b)}~$r=0,~l=-3,~k=j=|\alpha|=0$.

 By (3.4), we get
\begin{align}\label{42}
\Phi_4&=-i\int_{|\xi'|=1}\int^{+\infty}_{-\infty}\mathrm{Tr} [\pi^+_{\xi_n}
\sigma_{0}(\widetilde{\nabla}^J_{X}\widetilde{\nabla}^J_{Y} D^{-2}_{J})\times
\partial_{\xi_n}\sigma_{-3}(D^{-2}_{J})](x_0)d\xi_n\sigma(\xi')dx'\nonumber\\
&=i\int_{|\xi'|=1}\int^{+\infty}_{-\infty}\mathrm{Tr} [\partial_{\xi_n}\pi^+_{\xi_n}
\sigma_{0}(\widetilde{\nabla}^J_{X}\widetilde{\nabla}^J_{Y} D^{-2}_{J})\times
\sigma_{-3}(D^{-2}_{J})](x_0)d\xi_n\sigma(\xi')dx'.
\end{align}
 By Lemma 3.5, we have
\begin{align}\label{43}
\sigma_{-3}\left(D_{J}^{-2}\right)(x_0)|_{|\xi'|=1}=& \frac{-2 i h^{\prime}(0) \xi_{n}}{\left(1+\xi_{n}^{2}\right)^{3}}-\frac{i}{\left(1+\xi_{n}^{2}\right)^{2}}\bigg(\frac{3}{2}h^{\prime}(0) \xi_{n}-\frac{1}{2} \sum_{k=1}^{n-1} h^{\prime}(0) c\left(e_{k}\right) c\left(e_{n}\right) \xi_{k}\nonumber\\
&+\sum_{\alpha=1}^{n}\sum_{k=1}^{n-1}c[J(e_{\alpha})]c[(\nabla^{L}_{e_{\alpha}}J)(e_k)]\xi_{k}+\sum_{\alpha=1}^{n}c[J(e_{\alpha})]c[(\nabla^{L}_{e_{\alpha}}J)(e_n)]\xi_{n}\bigg) .
\end{align}

\begin{align}\label{45}
\partial_{\xi_n}\pi^+_{\xi_n}\sigma_{0}(\widetilde{\nabla}^J_{X}\widetilde{\nabla}^J_{Y} D^{-2}_{J})(x_0)|_{|\xi'|=1}
=&-\frac{i}{2(\xi_n-i)^2}\sum_{j,l=1}^{n-1}X_jY_l\xi_j\xi_l+\frac{i}{2(\xi_n-i)^2}X_nY_n\nonumber\\
&+\frac{1}{2(\xi_n-i)^2}\sum_{j=1}^{n-1}X_jY_n\xi_j+\frac{1}{2(\xi_n-i)^2}\sum_{l=1}^{n-1}X_nY_l\xi_l.
\end{align}
We note that $\beta<n,~\int_{|\xi'|=1}\xi_{\beta_{1}}\xi_{\beta_{2}}\cdots\xi_{\beta_{2d+1}}\sigma(\xi')=0$,
so we omit some items that have no contribution for computing {\bf case b)}.
By the trace identity  $\operatorname{Tr}(A B)=\operatorname{Tr}(B A)$ and the relation of the Clifford action, we get
\begin{align}
	\Phi_4
	={}&
	i\int_{|\xi'|=1}\int_{-\infty}^{+\infty}
	\operatorname{Tr}\!\left[
	\partial_{\xi_n}\pi_{\xi_n}^{+}
	\sigma_0\!\left(
	\widetilde{\nabla}^J_X
	\widetilde{\nabla}^J_Y
	D_J^{-2}
	\right)
	\sigma_{-3}(D_J^{-2})
	\right](x_0)
	\,d\xi_n\,\sigma(\xi')\,dx'
	\notag\\
	={}&
	\left(
	\frac{11}{48}\sum_{j=1}^{n-1}X_jY_j
	-\frac{11}{16}X_nY_n
	\right)
	h'(0)\pi\Omega_3\,dx'
	\notag\\
	&+
	\frac16
	\sum_{\alpha=1}^{n}\sum_{j=1}^{n-1}
	\left(X_jY_n+X_nY_j\right)
	g\!\left(
	J(e_\alpha),(\nabla^L_{e_\alpha}J)(e_j)
	\right)
	\pi\Omega_3\,dx'
	\notag\\
	&+
	\left(
	-\frac1{12}\sum_{j=1}^{n-1}X_jY_j
	+\frac14X_nY_n
	\right)
	\sum_{\alpha=1}^{n}
	g\!\left(
	J(e_\alpha),(\nabla^L_{e_\alpha}J)(e_n)
	\right)
	\pi\Omega_3\,dx'.
	\notag
\end{align}

 {\bf  case c)}~$r=-1,~\ell=-2,~k=j=|\alpha|=0$.

By (3.4), we get
\begin{align}\label{61}
\Phi_5=-i\int_{|\xi'|=1}\int^{+\infty}_{-\infty}\mathrm{Tr} [\pi^+_{\xi_n}
\sigma_{-1}(\widetilde{\nabla}^J_{X}\widetilde{\nabla}^J_{Y} D^{-2}_{J})\times
\partial_{\xi_n}\sigma_{-2}(D_{J}^{-2})](x_0)d\xi_n\sigma(\xi')dx'.
\end{align}
By Lemma 3.5, we have
\begin{align}\label{62}
\partial_{\xi_n}\sigma_{-2}(D_{J}^{-2})(x_0)|_{|\xi'|=1}=-\frac{2\xi_n}{(\xi_n^2+1)^2}.
\end{align}
Since
\begin{align}\label{3.38}
\sigma_{-1}(\widetilde{\nabla}^J_{X}\widetilde{\nabla}^J_{Y} D^{-2}_{J})(x_0)|_{|\xi'|=1}
=&\sigma_{2}(\widetilde{\nabla}^J_{X}\widetilde{\nabla}^J_{Y})\sigma_{-3}(D_{J}^{-2})
+\sigma_{1}(\widetilde{\nabla}^J_{X}\widetilde{\nabla}^J_{Y})\sigma_{-2}(D_{J}^{-2})\nonumber\\
&+\sum_{j=1}^{n}\partial_{\xi_{j}}\big[\sigma_{2}(\widetilde{\nabla}^J_{X}\widetilde{\nabla}^J_{Y})\big]
D_{x_{j}}\big[\sigma_{-2}(D_{J}^{-2})\big].
\end{align}
 Explicit representation the first item of \eqref{3.38},
\begin{align}
&\sigma_{2}(\widetilde{\nabla}^J_{X}\widetilde{\nabla}^J_{Y} )\sigma_{-3}(D_{J}^{-2})(x_0)|_{|\xi'|=1}\nonumber\\
=&-\sum_{j,l=1}^{n}X_jY_l\xi_j\xi_l\times \bigg \{\frac{-2 i h^{\prime}(0) \xi_{n}}{\left(1+\xi_{n}^{2}\right)^{3}}-\frac{i}{\left(1+\xi_{n}^{2}\right)^{2}}\bigg(\frac{3}{2}h^{\prime}(0) \xi_{n}-\frac{1}{2} \sum_{k=1}^{n-1} h^{\prime}(0) c\left(e_{k}\right) c\left(e_{n}\right) \xi_{k}\nonumber\\
&+\sum_{\alpha=1}^{n}\sum_{k=1}^{n-1}c[J(e_{\alpha})]c[(\nabla^{L}_{e_{\alpha}}J)(e_k)]\xi_{k}+\sum_{\alpha=1}^{n}c[J(e_{\alpha})]c[(\nabla^{L}_{e_{\alpha}}J)(e_n)]\xi_{n}\bigg) \bigg \}.
\end{align}
Explicit representation the second item of \eqref{3.38},
\begin{align}
&\sigma_{1}(\widetilde{\nabla}^J_{X}\widetilde{\nabla}^J_{Y}) \sigma_{-2}(D_{J}^{-2})(x_0)|_{|\xi'|=1}\nonumber\\
=&\frac{i}{1+\xi_{n}^{2}}\Big(\sum_{j,l=1}^{n}X_j\frac{\partial{Y_l}}{\partial{x_j}}\xi_l
+\sum_{j=1}^{n}\big[M(Y)+N(Y)\big]X_j\xi_j+\sum_{l=1}^{n}\big[M(X)+N(X)\big]Y_l\xi_l\Big).
\end{align}
Explicit representation the third item of \eqref{3.38},
\begin{align}
&\sum_{j=1}^{n}\partial_{\xi_{j}}\big[\sigma_{2}(\widetilde{\nabla}^J_{X}\widetilde{\nabla}^J_{Y})\big]
D_{x_{j}} \big[ \sigma_{-2}(D_{J}^{-2}) \big] (x_0)|_{|\xi'|=1} \nonumber\\
=&\sum_{j=1}^{n}\partial_{\xi_j}
\biggl[-\sum_{a,l=1}^{n}X_aY_l\xi_a\xi_l\biggr]
(-i)\partial_{x_j}\bigl(|\xi|^{-2}\bigr)\nonumber\\
=&-\frac{\sqrt{-1}}{(1+\xi_{n}^{2})^2}h^{\prime}(0)\sum_{l=1}^{n}(X_{n}Y_l+X_{l}Y_n)\xi_{l}.
\end{align}
We note that $\beta<n,~\int_{|\xi'|=1}\xi_{\beta_{1}}\xi_{\beta_{2}}\cdots\xi_{\beta_{2d+1}}\sigma(\xi')=0$,
so we omit some items that have no contribution for computing {\bf case c)}.
Also, straightforward computations yield

\begin{align}
	\Phi_5
	={}&
	\Biggl[
	\left(
	-\frac{11}{48}\sum_{j=1}^{n-1}X_jY_j
	-\frac1{16}X_nY_n
	\right)h'(0)
	-\frac12\sum_{j=1}^{n}
	X_j\frac{\partial Y_n}{\partial x_j}
	\notag\\
	&\qquad
	+\left(
	\frac1{12}\sum_{j=1}^{n-1}X_jY_j
	-\frac14X_nY_n
	\right)
	\sum_{\alpha=1}^{n}
	g\!\left(
	J(e_\alpha),(\nabla^L_{e_\alpha}J)(e_n)
	\right)
	\notag\\
	&\qquad
	-\frac16
	\sum_{\alpha=1}^{n}\sum_{j=1}^{n-1}
	(X_jY_n+X_nY_j)
	g\!\left(
	J(e_\alpha),(\nabla^L_{e_\alpha}J)(e_j)
	\right)
	\Biggr]\pi\Omega_3\,dx'.
\end{align}
Now $\Phi$ is the sum of the cases (a), (b) and (c). 

Write
$X=X^T+X_n\partial_{x_n}$ and
$Y=Y^T+Y_n\partial_{x_n}$,
where $X^T$ and $Y^T$ are the tangential components
with respect to the collar coordinates.
Since
$g=h(x_n)^{-1}g^{\partial M}+dx_n^2$
and $h(0)=1$, we have
\begin{align*}
	\sum_{j=1}^{n-1}(X_jY_j)(x_0)
	&=g(X^T,Y^T)(x_0),\\
	\sum_{j=1}^{n-1}
	\partial_{x_n}(X_jY_j)(x_0)
	&=
	\left.
	\partial_{x_n}\bigl(g(X^T,Y^T)\bigr)
	\right|_{x_0}
	+h'(0)g(X^T,Y^T)(x_0).
\end{align*}
Using these identities and summing the contributions
of the five cases, we obtain
\begin{align}
	\Phi
	=\sum_{\beta=1}^{5}\Phi_\beta
	={}&
	\Biggl[
	-\frac12h'(0)X_nY_n
	-\frac12\sum_{j=1}^{n}
	X_j\frac{\partial Y_n}{\partial x_j}
	\notag\\
	&\qquad
	+\frac14\partial_{x_n}(X_nY_n)
	-\frac1{12}\sum_{j=1}^{n-1}
	\partial_{x_n}(X_jY_j)
	\Biggr]\pi\Omega_3\,dx'.
	\notag
\end{align}
Then we obtain following theorem
\begin{thm}\label{thmb1}
Under the geometric assumptions of this section, let
$\widetilde{\nabla}^{J}$ be the connection defined
in \eqref{eq:connection-J}. Then, for
$X,Y\in\Gamma(TM)$, the following identity holds:
\begin{align}
&\widetilde{\operatorname{Wres}}[\pi^+(\widetilde{\nabla}_{X}^{J}\widetilde{\nabla}_{Y}^{J}D_{J}^{-2})\circ\pi^+(D_{J}^{-2})]\nonumber\\
=&\frac{4\pi^2}{3}
\int_M
\left[
\operatorname{Ric}(X,Y)-\frac{s}{2}g(X,Y)
\right]\,d\mathrm{vol}_g
\nonumber\\[2pt]
&+2\pi^2
\int_M\sum_{l=1}^{4}
\Big\{
X\!\left[g\!\left(Je_l,(\nabla^L_{e_l}J)Y\right)\right]
-Y\!\left[g\!\left(Je_l,(\nabla^L_{e_l}J)X\right)\right]
-g\!\left(Je_l,(\nabla^L_{e_l}J)[X,Y]\right)
\Big\}\,d\mathrm{vol}_g
\nonumber\\[2pt]
&+\pi^2\int_M g(X,Y)
\Bigg[
-\sum_{i,j=1}^{4}R(Je_i,Je_j,e_j,e_i)
+2\sum_{\nu,j=1}^{4}
g\!\left(
(\nabla^L_{e_j}J)e_\nu,
(\nabla^L_{e_\nu}J)e_j
\right)
\nonumber\\
&\qquad
+2\sum_{\nu,j=1}^{4}
g\!\left(
Je_\nu,
\left[
\nabla^L_{e_j}(\nabla^L_{e_\nu}J)
-\nabla^L_{\nabla^L_{e_j}e_\nu}J
\right]e_j
\right)
+\sum_{\alpha,\nu,j=1}^{4}
g\!\left(Je_\alpha,(\nabla^L_{e_\alpha}J)e_j\right)
g\!\left(Je_\nu,(\nabla^L_{e_\nu}J)e_j\right)
\nonumber\\
&\qquad
-\sum_{\nu,j=1}^{4}
g\!\left(
(\nabla^L_{e_\nu}J)e_j,
(\nabla^L_{e_\nu}J)e_j
\right)
+\sum_{\alpha,\nu,j=1}^{4}
g\!\left(Je_\alpha,(\nabla^L_{e_\nu}J)e_j\right)
g\!\left((\nabla^L_{e_\alpha}J)e_j,Je_\nu\right)
\Bigg]\,d\mathrm{vol}_g\nonumber\\
&+
\int_{\partial M}
\Biggl\{
-\left(
\frac1{12}g(X^T,Y^T)+\frac12X_nY_n
\right)h'(0)
-\frac12X(Y_n)
+\frac14\partial_{x_n}(X_nY_n)
\notag\\
&\qquad
-\frac1{12}\partial_{x_n}g(X^T,Y^T)
\Biggr\}\pi\Omega_3\,d\mathrm{vol}_{\partial M}.
\end{align}
where  $s$ is the scalar curvature.
\end{thm}

 \section{The residue for the $J$-twist of the Dirac operator $\widetilde{\nabla}_{X}^{J}\widetilde{\nabla}_{Y}^{J} D_{J}^{-1}$ and $D_{J}^{-3}$ }

In this section, we study the boundary residue functional
associated with
$\widetilde{\nabla}^{J}_X
\widetilde{\nabla}^{J}_Y D_J^{-1}$
and $D_J^{-3}$.
Its interior contribution is given by the
four-dimensional specialization of Theorem~\ref{thm1}.
It remains to compute the boundary contribution.

In addition to the collar assumptions of the preceding
section, assume that
\[
J(\partial_{x_n})=\partial_{x_n}
\]
throughout the collar neighborhood.
\begin{cor}
 Let $(M,g,J)$ be a closed oriented four-dimensional
 almost product Riemannian spin manifold, and let
 $\widetilde{\nabla}^{J}$ be the connection defined
 in \eqref{eq:connection-J}.
 Then, for $X,Y\in\Gamma(TM)$, the following identity holds:
 \begin{align}
&\operatorname{Wres}(\widetilde{\nabla}_{X}^{J}\widetilde{\nabla}_{Y}^{J}D_{J}^{-1}
  \circ D_{J}^{-3})\nonumber\\
=&\frac{4\pi^2}{3}
\int_M
\left[
\operatorname{Ric}(X,Y)-\frac{s}{2}g(X,Y)
\right]\,d\mathrm{vol}_g
\nonumber\\[2pt]
&+2\pi^2
\int_M\sum_{l=1}^{4}
\Big\{
X\!\left[g\!\left(Je_l,(\nabla^L_{e_l}J)Y\right)\right]
-Y\!\left[g\!\left(Je_l,(\nabla^L_{e_l}J)X\right)\right]
-g\!\left(Je_l,(\nabla^L_{e_l}J)[X,Y]\right)
\Big\}\,d\mathrm{vol}_g
\nonumber\\[2pt]
&+\pi^2\int_M g(X,Y)
\Bigg[
-\sum_{i,j=1}^{4}R(Je_i,Je_j,e_j,e_i)
+2\sum_{\nu,j=1}^{4}
g\!\left(
(\nabla^L_{e_j}J)e_\nu,
(\nabla^L_{e_\nu}J)e_j
\right)
\nonumber\\
&\qquad
+2\sum_{\nu,j=1}^{4}
g\!\left(
Je_\nu,
\left[
\nabla^L_{e_j}(\nabla^L_{e_\nu}J)
-\nabla^L_{\nabla^L_{e_j}e_\nu}J
\right]e_j
\right)
+\sum_{\alpha,\nu,j=1}^{4}
g\!\left(Je_\alpha,(\nabla^L_{e_\alpha}J)e_j\right)
g\!\left(Je_\nu,(\nabla^L_{e_\nu}J)e_j\right)
\nonumber\\
&\qquad
-\sum_{\nu,j=1}^{4}
g\!\left(
(\nabla^L_{e_\nu}J)e_j,
(\nabla^L_{e_\nu}J)e_j
\right)
+\sum_{\alpha,\nu,j=1}^{4}
g\!\left(Je_\alpha,(\nabla^L_{e_\nu}J)e_j\right)
g\!\left((\nabla^L_{e_\alpha}J)e_j,Je_\nu\right)
\Bigg]\,d\mathrm{vol}_g,
\end{align}
where $s$ is the scalar curvature.
\end{cor}

Applying Theorem~\ref{thmb0} to this composition, we obtain
\begin{align}\label{c5}
	\widetilde\Phi=&\int_{|\xi'|=1}\int_{-\infty}^{+\infty}\sum_{j,k=0}^{\infty}\sum \frac{(-i)^{|\alpha|+j+k+1}}{\alpha!(j+k+1)!}
	\mathrm{Tr}_{S(TM)}[\partial_{x_{n}}^{j}\partial_{\xi'}^{\alpha}\partial_{\xi_{n}}^{k}\sigma_{r}^{+}
	(\widetilde{\nabla}^J_{X}\widetilde{\nabla}^J_{Y} D_{J}^{-1})(x',0,\xi',\xi_{n})\nonumber\\
	&\times\partial_{x'}^{\alpha}\partial_{\xi_{n}}^{j+1}\partial_{x_{n}}^{k}\sigma_{l}(D_{J}^{-3})(x',0,\xi',\xi_{n})]
	\texttt{d}\xi_{n}\sigma(\xi')\texttt{d}x' ,
\end{align}

From lemma \ref{lem3} and lemma \ref{lem4}, we have
\begin{lem} The following identities hold:
\begin{align}
\sigma_{1}(\widetilde{\nabla}^J_{X}\widetilde{\nabla}^J_{Y} D_{J}^{-1})=&
-\sqrt{-1}\sum_{j,l=1}^nX_jY_l\xi_j\xi_l c[J(\xi)] |\xi|^{-2};\\
\sigma_{0}(\widetilde{\nabla}^J_{X}\widetilde{\nabla}^J_{Y} D_{J}^{-1})=&
\sigma_{2}(\widetilde{\nabla}^J_{X}\widetilde{\nabla}^J_{Y})\sigma_{-2}(D_{J}^{-1})
+\sigma_{1}(\widetilde{\nabla}^J_{X}\widetilde{\nabla}^J_{Y}) \sigma_{-1}(D_{J}^{-1})\nonumber\\
&+\sum_{j=1}^{n}\partial _{\xi_{j}}\big[\sigma_{2}(\widetilde{\nabla}^J_{X}\widetilde{\nabla}^J_{Y})\big]
D_{x_{j}}\big[\sigma_{-1}(D_{J}^{-1})\big].
\end{align}
\end{lem}

Write
\begin{equation}
	\sigma(D_J^3)=p_3+p_2+p_1+p_0,
	\qquad
	\sigma(D_J^{-3})
	=\sum_{j=3}^{\infty}q_{-j}.
\end{equation}
By the composition formula of pseudodifferential operators, we have
\begin{align}
	1
	&\sim \sigma(D_J^3\circ D_J^{-3})
	\nonumber\\
	&\sim
	\sum_{\alpha}\frac{1}{\alpha!}
	\partial_{\xi}^{\alpha}\bigl[\sigma(D_J^3)\bigr]
	D_x^{\alpha}\bigl[\sigma(D_J^{-3})\bigr]
	\nonumber\\
	&\sim
	(p_3+p_2+p_1+p_0)
	(q_{-3}+q_{-4}+q_{-5}+\cdots)
	\nonumber\\
	&\quad+
	\sum_{j=1}^{n}
	\bigl(
	\partial_{\xi_j}p_3+\partial_{\xi_j}p_2
	+\partial_{\xi_j}p_1+\partial_{\xi_j}p_0
	\bigr)
	\bigl(
	D_{x_j}q_{-3}+D_{x_j}q_{-4}
	+D_{x_j}q_{-5}+\cdots
	\bigr)
	+\cdots
	\nonumber\\
	&\sim
	p_3q_{-3}
	+
	\bigg(
	p_3q_{-4}+p_2q_{-3}
	+\sum_{j=1}^{n}
	(\partial_{\xi_j}p_3)D_{x_j}q_{-3}
	\bigg)
	+\cdots.
\end{align}
so
\begin{equation}
q_{-3}=p_3^{-1};~q_{-4}=-p_3^{-1}[p_2p_3^{-1}+\sum_j\partial_{\xi_j}p_3D_{x_j}(p_3^{-1})].
\end{equation}
Hence by Lemma 4.1 in \cite{WW}, (2.4) and (4.7), we have
\begin{lem} The following identities hold:
\begin{align}
\sigma_{-3}(D_{J}^{-3})&=\sqrt{-1} c[J(\xi)] |\xi|^{-4};\\
\sigma_{-4}({D}_{J}^{-3})&=\frac{c[J(\xi)]\sigma_{2}({D}_{J}^{3})c[J(\xi)]}{|\xi|^{8}}+\frac{c[J(\xi)]}{|\xi|^{10}}\sum_ {j=1}^{n} \Big(c[J(dx_j)]|\xi|^{2}+2\xi^{j}c[J(\xi)]\Big)\\
&\times \Big[\partial_{x_j}(c[J(\xi)])|\xi|^2-2c[J(\xi)]\partial_{x_j}(|\xi|^2)\Big].\nonumber
\end{align}
where,
\begin{align}
\sigma_2({D}_{J}^{3})&=\sum^{n}_{i,j,l=1}c[J(dx_{l})]\partial_{l}(g^{ij})\xi_{i}\xi_{j}+c[J(\xi)](2\sigma^k-\Gamma^k)\xi_{k}\\
&-\sum^{n}_{\alpha=1}c[J(\xi)]c[J(e_{\alpha})]c[(\nabla^{L}_{e_{\alpha}}J)(\xi^{*})]-\frac{1}{4}|\xi|^2\sum^{n}_{s,t,l=1}\omega_{s,t}(e_{l})c[J(e_{l})]c(e_{s})c(e_{t}).\nonumber
\end{align}
\end{lem}
Now we  need to compute $\int_{\partial M} \widetilde{\Phi}$. When $n=4$, then ${\rm Tr}_{S(TM)}[\operatorname{Id}]=4$, the sum is taken over $
r+l-k-j-|\alpha|=-3,~~r\leq 1,~~l\leq-3,$ then we have the following five cases:

 {\bf case a)~I)}~$r=1,~l=-3,~k=j=0,~|\alpha|=1$.

By (4.2), we get
\begin{equation}
\label{b241}
\widetilde{\Phi}_1=-\int_{|\xi'|=1}\int^{+\infty}_{-\infty}\sum_{|\alpha|=1}
\mathrm{Tr}[\partial^\alpha_{\xi'}\pi^+_{\xi_n}\sigma_{1}(\widetilde{\nabla}^J_{X}\widetilde{\nabla}^J_{Y} D^{-1}_{J})\times
 \partial^\alpha_{x'}\partial_{\xi_n}\sigma_{-3}(D^{-3}_{J})](x_0)d\xi_n\sigma(\xi')dx'.
\end{equation}
Write $J(e_p)=\sum_{h=1}^{n}a_h^p e_h$
in the adapted orthonormal frame.
By the collar assumption and the orthogonality of $J$,
we have
\[
a_h^n=a_n^h=0\quad(h<n),
\qquad a_n^n=1.
\]
For $\beta<n$, Lemma 2.2 in \cite{Wa3} yields
\begin{align}
\partial_{x_\beta}\sigma_{-3}(D^{-3}_{J})(x_0)=& \partial_{x_\beta}(\sqrt{-1} c[J(\xi)] |\xi|^{-4})(x_0)\nonumber\\
=&\frac{\sqrt{-1}}{(1+\xi_{n}^2)^2}\sum^{n}_{h=1}\sum^{n-1}_{p=1}\xi_{p}\partial_{x_\beta}(a^{p}_{h})c(e_{h}).
\end{align}
then
\begin{align}
	\partial_{\xi_n} \partial_{x_\beta}\sigma_{-3}(D^{-3}_{J})(x_0)|_{|\xi'|=1}=-\frac{4\sqrt{-1}{\xi_n}}{(1+\xi_{n}^2)^3}\sum^{n}_{h=1}\sum^{n-1}_{p=1}\xi_{p}\partial_{x_\beta}(a^{p}_{h})c(e_{h}).
\end{align}
and
\begin{align}
	&\partial_{\xi_\beta}
	\sigma_1\bigl(
	\widetilde{\nabla}^J_{X}\widetilde{\nabla}^J_{Y} D_J^{-1}
	\bigr)(x_0)\big|_{|\xi'|=1}
	\nonumber\\
	&=-\frac{i}{1+\xi_n^2}
	\sum_{l=1}^{n}
	(X_\beta Y_l+X_lY_\beta)
	\xi_lc[J(\xi)]
	-\frac{i}{1+\xi_n^2}
	\sum_{j,l=1}^{n}X_jY_l\xi_j\xi_l
	\sum_{h=1}^{n}a_h^\beta c(e_h)
	\nonumber\\
	&\quad
	+\frac{2i\xi_\beta}{(1+\xi_n^2)^2}
	\sum_{j,l=1}^{n}
	X_jY_l\xi_j\xi_lc[J(\xi)].
\end{align}
then
\begin{align}
	&\pi^+_{\xi_n} \partial_{\xi_\beta}\sigma_{1}(\widetilde{\nabla}^J_{X}\widetilde{\nabla}^J_{Y} D^{-1}_{J})(x_0)|_{|\xi'|=1}\nonumber\\
	&=-\frac{c[J(\xi')]+ic(e_n)}{2(\xi_n-i)}
	\Biggl[
	\sum_{l=1}^{n-1}
	(X_\beta Y_l+X_lY_\beta)\xi_l
	-\xi_\beta
	\left(
	\sum_{j,l=1}^{n-1}X_jY_l\xi_j\xi_l+X_nY_n
	\right)
	+i(X_\beta Y_n+X_nY_\beta)
	\Biggr]
	\nonumber\\
	&\quad
	-\frac{
		\displaystyle\sum_{h=1}^{n}a_h^\beta c(e_h)
		+i\xi_\beta c(e_n)}
	{2(\xi_n-i)}
	\Biggl[
	\sum_{j,l=1}^{n-1}X_jY_l\xi_j\xi_l-X_nY_n
	+i\sum_{j=1}^{n-1}
	(X_jY_n+X_nY_j)\xi_j
	\Biggr]
	\nonumber\\
	&\quad
	-\frac{i\xi_\beta\bigl(c[J(\xi')]+ic(e_n)\bigr)}
	{2(\xi_n-i)^2}
	\Biggl[
	\sum_{j,l=1}^{n-1}X_jY_l\xi_j\xi_l-X_nY_n
	+i\sum_{j=1}^{n-1}
	(X_jY_n+X_nY_j)\xi_j
	\Biggr].
\end{align}

\begin{align}
	&\mathrm{Tr}\Bigl[
	\pi_{\xi_n}^{+}\partial_{\xi_\beta}
	\sigma_1\bigl(
	\widetilde{\nabla}^J_X
	\widetilde{\nabla}^J_Y D_J^{-1}
	\bigr)
	\times
	\partial_{\xi_n}\partial_{x_\beta}
	\sigma_{-3}(D_J^{-3})
	\Bigr](x_0)\big|_{|\xi'|=1}
	\nonumber\\
	&\sim
	\frac{8\xi_n}{(\xi_n-i)^4(\xi_n+i)^3}
	\sum_{h=1}^{n}\sum_{p,q=1}^{n-1}
	(X_qY_n+X_nY_q)\xi_p\xi_q
	a_h^\beta\partial_{x_\beta}(a_h^p).
\end{align}

\begin{align}
	\widetilde{\Phi}_1=&-\int_{|\xi'|=1}\int^{+\infty}_{-\infty}\sum_{|\alpha|=1}
	\mathrm{Tr} [\pi^+_{\xi_n} \partial_{\xi_\beta}\sigma_{1}(\widetilde{\nabla}^J_{X}\widetilde{\nabla}^J_{Y} D^{-1}_{J}) \times \partial_{\xi_n} \partial_{x_\beta}\sigma_{-3}(D^{-3}_{J})](x_0)d\xi_n\sigma(\xi')dx'\nonumber\\
	=&-\frac{1}{6}
	\sum_{\beta,p=1}^{n-1}\sum_{h=1}^{n}
	(X_pY_n+X_nY_p)
	a_h^\beta\partial_{x_\beta}(a_h^p)\,\pi\Omega_3 dx'.
\end{align}

  {\bf case a)~II)}~$r=1,~l=-3,~k=|\alpha|=0,~j=1$.

  By (4.2), we get
\begin{equation}
\label{b261}
\widetilde{\Phi}_2=-\frac{1}{2}\int_{|\xi'|=1}\int^{+\infty}_{-\infty}
\mathrm{Tr} [\partial_{x_n}\pi^+_{\xi_n}\sigma_{1}(\widetilde{\nabla}^J_{X}\widetilde{\nabla}^J_{Y} D_{J}^{-1})\times
\partial_{\xi_n}^2\sigma_{-3}(D_{J}^{-3})](x_0)d\xi_n\sigma(\xi')dx'.
\end{equation}
By Lemma 4.3, we have
\begin{align}
	\partial_{\xi_n}^2\sigma_{-3}(D_{J}^{-3})(x_0)|_{|\xi'|=1}=&\partial_{\xi_n}^2(\sqrt{-1} c[J(\xi)] |\xi|^{-4})(x_0)|_{|\xi'|=1}\nonumber\\
	=&4\sqrt{-1}\frac{(5\xi_n^2-1)}{(1+\xi_n^2)^4}\sum_{\beta=1}^{n} \sum_{j=1}^{n-1}\xi_{j}a_{\beta}^{j}c(e_{\beta})+12\sqrt{-1}
	\frac{\xi_n^3-\xi_n}
	{(1+\xi_n^2)^4}c(dx_n).\nonumber
\end{align}
and
\begin{align}
	&\partial_{x_n}\sigma_{1}(\widetilde{\nabla}^J_{X}\widetilde{\nabla}^J_{Y} D_{J}^{-1}) (x_0) |_{|\xi'|=1}\nonumber\\
	&=\partial_{x_n}(-\sqrt{-1}\sum_{k,l=1}^n X_k Y_l \xi_k \xi_l c[J(\xi)] |\xi|^{-2})\nonumber\\
	&=-\sqrt{-1}\sum_{k,l,p=1}^{n} X_k Y_l \xi_k \xi_l \xi_p \bigg[ \sum_{h=1}^{n} \frac{\partial_{x_n}(a_h^p)c(e_h)}{1+\xi_n^2}+\sum_{h=1}^{n-1}\frac{a_h^p \left.\partial_{x_n}c(dx_h)\right|_{x_0}}{1+\xi_n^2}-\sum_{h=1}^{n} \frac{(a_h^p)c(e_h)}{(1+\xi_n^2)^2}h'(0)\bigg]\nonumber\\
	&\quad-\sqrt{-1}\sum_{k,l,p=1}^{n}
	\partial_{x_n}(X_kY_l)\xi_k\xi_l\xi_p
	\sum_{h=1}^{n}
	\frac{a_h^pc(e_h)}{1+\xi_n^2}.\nonumber
\end{align}

We note that $\beta<n,~\int_{|\xi'|=1}\xi_{\beta_{1}}\xi_{\beta_{2}}\cdots\xi_{\beta_{2d+1}}\sigma(\xi')=0$,
so we omit some items that have no contribution for computing {\bf case a)~II)}.
Then there is the following formula

\begin{align}
	&\left.
	\mathrm{Tr}\!\left[
	\partial_{x_n}\pi_{\xi_n}^{+}
	\sigma_1\!\left(
	\widetilde{\nabla}^J_X
	\widetilde{\nabla}^J_Y
	D_J^{-1}
	\right)
	\partial_{\xi_n}^{2}\sigma_{-3}(D_J^{-3})
	\right](x_0)
	\right|_{|\xi'|=1}
	\nonumber\\
	&\sim
	-\frac{4i(3\xi_n^2-2i\xi_n-1)}
	{(\xi_n-i)^5(\xi_n+i)^4}
	\sum_{k,l=1}^{n-1}
	X_kY_l\xi_k\xi_l\,h'(0)+
	\frac{4(3\xi_n^2-2i\xi_n-1)(2\xi_n-i)}
	{(\xi_n-i)^5(\xi_n+i)^4}
	X_nY_n\,h'(0)
	\nonumber\\
	&\quad-
	\frac{8(3\xi_n^2-2i\xi_n-1)}
	{(1+\xi_n^2)^4}
	\sum_{k,l=1}^{n-1}
	\partial_{x_n}(X_kY_l)\xi_k\xi_l+
	\frac{8(3\xi_n^2-2i\xi_n-1)}
	{(1+\xi_n^2)^4}
	\partial_{x_n}(X_nY_n).
\end{align}

Therefore, we get
\begin{align}
	\widetilde{\Phi}_2
	&=\Biggl[
	\left(
	\frac{5}{32}\sum_{j=1}^{n-1}X_jY_j
	+\frac{1}{32}X_nY_n
	\right)h'(0)
	\nonumber\\
	&\qquad
	-\frac16\sum_{j=1}^{n-1}
	\partial_{x_n}(X_jY_j)
	+\frac12\partial_{x_n}(X_nY_n)
	\Biggr]\pi\Omega_3\,dx'.
\end{align}

 {\bf case a)~III)}~$r=1,~l=-3,~j=|\alpha|=0,~k=1$.

 By (4.2), we get
\begin{align}\label{36i}
\widetilde{\Phi}_3&=-\frac{1}{2}\int_{|\xi'|=1}\int^{+\infty}_{-\infty}
\mathrm{Tr} [\partial_{\xi_n}\pi^+_{\xi_n}\sigma_{1}({\nabla}^J_{X}{\nabla}^J_{Y} D_{J}^{-1})\times
\partial_{\xi_n}\partial_{x_n}\sigma_{-3}(D_{J}^{-3})](x_0)d\xi_n\sigma(\xi')dx'.
\end{align}
\noindent By Lemma 4.3, we have
\begin{align}\label{37i}
\partial_{\xi_n}\partial_{x_n}\sigma_{-3}(D_{J}^{-3})(x_0)|_{|\xi'|=1}
=&\partial_{\xi_n}\partial_{x_n}\left(\frac{ic[J(\xi)]}{|\xi|^4}\right)(x_0)|_{|\xi'|=1}\\
=&-\frac{4 i \xi _n}{\left(\xi _n^2+1\right)^3}\sum^{n}_{h=1}\sum^{n-1}_{p=1}\xi_{p}\partial_{x_n}(a^{p}_{h})c(e_{h})
+\frac{i(1-3 \xi _n^2)}{\left(\xi _n^2+1\right)^3}\sum^{n}_{h=1}\partial_{x_n}(a^{n}_{h})c(e_{h})\nonumber\\
&-\frac{4 i \xi _n}{\left(\xi _n^2+1\right)^3}\sum^{n-1}_{h,p=1}\xi_{p}a^{p}_{h}\left.\partial_{x_n}c(dx_h)\right|_{x_0}-\frac{2i(1-5 \xi _n^2)}{\left(\xi _n^2+1\right)^4}h'(0)\sum^{n}_{h=1}a^{n}_{h}c(e_{h})\nonumber\\
&+\frac{12 i \xi _n}{\left(\xi _n^2+1\right)^4}h'(0)\sum^{n}_{h=1}\sum^{n-1}_{p=1}\xi_{p}a^{p}_{h}c(e_{h}),\nonumber
\end{align}

and we have
\begin{align}
\sigma_{1}(\widetilde{\nabla}^J_{X}\widetilde{\nabla}^J_{Y} D_{J}^{-1})(x_0)|_{|\xi'|=1}
=&-i\sum_{j,l=1}^nX_jY_l\xi_j\xi_l c[J(\xi)] |\xi|^{-2}\nonumber\\
=&-\frac{i}{1+\xi_{n}^{2}}\sum_{j,l,k,u=1}^{n} X_j Y_l \xi_j \xi_l \xi_k a_{u}^{k} c(e_u)\nonumber\\
=&-\frac{i}{1+\xi_{n}^{2}}\sum_{u=1}^{n}\sum_{j,l,k=1}^{n-1} X_j Y_l \xi_j \xi_l \xi_k a_{u}^{k} c(e_u) -\frac{i\xi_{n}}{1+\xi_{n}^{2}}\sum_{u=1}^{n}\sum_{j,l=1}^{n-1} X_j Y_l \xi_j \xi_l a_{u}^{n} c(e_u)\nonumber\\
&-\frac{i\xi_{n}}{1+\xi_{n}^{2}}\sum_{u=1}^{n}\sum_{j,k=1}^{n-1} X_j Y_n \xi_j \xi_k a_{u}^{k} c(e_u) -\frac{i\xi_{n}}{1+\xi_{n}^{2}}\sum_{u=1}^{n}\sum_{l,k=1}^{n-1} X_n Y_l \xi_l \xi_k a_{u}^{k} c(e_u)\nonumber\\
&-\frac{i\xi_{n}^{2}}{1+\xi_{n}^{2}}\sum_{u=1}^{n}\sum_{j=1}^{n-1} X_j Y_n \xi_j a_{u}^{n} c(e_u) -\frac{i\xi_{n}^2}{1+\xi_{n}^{2}}\sum_{u=1}^{n}\sum_{l=1}^{n-1} X_n Y_l \xi_l a_{u}^{n} c(e_u)\nonumber\\
&-\frac{i\xi_{n}^2}{1+\xi_{n}^{2}}\sum_{u=1}^{n}\sum_{k=1}^{n-1} X_n Y_n\xi_k a_{u}^{k} c(e_u) -\frac{i\xi_{n}^3}{1+\xi_{n}^{2}}\sum_{u=1}^{n} X_n Y_n a_{u}^{n} c(e_u),
\end{align}
then
\begin{align}
	\partial_{\xi_n}\pi^+_{\xi_n}\sigma_{1}(\widetilde{\nabla}^J_{X}\widetilde{\nabla}^J_{Y} D_{J}^{-1})(x_0)|_{|\xi'|=1}
	=\;\;&\frac{1}{2(\xi_n-i)^2}\sum_{u=1}^{n}\sum_{j,l,k=1}^{n-1} X_j Y_l \xi_j \xi_l \xi_k a_{u}^{k} c(e_u) +\frac{i}{2(\xi_n-i)^2}\sum_{u=1}^{n}\sum_{j,l=1}^{n-1} X_j Y_l \xi_j \xi_l a_{u}^{n} c(e_u)\nonumber\\
	&+\frac{i}{2(\xi_n-i)^2}\sum_{u=1}^{n}\sum_{j,k=1}^{n-1} X_j Y_n \xi_j \xi_k a_{u}^{k} c(e_u) +\frac{i}{2(\xi_n-i)^2}\sum_{u=1}^{n}\sum_{l,k=1}^{n-1} X_n Y_l \xi_l \xi_k a_{u}^{k} c(e_u)\nonumber\\
	&-\frac{1}{2(\xi_n-i)^2}\sum_{u=1}^{n}\sum_{j=1}^{n-1} X_j Y_n \xi_j a_{u}^{n} c(e_u) -\frac{1}{2(\xi_n-i)^2}\sum_{u=1}^{n}\sum_{l=1}^{n-1} X_n Y_l \xi_l a_{u}^{n} c(e_u)\nonumber\\
	&-\frac{1}{2(\xi_n-i)^2}\sum_{u=1}^{n}\sum_{k=1}^{n-1} X_n Y_n\xi_k a_{u}^{k} c(e_u) -\frac{i}{2(\xi_n-i)^2}\sum_{u=1}^{n} X_n Y_n a_{u}^{n} c(e_u).
\end{align}

We note that $\beta<n,~\int_{|\xi'|=1}\xi_{\beta_{1}}\xi_{\beta_{2}}\cdots\xi_{\beta_{2d+1}}\sigma(\xi')=0$,
so we omit some items that have no contribution for computing {\bf case a)~III)}, then
\begin{align}
	&\left.
	\mathrm{Tr}\!\left[
	\partial_{\xi_n}\pi_{\xi_n}^{+}
	\sigma_1\!\left(
	\widetilde{\nabla}^J_X
	\widetilde{\nabla}^J_Y
	D_J^{-1}
	\right)
	\partial_{\xi_n}\partial_{x_n}
	\sigma_{-3}(D_J^{-3})
	\right](x_0)
	\right|_{|\xi'|=1}
	\nonumber\\
	&\sim
	\frac{4i(\xi_n^2-4i\xi_n-1)h'(0)}
	{(\xi_n-i)^5(\xi_n+i)^4}
	\left[
	\sum_{j,l=1}^{n-1}
	X_jY_l\xi_j\xi_l-X_nY_n
	\right].
	\nonumber
\end{align}

Therefore, we get
\begin{align}\label{41i}
\widetilde{\Phi}_3=\left(-\frac{25 }{96}\sum_{j=1}^{n-1}X_jY_j+\frac{25}{32}X_nY_n\right)h'(0)\pi\Omega_3dx'.
\end{align}

 {\bf case b)}~$r=0,~l=-3,~k=j=|\alpha|=0$.

By (4.2), we get
\begin{align}\label{42i}
\widetilde{\Phi}_4&=-i\int_{|\xi'|=1}\int^{+\infty}_{-\infty}
\mathrm{Tr}[\pi^+_{\xi_n}\sigma_{0}(\widetilde{\nabla}^J_X
\widetilde{\nabla}^J_Y D_{J}^{-1})\times
\partial_{\xi_n}\sigma_{-3}(D_{J}^{-3})](x_0)d\xi_n\sigma(\xi')dx'.
\end{align}
By Lemma 4.3, we obtain
\begin{align}\label{43i}
\partial_{\xi_n}\sigma_{-3}(D_{J}^{-3})(x_0)|_{|\xi'|=1}
&=\partial_{\xi_n}\left(\frac{ic[J(\xi)]}{|\xi|^{4}}\right)(x_0)|_{|\xi'|=1}\\
&=-\frac{4 i \xi_n}{\left(\xi_n^2+1\right)^3}\sum^{n}_{\beta=1}\sum^{n-1}_{k=1}\xi_{k}a_{\beta}^{k}c(e_{\beta})+\frac{i(1-3 \xi _n^2)}{\left(\xi _n^2+1\right){}^3}\sum^{n}_{\beta=1}a_{\beta}^nc(e_{\beta}),
\end{align}

By Lemma 4.2, we have
\begin{align}\label{b41}
\sigma_{0}(\widetilde{\nabla}^J_{X}\widetilde{\nabla}^J_{Y} D_{J}^{-1})=&
\sigma_{2}(\widetilde{\nabla}^J_{X}\widetilde{\nabla}^J_{Y})\sigma_{-2}(D_{J}^{-1})
+\sigma_{1}(\widetilde{\nabla}^J_{X}\widetilde{\nabla}^J_{Y})\sigma_{-1}(D_{J}^{-1})\nonumber\\
&+\sum_{j=1}^{n}\partial _{\xi_{j}}\big[\sigma_{2}(\widetilde{\nabla}^J_{X}\widetilde{\nabla}^J_{Y})\big]
D_{x_{j}}\big[\sigma_{-1}(D_{J}^{-1})\big]\nonumber\\
&:=\mathcal{A}+\mathcal{B}+\mathcal{C}.
\end{align}

(1) For the first term in  \eqref{b41}, we have
\begin{align}
\mathcal{A}(x_0):=&\sigma_{2}(\widetilde{\nabla}^J_{X}\widetilde{\nabla}^J_{Y})\sigma_{-2}(D_{J}^{-1})(x_0)|_{|\xi'|=1}
=-\sum_{j,l=1}^{n}X_jY_l\xi_j\xi_l\sigma_{-2}(D_{J}^{-1})(x_0)|_{|\xi'|=1}
\end{align}
where
\begin{align}
	\sigma_{-2}({D}_{J}^{-1})(x_0)|_{|\xi'|=1}
	&=\frac{c[J(\xi)]\sigma_{0}({D}_{J})(x_0)c[J(\xi)]}{|\xi|^4}-\frac{c[J(\xi)]}{|\xi|^6}h'(0)c(e_n)c[J(\xi)]\nonumber\\
	&+\frac{c[J(\xi)]}{|\xi|^4}\bigg(\sum_{j,p,h=1}^{n}\xi_p\partial_{x_j}(a_{h}^{p})c[J(e_j)]c(e_h)+\sum_{p=1}^{n}\sum_{h=1}^{n-1}\xi_pa_{h}^{p}c(e_n)\left.\partial_{x_n}c(dx_h)\right|_{x_0}\bigg),
\end{align}
and
\begin{align}
	\sigma_{0}({{D}_{J}})(x_0)=-\frac{1}{4}h'(0)\sum_{\mu=1}^{n}\sum_{\nu=1}^{n-1}a_{\nu}^{\mu}c(e_{\mu})c(e_{n})c(e_{\nu}).
\end{align}

We decompose  $\mathcal{A}(x_0)$ as $\mathcal{A}(x_0)=\mathcal{A}_1 (x_0) -h'(0) \mathcal{A}_2 (x_0) + \mathcal{A}_3 (x_0) + \mathcal{A}_4 (x_0)$, where
\begin{align}
	\mathcal{A}_1 (x_0):=&-\frac{1}{(1+\xi_n^2)^2}\sum_{j,l=1}^{n}X_jY_l\xi_j\xi_lc[J(\xi)]\sigma_{0}({D}_{J})(x_0)c[J(\xi)]\nonumber\\
	=&-\sum_{j,l,q,\alpha,p,\gamma=1}^{n}\frac{\xi_j\xi_l\xi_q\xi_{\alpha}}{(1+\xi_n^2)^2}X_jY_l a_p^q a_{\gamma}^{\alpha} c(e_p)\sigma_{0}({D}_{J})(x_0)c(e_\gamma)\\
	\mathcal{A}_2 (x_0):=&-\frac{1}{(1+\xi_n^2)^3}\sum_{j,l=1}^{n}X_jY_l\xi_j\xi_lc[J(\xi)]c(e_n)c[J(\xi)]\nonumber\\
	=&-\sum_{j,l,q,\alpha,p,\gamma=1}^{n}\frac{\xi_j\xi_l\xi_q\xi_{\alpha}}{(1+\xi_n^2)^3}X_jY_l a_p^q a_{\gamma}^{\alpha} c(e_p)c(e_n)c(e_\gamma)\\
	\mathcal{A}_3 (x_0):=&-\frac{1}{(1+\xi_n^2)^2}\sum_{j,l,\gamma,p,h=1}^{n}X_jY_l\xi_j\xi_lc[J(\xi)]\xi_p\partial_{x_{\gamma}}(a_{h}^{p})c[J(e_{\gamma})]c(e_h)\nonumber\\
	=&-\sum_{j,l,\gamma,p,h,q,\theta,y=1}^{n}\frac{\xi_j\xi_l\xi_q\xi_p}{(1+\xi_n^2)^2}X_jY_l a_{\theta}^{q} a_{y}^{\gamma}\partial_{x_{\gamma}}(a_{h}^{p})c(e_{\theta})c(e_y)c(e_h)\\
	\mathcal{A}_4 (x_0):=&-\frac{1}{(1+\xi_n^2)^2}\sum_{j,l,p=1}^{n} \sum_{h=1}^{n-1} X_jY_l\xi_j\xi_l c[J(\xi)] \xi_p a_{h}^{p}c(e_n)\left.\partial_{x_n}c(dx_h)\right|_{x_0}\nonumber\\
	=&-\frac{1}{2} h'(0) \sum_{j,l,q,\theta=1}^{n} \sum_{p,h=1}^{n-1}\frac{\xi_j\xi_l\xi_q\xi_p}{(1+\xi_n^2)^2} X_jY_l a_{\theta}^{q} a_{h}^{p} c(e_{\theta}) c(e_n)c(e_h).
\end{align}
Using the identities
\begin{align}
	&\pi_{\xi_n}^{+}\frac{1}{(1+\xi_n^2)^2}=-\frac{i\xi_n+2}{4(\xi_n-i)^2},\;\;\;\;\;\;\;\pi_{\xi_n}^{+}\frac{\xi_n}{(1+\xi_n^2)^2}=-\frac{i}{4(\xi_n-i)^2},\nonumber\\
	&\pi_{\xi_n}^{+}\frac{\xi_n^2}{(1+\xi_n^2)^2}=-\frac{i\xi_n}{4(\xi_n-i)^2},\;\;\;\;\;\;\;\pi_{\xi_n}^{+}\frac{\xi_n^3}{(1+\xi_n^2)^2}=\frac{2\xi_n-i}{4(\xi_n-i)^2},\nonumber\\
	&\pi_{\xi_n}^{+}\frac{\xi_n^4}{(1+\xi_n^2)^2}=\frac{3i\xi_n+2}{4(\xi_n-i)^2},\nonumber
\end{align}
we obtain
\begin{align}
	\pi_{\xi_n}^{+}\mathcal{A}_1 (x_0)=&-\frac{3i\xi_n+2}{4(\xi_n-i)^2}\sum_{p,\gamma=1}^{n}X_n Y_n a_p^n a_{\gamma}^{n} c(e_p)\sigma_{0}({D}_{J})(x_0)c(e_\gamma)\nonumber\\
	&-\frac{2\xi_n-i}{4(\xi_n-i)^2}	\sum_{p,\gamma=1}^{n}\sum_{j=1}^{n-1}\xi_j X_j Y_n a_p^n a_{\gamma}^{n} c(e_p)\sigma_{0}({D}_{J})(x_0)c(e_\gamma)\nonumber\\
	&-\frac{2\xi_n-i}{4(\xi_n-i)^2}	\sum_{p,\gamma=1}^{n}\sum_{l=1}^{n-1}\xi_l X_n Y_l a_p^n a_{\gamma}^{n} c(e_p)\sigma_{0}({D}_{J})(x_0)c(e_\gamma)\nonumber\\
	&-\frac{2\xi_n-i}{4(\xi_n-i)^2}	\sum_{p,\gamma=1}^{n}\sum_{q=1}^{n-1}\xi_q X_n Y_n a_p^q a_{\gamma}^{n}c(e_p)\sigma_{0}({D}_{J})(x_0)c(e_\gamma)\nonumber\\
	&-\frac{2\xi_n-i}{4(\xi_n-i)^2}	\sum_{p,\gamma=1}^{n}\sum_{\alpha=1}^{n-1}\xi_\alpha X_n Y_n a_p^n a_{\gamma}^{\alpha} c(e_p)\sigma_{0}({D}_{J})(x_0)c(e_\gamma)\nonumber\\
    &+\frac{i\xi_n}{4(\xi_n-i)^2} \sum_{p,\gamma=1}^{n}\sum_{j,l=1}^{n-1}\xi_j \xi_l X_j Y_l a_p^n a_{\gamma}^{n} c(e_p)\sigma_{0}({D}_{J})(x_0)c(e_\gamma)\nonumber\\
    &+\frac{i\xi_n}{4(\xi_n-i)^2} \sum_{p,\gamma=1}^{n}\sum_{j,q=1}^{n-1}\xi_j \xi_q X_j Y_n a_p^q a_{\gamma}^{n} c(e_p)\sigma_{0}({D}_{J})(x_0)c(e_\gamma)\nonumber\\
    &+\frac{i\xi_n}{4(\xi_n-i)^2} \sum_{p,\gamma=1}^{n}\sum_{j,\alpha=1}^{n-1}\xi_j \xi_\alpha X_j Y_n a_p^n a_{\gamma}^{\alpha} c(e_p)\sigma_{0}({D}_{J})(x_0)c(e_\gamma)\nonumber\\
    &+\frac{i\xi_n}{4(\xi_n-i)^2} \sum_{p,\gamma=1}^{n}\sum_{l,q=1}^{n-1} \xi_l \xi_q X_n Y_l a_p^q a_{\gamma}^{n} c(e_p)\sigma_{0}({D}_{J})(x_0)c(e_\gamma)\nonumber\\
    &+\frac{i\xi_n}{4(\xi_n-i)^2} \sum_{p,\gamma=1}^{n}\sum_{l,\alpha=1}^{n-1} \xi_l \xi_\alpha X_n Y_l a_p^n a_{\gamma}^{\alpha} c(e_p)\sigma_{0}({D}_{J})(x_0)c(e_\gamma)\nonumber\\
    &+\frac{i\xi_n}{4(\xi_n-i)^2} \sum_{p,\gamma=1}^{n}\sum_{q,\alpha=1}^{n-1} \xi_q \xi_\alpha X_n Y_n a_p^q a_{\gamma}^{\alpha} c(e_p)\sigma_{0}({D}_{J})(x_0)c(e_\gamma)\nonumber\\
    &+\frac{i}{4(\xi_n-i)^2} \sum_{p,\gamma=1}^{n} \sum_{l,q,\alpha=1}^{n-1} \xi_l \xi_q \xi_\alpha X_n Y_l a_p^q a_{\gamma}^{\alpha} c(e_p)\sigma_{0}({D}_{J})(x_0)c(e_\gamma)\nonumber\\
    &+\frac{i}{4(\xi_n-i)^2} \sum_{p,\gamma=1}^{n} \sum_{j,q,\alpha=1}^{n-1} \xi_j \xi_q \xi_\alpha X_j Y_n a_p^q a_{\gamma}^{\alpha} c(e_p)\sigma_{0}({D}_{J})(x_0)c(e_\gamma)\nonumber\\
    &+\frac{i}{4(\xi_n-i)^2} \sum_{p,\gamma=1}^{n}\sum_{j,l,\alpha=1}^{n-1}\xi_j \xi_l  \xi_\alpha X_j Y_l a_p^n a_{\gamma}^{\alpha} c(e_p)\sigma_{0}({D}_{J})(x_0)c(e_\gamma)\nonumber\\
    &+\frac{i}{4(\xi_n-i)^2} \sum_{p,\gamma=1}^{n} \sum_{j,l,q=1}^{n-1} \xi_j \xi_l \xi_q  X_j Y_l a_p^q a_{\gamma}^{n} c(e_p)\sigma_{0}({D}_{J})(x_0)c(e_\gamma)\nonumber\\
    &+\frac{i\xi_n+2}{4(\xi_n-i)^2} \sum_{p,\gamma=1}^{n} \sum_{j,l,q,\alpha=1}^{n-1} \xi_j \xi_l \xi_q \xi_{\alpha}  X_j Y_l a_p^q a_{\gamma}^{\alpha} c(e_p)\sigma_{0}({D}_{J})(x_0)c(e_\gamma),
\end{align}
Consequently,
\begin{align}
	&\bigg(\pi^+_{\xi_n}\mathcal{A}_1\times\partial_{\xi_n}\sigma_{-3}(D_{J}^{-3})\bigg)(x_0)|_{|\xi'|=1}\nonumber\\	
	\sim &\mathbb{C}_1 h'(0) \sum_{v=1}^{n-1}\sum_{p,\gamma,\beta,\mu=1}^{n} X_n Y_n a_p^n a_{\gamma}^{n} a_v^{\mu} a_{\beta}^{n} c(e_p) c(e_{\mu}) c(e_n) c(e_v) c(e_\gamma) c(e_{\beta})\nonumber\\
	&+\mathbb{C}_2 h'(0) \sum_{j,k,v=1}^{n-1}\sum_{p,\gamma,\beta,\mu=1}^{n} X_j Y_n a_p^n a_{\gamma}^{n} a_v^{\mu} a_{\beta}^{k} \xi_j \xi_k c(e_p) c(e_{\mu}) c(e_n) c(e_v) c(e_\gamma) c(e_{\beta})\nonumber\\
	&+\mathbb{C}_2 h'(0) \sum_{l,k,v=1}^{n-1}\sum_{p,\gamma,\beta,\mu=1}^{n} X_n Y_l a_p^n a_{\gamma}^{n} a_v^{\mu} a_{\beta}^{k} \xi_l \xi_k c(e_p) c(e_{\mu}) c(e_n) c(e_v) c(e_\gamma) c(e_{\beta})\nonumber\\
	&+\mathbb{C}_2 h'(0) \sum_{k,q,v=1}^{n-1}\sum_{p,\gamma,\beta,\mu=1}^{n} X_n Y_n a_p^q a_{\gamma}^{n} a_v^{\mu} a_{\beta}^{k} \xi_q \xi_k c(e_p) c(e_{\mu}) c(e_n) c(e_v) c(e_\gamma) c(e_{\beta})\nonumber\\
	&+\mathbb{C}_2 h'(0) \sum_{k,\alpha,v=1}^{n-1}\sum_{p,\gamma,\beta,\mu=1}^{n} X_n Y_n a_p^n a_{\gamma}^{\alpha} a_v^{\mu} a_{\beta}^{k} \xi_{\alpha} \xi_k c(e_p) c(e_{\mu}) c(e_n) c(e_v) c(e_\gamma) c(e_{\beta})\nonumber\\
	&+\mathbb{C}_3 h'(0) \sum_{j,l,v=1}^{n-1}\sum_{p,\gamma,\beta,\mu=1}^{n} X_j Y_l a_p^n a_{\gamma}^{n} a_v^{\mu} a_{\beta}^{n} \xi_j \xi_l c(e_p) c(e_{\mu}) c(e_n) c(e_v) c(e_\gamma) c(e_{\beta})\nonumber\\
	&+\mathbb{C}_3 h'(0) \sum_{j,q,v=1}^{n-1}\sum_{p,\gamma,\beta,\mu=1}^{n} X_j Y_n a_p^q a_{\gamma}^{n} a_v^{\mu} a_{\beta}^{n} \xi_j \xi_q c(e_p) c(e_{\mu}) c(e_n) c(e_v) c(e_\gamma) c(e_{\beta})\nonumber\\
	&+\mathbb{C}_3 h'(0) \sum_{j,\alpha,v=1}^{n-1}\sum_{p,\gamma,\beta,\mu=1}^{n} X_j Y_n a_p^n a_{\gamma}^{\alpha} a_v^{\mu} a_{\beta}^{n} \xi_j \xi_\alpha c(e_p) c(e_{\mu}) c(e_n) c(e_v) c(e_\gamma) c(e_{\beta})\nonumber\\
	&+\mathbb{C}_3 h'(0) \sum_{l,q,v=1}^{n-1}\sum_{p,\gamma,\beta,\mu=1}^{n} X_n Y_l a_p^q a_{\gamma}^{n} a_v^{\mu} a_{\beta}^{n} \xi_l \xi_q c(e_p) c(e_{\mu}) c(e_n) c(e_v) c(e_\gamma) c(e_{\beta})\nonumber\\
	&+\mathbb{C}_3 h'(0) \sum_{l,\alpha,v=1}^{n-1}\sum_{p,\gamma,\beta,\mu=1}^{n} X_n Y_l a_p^n a_{\gamma}^{\alpha} a_v^{\mu} a_{\beta}^{n} \xi_l \xi_\alpha c(e_p) c(e_{\mu}) c(e_n) c(e_v) c(e_\gamma) c(e_{\beta})\nonumber\\
	&+\mathbb{C}_3 h'(0) \sum_{q,\alpha,v=1}^{n-1}\sum_{p,\gamma,\beta,\mu=1}^{n} X_n Y_n a_p^q a_{\gamma}^{\alpha} a_v^{\mu} a_{\beta}^{n} \xi_q \xi_\alpha c(e_p) c(e_{\mu}) c(e_n) c(e_v) c(e_\gamma) c(e_{\beta})\nonumber\\
	&+\mathbb{C}_4 h'(0) \sum_{l,q,\alpha,k,v=1}^{n-1} \sum_{p,\gamma,\beta,\mu=1}^{n} X_n Y_l a_p^q a_{\gamma}^{\alpha} a_v^{\mu} a_{\beta}^{k} \xi_l \xi_q \xi_\alpha \xi_k c(e_p) c(e_{\mu}) c(e_n) c(e_v) c(e_\gamma) c(e_{\beta})\nonumber\\
	&+\mathbb{C}_4 h'(0) \sum_{j,q,\alpha,k,v=1}^{n-1} \sum_{p,\gamma,\beta,\mu=1}^{n} X_j Y_n a_p^q a_{\gamma}^{\alpha} a_v^{\mu} a_{\beta}^{k} \xi_j \xi_q \xi_\alpha \xi_k c(e_p) c(e_{\mu}) c(e_n) c(e_v) c(e_\gamma) c(e_{\beta})\nonumber\\
	&+\mathbb{C}_4 h'(0) \sum_{j,l,\alpha,k,v=1}^{n-1} \sum_{p,\gamma,\beta,\mu=1}^{n} X_j Y_l a_p^n a_{\gamma}^{\alpha} a_v^{\mu} a_{\beta}^{k} \xi_j \xi_l \xi_\alpha \xi_k c(e_p) c(e_{\mu}) c(e_n) c(e_v) c(e_\gamma) c(e_{\beta})\nonumber\\
	&+\mathbb{C}_4 h'(0) \sum_{j,l,q,k,v=1}^{n-1} \sum_{p,\gamma,\beta,\mu=1}^{n} X_j Y_l a_p^q a_{\gamma}^{n} a_v^{\mu} a_{\beta}^{k} \xi_j \xi_l \xi_q \xi_k c(e_p) c(e_{\mu}) c(e_n) c(e_v) c(e_\gamma) c(e_{\beta})\nonumber\\
	&+\mathbb{C}_5 h'(0) \sum_{j,l,q,\alpha,v=1}^{n-1} \sum_{p,\gamma,\beta,\mu=1}^{n} X_j Y_l a_p^q a_{\gamma}^{\alpha} a_v^{\mu} a_{\beta}^{n} \xi_j \xi_l \xi_q \xi_\alpha c(e_p) c(e_{\mu}) c(e_n) c(e_v) c(e_\gamma) c(e_{\beta}),
\end{align}
where
\begin{align}
	&\mathbb{C}_1=\frac{(3\xi_n-2i)(3\xi_n^2-1)}{16(\xi_n-i)^5 (\xi_n+i)^3},\;\;\mathbb{C}_2=-\frac{2i\xi_n^2+\xi_n}{4(\xi_n-i)^5 (\xi_n+i)^3},\nonumber\\
	&\mathbb{C}_3=\frac{\xi_n-3\xi_n^3}{16(\xi_n-i)^5 (\xi_n+i)^3},\;\;\mathbb{C}_4=-\frac{\xi_n}{4(\xi_n-i)^5 (\xi_n+i)^3},\nonumber\\
	&\mathbb{C}_5=-\frac{(\xi_n-2i)(3\xi_n^2-1)}{16(\xi_n-i)^5 (\xi_n+i)^3}.
\end{align}
Using the trace identity 
$\operatorname{Tr}(A B)=\operatorname{Tr}(B A)$ 
and the Clifford relations, we obtain
\begin{align}
	\sum_{v=1}^{n-1}\sum_{p,\mu,\gamma,\beta=1}^{n}\mathrm{Tr}\bigg(&c(e_p) c(e_{\mu}) c(e_n) c(e_v) c(e_\gamma) c(e_{\beta})\bigg)\nonumber\\
	=\sum_{v=1}^{n-1}\sum_{p,\mu,\gamma,\beta=1}^{n}\bigg(&\delta_{p}^{\beta}\delta_{\mu}^{v}\delta_{\gamma}^{n} - \delta_{p}^{\beta}\delta_{\mu}^{n}\delta_{\gamma}^{v} - \delta_{p}^{\gamma}\delta_{\mu}^{v}\delta_{\beta}^{n} + \delta_{p}^{\gamma}\delta_{\mu}^{n}\delta_{\beta}^{v} - \delta_{p}^{v}\delta_{\mu}^{\beta}\delta_{\gamma}^{n} + \delta_{p}^{v}\delta_{\mu}^{\gamma}\delta_{\beta}^{n} \nonumber\\
	&- \delta_{p}^{v}\delta_{\mu}^{n}\delta_{\beta}^{\gamma} + \delta_{p}^{n}\delta_{\mu}^{\beta}\delta_{v}^{\gamma} - \delta_{p}^{n}\delta_{\mu}^{\gamma}\delta_{\beta}^{v} + \delta_{p}^{n}\delta_{\mu}^{v}\delta_{\beta}^{\gamma} -\delta_{p}^{\mu}\delta_{\beta}^{n}\delta_{v}^{\gamma} + \delta_{p}^{\mu}\delta_{n}^{\gamma}\delta_{\beta}^{v} \bigg)\mathrm{Tr}[id],
\end{align}
it follows that
\begin{align}
&\mathrm{Tr}\bigg(\pi^+_{\xi_n}\mathcal{A}_1\times\partial_{\xi_n}\sigma_{-3}(D_{J}^{-3})\bigg)(x_0)|_{|\xi'|=1}\nonumber\\
\sim &(\mathbb{C}_1+2\mathbb{C}_2-\mathbb{C}_3)h'(0)\sum_{v=1}^{n-1} X_n Y_n a_v^v \mathrm{Tr}[id] + (\mathbb{C}_3+2\mathbb{C}_4-\mathbb{C}_5)h'(0) \sum_{j,l,v=1}^{n-1} X_j Y_l \xi_j \xi_l a_v^v \mathrm{Tr}[id].
\end{align}
Integrating over \(\xi_n\) and \(|\xi'|=1\), we obtain
\begin{align}
&-i\int_{|\xi'|=1}\int^{+\infty}_{-\infty}
\mathrm{Tr} \bigg(\pi^+_{\xi_n}\mathcal{A}_1 \times
\partial_{\xi_n}\sigma_{-3}(D_J^{-3})\bigg)(x_0)d\xi_n\sigma(\xi')dx'\nonumber\\
=&-\frac{1}{8} \sum_{v=1}^{n-1} a_v^v X_nY_nh'(0)\pi\Omega_3 dx'-\frac{1 }{12}\sum_{j,v=1}^{n-1} a_v^v X_jY_jh'(0)\pi\Omega_3 dx'.
\end{align}

Similarly, for 
$\mathcal{A}_2 (x_0)$, we have
\begin{align}
	\pi_{\xi_n}^{+}\mathcal{A}_2 (x_0)=&\;\;\;\;\frac{3i\xi_n^2+\xi_n}{16(\xi_n-i)^3}\sum_{p,\gamma=1}^{n}X_n Y_n a_p^n a_{\gamma}^{n} c(e_p)c(e_n)c(e_\gamma)\nonumber\\
	&+\frac{3i\xi_n+1}{16(\xi_n-i)^3} \sum_{p,\gamma=1}^{n}\sum_{j=1}^{n-1}\xi_j X_j Y_n a_p^n a_{\gamma}^{n} c(e_p)c(e_n)c(e_\gamma)\nonumber\\
	&+\frac{3i\xi_n+1}{16(\xi_n-i)^3}	\sum_{p,\gamma=1}^{n}\sum_{l=1}^{n-1}\xi_l X_n Y_l a_p^n a_{\gamma}^{n} c(e_p)c(e_n)c(e_\gamma)\nonumber\\
	&+\frac{3i\xi_n+1}{16(\xi_n-i)^3} \sum_{p,\gamma=1}^{n}\sum_{q=1}^{n-1}\xi_q X_n Y_n a_p^q a_{\gamma}^{n} c(e_p)c(e_n)c(e_\gamma)\nonumber\\
	&+\frac{3i\xi_n+1}{16(\xi_n-i)^3}	\sum_{p,\gamma=1}^{n}\sum_{\alpha=1}^{n-1}\xi_\alpha X_n Y_n a_p^n a_{\gamma}^{\alpha} c(e_p)c(e_n)c(e_\gamma)\nonumber\\
	&+\frac{i\xi_n^2+3\xi_n}{16(\xi_n-i)^3} \sum_{p,\gamma=1}^{n}\sum_{j,l=1}^{n-1}\xi_j \xi_l X_j Y_l a_p^n a_{\gamma}^{n} c(e_p)c(e_n)c(e_\gamma)\nonumber\\
	&+\frac{i\xi_n^2+3\xi_n}{16(\xi_n-i)^3} \sum_{p,\gamma=1}^{n}\sum_{j,q=1}^{n-1}\xi_j \xi_q X_j Y_n a_p^q a_{\gamma}^{n} c(e_p)c(e_n)c(e_\gamma)\nonumber\\
	&+\frac{i\xi_n^2+3\xi_n}{16(\xi_n-i)^3} \sum_{p,\gamma=1}^{n}\sum_{j,\alpha=1}^{n-1}\xi_j \xi_\alpha X_j Y_n a_p^n a_{\gamma}^{\alpha} c(e_p)c(e_n)c(e_\gamma)\nonumber\\
	&+\frac{i\xi_n^2+3\xi_n}{16(\xi_n-i)^3}\sum_{p,\gamma=1}^{n}\sum_{l,q=1}^{n-1} \xi_l \xi_q X_n Y_l a_p^q a_{\gamma}^{n} c(e_p)c(e_n)c(e_\gamma)\nonumber\\
	&+\frac{i\xi_n^2+3\xi_n}{16(\xi_n-i)^3}\sum_{p,\gamma=1}^{n}\sum_{l,\alpha=1}^{n-1} \xi_l \xi_\alpha X_n Y_l a_p^n a_{\gamma}^{\alpha} c(e_p)c(e_n)c(e_\gamma)\nonumber\\
	&+\frac{i\xi_n^2+3\xi_n}{16(\xi_n-i)^3}\sum_{p,\gamma=1}^{n}\sum_{q,\alpha=1}^{n-1} \xi_q \xi_\alpha X_n Y_n a_p^q a_{\gamma}^{\alpha} c(e_p)c(e_n)c(e_\gamma)\nonumber\\
	&+\frac{i\xi_n+3}{16(\xi_n-i)^3}  \sum_{p,\gamma=1}^{n} \sum_{l,q,\alpha=1}^{n-1} \xi_l \xi_q \xi_\alpha X_n Y_l a_p^q a_{\gamma}^{\alpha} c(e_p)c(e_n)c(e_\gamma)\nonumber\\
	&+\frac{i\xi_n+3}{16(\xi_n-i)^3} \sum_{p,\gamma=1}^{n} \sum_{j,q,\alpha=1}^{n-1} \xi_j \xi_q \xi_\alpha X_j Y_n a_p^q a_{\gamma}^{\alpha} c(e_p)c(e_n)c(e_\gamma)\nonumber\\
	&+\frac{i\xi_n+3}{16(\xi_n-i)^3} \sum_{p,\gamma=1}^{n}\sum_{j,l,\alpha=1}^{n-1}\xi_j \xi_l  \xi_\alpha X_j Y_l a_p^n a_{\gamma}^{\alpha} c(e_p)c(e_n)c(e_\gamma)\nonumber\\
	&+\frac{i\xi_n+3}{16(\xi_n-i)^3} \sum_{p,\gamma=1}^{n} \sum_{j,l,q=1}^{n-1} \xi_j \xi_l \xi_q  X_j Y_l a_p^q a_{\gamma}^{n} c(e_p)c(e_n)c(e_\gamma)\nonumber\\
	&+\frac{3i\xi_n^2+9\xi_n-8i}{16(\xi_n-i)^3} \sum_{p,\gamma=1}^{n} \sum_{j,l,q,\alpha=1}^{n-1} \xi_j \xi_l \xi_q \xi_{\alpha}  X_j Y_l a_p^q a_{\gamma}^{\alpha} c(e_p)c(e_n)c(e_\gamma),
\end{align}
Consequently,
\begin{align}
	&\bigg(\pi^+_{\xi_n}\mathcal{A}_2\times\partial_{\xi_n}\sigma_{-3}(D_{J}^{-3})\bigg)(x_0)|_{|\xi'|=1}\nonumber\\	
	\sim &\;\;\;\mathbb{C}_6  \sum_{p,\gamma,\beta=1}^{n} X_n Y_n a_p^n a_{\gamma}^{n} a_{\beta}^{n} c(e_p)c(e_n)c(e_\gamma) c(e_{\beta})\nonumber\\
	&+\mathbb{C}_7  \sum_{j,k=1}^{n-1}\sum_{p,\gamma,\beta=1}^{n} X_j Y_n a_p^n a_{\gamma}^{n} a_{\beta}^{k} \xi_j \xi_k c(e_p)c(e_n)c(e_\gamma) c(e_{\beta})\nonumber\\
	&+\mathbb{C}_7  \sum_{l,k=1}^{n-1}\sum_{p,\gamma,\beta=1}^{n} X_n Y_l a_p^n a_{\gamma}^{n} a_{\beta}^{k} \xi_l \xi_k c(e_p)c(e_n) c(e_\gamma) c(e_{\beta})\nonumber\\
	&+\mathbb{C}_7  \sum_{k,q=1}^{n-1}\sum_{p,\gamma,\beta=1}^{n} X_n Y_n a_p^q a_{\gamma}^{n} a_{\beta}^{k} \xi_q \xi_k c(e_p)c(e_n)c(e_\gamma) c(e_{\beta})\nonumber\\
	&+\mathbb{C}_7  \sum_{k,\alpha=1}^{n-1}\sum_{p,\gamma,\beta=1}^{n} X_n Y_n a_p^n a_{\gamma}^{\alpha} a_{\beta}^{k} \xi_{\alpha} \xi_k c(e_p)c(e_n)c(e_\gamma) c(e_{\beta})\nonumber\\
	&+\mathbb{C}_8  \sum_{j,l=1}^{n-1}\sum_{p,\gamma,\beta=1}^{n} X_j Y_l a_p^n a_{\gamma}^{n} a_{\beta}^{n} \xi_j \xi_l c(e_p)c(e_n)c(e_\gamma) c(e_{\beta})\nonumber\\
	&+\mathbb{C}_8  \sum_{j,q=1}^{n-1}\sum_{p,\gamma,\beta=1}^{n} X_j Y_n a_p^q a_{\gamma}^{n} a_{\beta}^{n} \xi_j \xi_q c(e_p)c(e_n)c(e_\gamma) c(e_{\beta})\nonumber\\
	&+\mathbb{C}_8  \sum_{j,\alpha=1}^{n-1}\sum_{p,\gamma,\beta=1}^{n} X_j Y_n a_p^n a_{\gamma}^{\alpha} a_{\beta}^{n} \xi_j \xi_\alpha c(e_p)c(e_n)c(e_\gamma) c(e_{\beta})\nonumber\\
	&+\mathbb{C}_8  \sum_{l,q=1}^{n-1}\sum_{p,\gamma,\beta=1}^{n} X_n Y_l a_p^q a_{\gamma}^{n} a_{\beta}^{n} \xi_l \xi_q c(e_p)c(e_n)c(e_\gamma)c(e_{\beta})\nonumber\\
	&+\mathbb{C}_8  \sum_{l,\alpha=1}^{n-1}\sum_{p,\gamma,\beta=1}^{n} X_n Y_l a_p^n a_{\gamma}^{\alpha} a_{\beta}^{n} \xi_l \xi_\alpha c(e_p)c(e_n)c(e_\gamma)c(e_{\beta})\nonumber\\
	&+\mathbb{C}_8  \sum_{q,\alpha=1}^{n-1}\sum_{p,\gamma,\beta=1}^{n} X_n Y_n a_p^q a_{\gamma}^{\alpha} a_{\beta}^{n} \xi_q \xi_\alpha c(e_p)c(e_n)c(e_\gamma)c(e_{\beta})\nonumber\\
	&+\mathbb{C}_9  \sum_{l,q,\alpha,k=1}^{n-1} \sum_{p,\gamma,\beta=1}^{n} X_n Y_l a_p^q a_{\gamma}^{\alpha} a_{\beta}^{k} \xi_l \xi_q \xi_\alpha \xi_k c(e_p)c(e_n)c(e_\gamma) c(e_{\beta})\nonumber\\
	&+\mathbb{C}_9  \sum_{j,q,\alpha,k=1}^{n-1} \sum_{p,\gamma,\beta=1}^{n} X_j Y_n a_p^q a_{\gamma}^{\alpha} a_{\beta}^{k} \xi_j \xi_q \xi_\alpha \xi_k c(e_p)c(e_n)c(e_\gamma)c(e_{\beta})\nonumber\\
	&+\mathbb{C}_9  \sum_{j,l,\alpha,k=1}^{n-1} \sum_{p,\gamma,\beta=1}^{n} X_j Y_l a_p^n a_{\gamma}^{\alpha} a_{\beta}^{k} \xi_j \xi_l \xi_\alpha \xi_k c(e_p)c(e_n)c(e_\gamma) c(e_{\beta})\nonumber\\
	&+\mathbb{C}_9  \sum_{j,l,q,k=1}^{n-1} \sum_{p,\gamma,\beta=1}^{n} X_j Y_l a_p^q a_{\gamma}^{n} a_{\beta}^{k} \xi_j \xi_l \xi_q \xi_k c(e_p)c(e_n)c(e_\gamma) c(e_{\beta})\nonumber\\
	&+\mathbb{C}_{10}  \sum_{j,l,q,\alpha=1}^{n-1} \sum_{p,\gamma,\beta=1}^{n} X_j Y_l a_p^q a_{\gamma}^{\alpha} a_{\beta}^{n} \xi_j \xi_l \xi_q \xi_\alpha c(e_p)c(e_n)c(e_\gamma)c(e_{\beta}),
\end{align}
where
\begin{align}
	&\mathbb{C}_6=\frac{(3\xi_n^2-i\xi_n)(3\xi_n^2-1)}{16(\xi_n-i)^6 (\xi_n+i)^3},\;\;\mathbb{C}_7=\frac{3\xi_n^2-i\xi_n}{4(\xi_n-i)^6 (\xi_n+i)^3},\nonumber\\
	&\mathbb{C}_8=\frac{(\xi_n^2-3i\xi_n)(3\xi_n^2-1)}{16(\xi_n-i)^6 (\xi_n+i)^3},\;\;\mathbb{C}_9=\frac{\xi_n^2-3i\xi_n}{4(\xi_n-i)^6 (\xi_n+i)^3},\nonumber\\
	&\mathbb{C}_{10}=\frac{(3\xi_n^2-9i\xi_n-8)(3\xi_n^2-1)}{16(\xi_n-i)^6 (\xi_n+i)^3}.
\end{align}
Using the Clifford trace identity, we have
\begin{align}
	\sum_{p,\gamma,\beta=1}^{n}\mathrm{Tr}\bigg(c(e_p) c(e_n) c(e_\gamma) c(e_{\beta})\bigg)
	=\sum_{p,\gamma,\beta=1}^{n}\bigg( \delta_{p}^{\beta}\delta_{\gamma}^{n} -  \delta_{\gamma}^{p}\delta_{\beta}^{n} +  \delta_{p}^{n}\delta_{\gamma}^{\beta} \bigg)\mathrm{Tr}[id],
\end{align}
It follows that
\begin{align}
	&\mathrm{Tr}\bigg(\pi^+_{\xi_n}\mathcal{A}_2\times\partial_{\xi_n}\sigma_{-3}(D_{J}^{-3})\bigg)(x_0)|_{|\xi'|=1}\nonumber\\
	=&(\mathbb{C}_6+2\mathbb{C}_7-\mathbb{C}_8) X_n Y_n \mathrm{Tr}[id] + (\mathbb{C}_8+2\mathbb{C}_9-\mathbb{C}_{10}) \sum_{j,l=1}^{n-1} X_j Y_l \xi_j \xi_l \mathrm{Tr}[id] 
\end{align}
Integrating over $\xi_n$ and $|\xi'|=1$, we obtain
\begin{align}
	&-i\int_{|\xi'|=1}\int^{+\infty}_{-\infty}
	\mathrm{Tr} \bigg(\pi^+_{\xi_n}\mathcal{A}_2 \times
	\partial_{\xi_n}\sigma_{-3}(D_J^{-3})\bigg)(x_0)d\xi_n\sigma(\xi')dx'\nonumber\\
	=&\;\;\;\frac{7}{32} X_nY_n \pi\Omega_3 dx'+\frac{25}{96}\sum_{j=1}^{n-1} X_j Y_j \pi\Omega_3 dx'.
\end{align}

Similarly, for 
$\mathcal{A}_3 (x_0)$, we have
\begin{align}
	\pi_{\xi_n}^{+}\mathcal{A}_3 (x_0)=
	&-\frac{2\xi_n-i}{4(\xi_n-i)^2}	\sum_{\theta,y,\gamma,h=1}^{n} \sum_{p=1}^{n-1}\xi_p X_n Y_n a_{\theta}^{n} a_{y}^{\gamma}\partial_{x_{\gamma}}(a_h^p) c(e_\theta)c(e_y)c(e_h)\nonumber\\
	&+\frac{i\xi_n}{4(\xi_n-i)^2} \sum_{\theta,y,\gamma,h=1}^{n} \sum_{j,p=1}^{n-1}\xi_j \xi_p X_j Y_n a_{\theta}^{n} a_{y}^{\gamma}\partial_{x_{\gamma}}(a_h^p) c(e_\theta)c(e_y)c(e_h)\nonumber\\
	&+\frac{i\xi_n}{4(\xi_n-i)^2} \sum_{\theta,y,\gamma,h=1}^{n} \sum_{l,p=1}^{n-1} \xi_l \xi_p X_n Y_l a_{\theta}^{n} a_{y}^{\gamma}\partial_{x_{\gamma}}(a_h^p) c(e_\theta)c(e_y)c(e_h)\nonumber\\
	&+\frac{i\xi_n}{4(\xi_n-i)^2}\sum_{\theta,y,\gamma,h=1}^{n} \sum_{q,p=1}^{n-1} \xi_q \xi_p X_n Y_n a_{\theta}^{q} a_{y}^{\gamma}\partial_{x_{\gamma}}(a_h^p) c(e_\theta)c(e_y)c(e_h)\nonumber\\
	&+\frac{i}{4(\xi_n-i)^2} \sum_{\theta,y,\gamma,h=1}^{n} \sum_{l,q,p=1}^{n-1} \xi_l \xi_q \xi_p X_n Y_l a_{\theta}^{q} a_{y}^{\gamma}\partial_{x_{\gamma}}(a_h^p) c(e_\theta)c(e_y)c(e_h)\nonumber\\
	&+\frac{i}{4(\xi_n-i)^2} \sum_{\theta,y,\gamma,h=1}^{n} \sum_{j,q,p=1}^{n-1} \xi_j \xi_q \xi_p X_j Y_n a_{\theta}^{q} a_{y}^{\gamma}\partial_{x_{\gamma}}(a_h^p) c(e_\theta)c(e_y)c(e_h)\nonumber\\
	&+\frac{i}{4(\xi_n-i)^2} \sum_{\theta,y,\gamma,h=1}^{n} \sum_{j,l,p=1}^{n-1}\xi_j \xi_l  \xi_p X_j Y_l a_{\theta}^{n} a_{y}^{\gamma}\partial_{x_{\gamma}}(a_h^p) c(e_\theta)c(e_y)c(e_h)\nonumber\\
	&+\frac{i\xi_n+2}{4(\xi_n-i)^2} \sum_{\theta,y,\gamma,h=1}^{n} \sum_{j,l,q,p=1}^{n-1} \xi_j \xi_l \xi_q \xi_{p}  X_j Y_l a_{\theta}^{q} a_{y}^{\gamma}\partial_{x_{\gamma}}(a_h^p) c(e_\theta)c(e_y)c(e_h).
\end{align}
Consequently,
\begin{align}
	&\bigg(\pi^+_{\xi_n}\mathcal{A}_3\times\partial_{\xi_n}\sigma_{-3}(D_{J}^{-3})\bigg)(x_0)|_{|\xi'|=1}\nonumber\\	
	\sim &\;\;\;\mathbb{C}_{11}	\sum_{\theta,y,\gamma,h,\beta=1}^{n} \sum_{p,k=1}^{n-1}\xi_p \xi_k X_n Y_n a_{\theta}^{n} a_{y}^{\gamma}\partial_{x_{\gamma}}(a_h^p)a_{\beta}^{k} c(e_\theta)c(e_y)c(e_h)c(e_{\beta})\nonumber\\
	&+\mathbb{C}_{12} \sum_{\theta,y,\gamma,h,\beta=1}^{n} \sum_{j,p=1}^{n-1}\xi_j \xi_p X_j Y_n a_{\theta}^{n} a_{y}^{\gamma}\partial_{x_{\gamma}}(a_h^p)a_{\beta}^{n} c(e_\theta)c(e_y)c(e_h)c(e_{\beta})\nonumber\\
	&+\mathbb{C}_{12} \sum_{\theta,y,\gamma,h,\beta=1}^{n} \sum_{l,p=1}^{n-1} \xi_l \xi_p X_n Y_l a_{\theta}^{n} a_{y}^{\gamma}\partial_{x_{\gamma}}(a_h^p)a_{\beta}^{n} c(e_\theta)c(e_y)c(e_h)c(e_{\beta})\nonumber\\
	&+\mathbb{C}_{12} \sum_{\theta,y,\gamma,h,\beta=1}^{n} \sum_{q,p=1}^{n-1} \xi_q \xi_p X_n Y_n a_{\theta}^{q} a_{y}^{\gamma}\partial_{x_{\gamma}}(a_h^p)a_{\beta}^{n} c(e_\theta)c(e_y)c(e_h)c(e_{\beta})\nonumber\\
	&+\mathbb{C}_{13} \sum_{\theta,y,\gamma,h,\beta=1}^{n} \sum_{l,q,p,k=1}^{n-1} \xi_l \xi_q \xi_p \xi_k X_n Y_l a_{\theta}^{q} a_{y}^{\gamma}\partial_{x_{\gamma}}(a_h^p)a_{\beta}^{k} c(e_\theta)c(e_y)c(e_h)c(e_{\beta})\nonumber\\
	&+\mathbb{C}_{13} \sum_{\theta,y,\gamma,h,\beta=1}^{n} \sum_{j,q,p,k=1}^{n-1} \xi_j \xi_q \xi_p \xi_k X_j Y_n a_{\theta}^{q} a_{y}^{\gamma}\partial_{x_{\gamma}}(a_h^p)a_{\beta}^{k} c(e_\theta)c(e_y)c(e_h)c(e_{\beta})\nonumber\\
	&+\mathbb{C}_{13} \sum_{\theta,y,\gamma,h,\beta=1}^{n} \sum_{j,l,p,k=1}^{n-1}\xi_j \xi_l \xi_p \xi_k X_j Y_l a_{\theta}^{n} a_{y}^{\gamma}\partial_{x_{\gamma}}(a_h^p)a_{\beta}^{k} c(e_\theta)c(e_y)c(e_h)c(e_{\beta})\nonumber\\
	&+\mathbb{C}_{14} \sum_{\theta,y,\gamma,h,\beta=1}^{n} \sum_{j,l,q,p=1}^{n-1} \xi_j \xi_l \xi_q \xi_{p}  X_j Y_l a_{\theta}^{q} a_{y}^{\gamma}\partial_{x_{\gamma}}(a_h^p) a_{\beta}^{n} c(e_\theta)c(e_y)c(e_h)c(e_{\beta}),
\end{align}
where
\begin{align}
	&\mathbb{C}_{11}=\frac{2i\xi_n^2+\xi_n}{(\xi_n-i)^5 (\xi_n+i)^3},\;\;\mathbb{C}_{12}=\frac{3\xi_n^3-\xi_n}{4(\xi_n-i)^5 (\xi_n+i)^3},\nonumber\\
	&\mathbb{C}_{13}=\frac{\xi_n}{(\xi_n-i)^5 (\xi_n+i)^3},\;\;\mathbb{C}_{14}=\frac{(\xi_n-2i)(3\xi_n^2-1)}{4(\xi_n-i)^5 (\xi_n+i)^3}.
\end{align}
Using the Clifford trace identity, we have
\begin{align}
	\sum_{\theta,y,h,\beta=1}^{n}\mathrm{Tr}\bigg(c(e_\theta) c(e_y) c(e_h) c(e_{\beta})\bigg)
	=\sum_{\theta,y,h,\beta=1}^{n}\bigg( \delta_{\theta}^{\beta}\delta_{h}^{y} -  \delta_{h}^{\theta}\delta_{\beta}^{y} +  \delta_{\theta}^{y}\delta_{h}^{\beta} \bigg)\mathrm{Tr}[id],
\end{align}
It follows that
\begin{align}
	&\mathrm{Tr}\bigg(\pi^+_{\xi_n}\mathcal{A}_3\times\partial_{\xi_n}\sigma_{-3}(D_{J}^{-3})\bigg)(x_0)|_{|\xi'|=1}\nonumber\\
	=&\;\;\bigg( \mathbb{C}_{12}+ \mathbb{C}_{13} \bigg) \sum_{y,\gamma=1}^{n} \sum_{j,p=1}^{n-1} \xi_j \xi_p \bigg( X_j Y_n + X_n Y_j \bigg) a_y^{\gamma} \partial_{x_{\gamma}}(a_y^p) \mathrm{Tr}[id]\nonumber\\
	& - \mathbb{C}_{14}\sum_{\theta=1}^{n} \sum_{j,l,p,q=1}^{n-1} \xi_j \xi_l \xi_p \xi_q X_j Y_l a_{\theta}^q \partial_{x_n}(a_{\theta}^p) \mathrm{Tr}[id] 
\end{align}
Integrating over $\xi_n$ and $|\xi'|=1$, we obtain
\begin{align}
	&-i\int_{|\xi'|=1}\int^{+\infty}_{-\infty}
	\mathrm{Tr} \bigg(\pi^+_{\xi_n}\mathcal{A}_3 \times
	\partial_{\xi_n}\sigma_{-3}(D_J^{-3})\bigg)(x_0)d\xi_n\sigma(\xi')dx'\nonumber\\
	=&\;\;\;\frac{1}{8} \sum_{y,\gamma=1}^{n} \sum_{j=1}^{n-1} \bigg( X_j Y_n + X_n Y_j \bigg) a_y^{\gamma} \partial_{x_{\gamma}}(a_y^j)\pi\Omega_3 dx'.
\end{align}

Finally, for $\mathcal{A}_4 (x_0)$, we have
\begin{align}
	\pi_{\xi_n}^{+}\mathcal A_4(x_0)
	={}&-\frac{h'(0)(2\xi_n-i)}{8(\xi_n-i)^2}
	\sum_{\theta=1}^{n}\sum_{p,h=1}^{n-1}
	X_nY_n\xi_p a_\theta^n a_h^p
	c(e_\theta)c(e_n)c(e_h)
	\nonumber\\
	&+\frac{ih'(0)\xi_n}{8(\xi_n-i)^2}
	\sum_{\theta=1}^{n}\sum_{j,p,h=1}^{n-1}
	X_jY_n\xi_j\xi_p a_\theta^n a_h^p
	c(e_\theta)c(e_n)c(e_h)
	\nonumber\\
	&+\frac{ih'(0)\xi_n}{8(\xi_n-i)^2}
	\sum_{\theta=1}^{n}\sum_{l,p,h=1}^{n-1}
	X_nY_l\xi_l\xi_p a_\theta^n a_h^p
	c(e_\theta)c(e_n)c(e_h)
	\nonumber\\
	&+\frac{ih'(0)\xi_n}{8(\xi_n-i)^2}
	\sum_{\theta=1}^{n}\sum_{q,p,h=1}^{n-1}
	X_nY_n\xi_q\xi_p a_\theta^q a_h^p
	c(e_\theta)c(e_n)c(e_h)
	\nonumber\\
	&+\frac{ih'(0)}{8(\xi_n-i)^2}
	\sum_{\theta=1}^{n}\sum_{l,q,p,h=1}^{n-1}
	X_nY_l\xi_l\xi_q\xi_p a_\theta^q a_h^p
	c(e_\theta)c(e_n)c(e_h)
	\nonumber\\
	&+\frac{ih'(0)}{8(\xi_n-i)^2}
	\sum_{\theta=1}^{n}\sum_{j,q,p,h=1}^{n-1}
	X_jY_n\xi_j\xi_q\xi_p a_\theta^q a_h^p
	c(e_\theta)c(e_n)c(e_h)
	\nonumber\\
	&+\frac{ih'(0)}{8(\xi_n-i)^2}
	\sum_{\theta=1}^{n}\sum_{j,l,p,h=1}^{n-1}
	X_jY_l\xi_j\xi_l\xi_p a_\theta^n a_h^p
	c(e_\theta)c(e_n)c(e_h)
	\nonumber\\
	&+\frac{h'(0)(i\xi_n+2)}{8(\xi_n-i)^2}
	\sum_{\theta=1}^{n}\sum_{j,l,q,p,h=1}^{n-1}
	X_jY_l\xi_j\xi_l\xi_q\xi_p a_\theta^q a_h^p
	c(e_\theta)c(e_n)c(e_h).
\end{align}
Consequently,
\begin{align}\label{5.64}
	&\left(
	\pi_{\xi_n}^{+}\mathcal A_4
	\times\partial_{\xi_n}\sigma_{-3}(D_J^{-3})
	\right)(x_0)\big|_{|\xi'|=1}
	\nonumber\\
	&\sim
	-2h'(0)\mathbb C_2
	\sum_{p,\gamma,\beta=1}^{n}
	\sum_{\alpha,k=1}^{n-1}
	\xi_\alpha\xi_kX_nY_n
	a_p^n a_\gamma^\alpha a_\beta^k
	c(e_p)c(e_n)c(e_\gamma)c(e_\beta)
	\nonumber\\
	&\quad
	-2h'(0)\mathbb C_3
	\sum_{p,\gamma,\beta=1}^{n}
	\sum_{j,\alpha=1}^{n-1}
	\xi_j\xi_\alpha X_jY_n
	a_p^n a_\gamma^\alpha a_\beta^n
	c(e_p)c(e_n)c(e_\gamma)c(e_\beta)
	\nonumber\\
	&\quad
	-2h'(0)\mathbb C_3
	\sum_{p,\gamma,\beta=1}^{n}
	\sum_{l,\alpha=1}^{n-1}
	\xi_l\xi_\alpha X_nY_l
	a_p^n a_\gamma^\alpha a_\beta^n
	c(e_p)c(e_n)c(e_\gamma)c(e_\beta)
	\nonumber\\
	&\quad
	-2h'(0)\mathbb C_3
	\sum_{p,\gamma,\beta=1}^{n}
	\sum_{q,\alpha=1}^{n-1}
	\xi_q\xi_\alpha X_nY_n
	a_p^q a_\gamma^\alpha a_\beta^n
	c(e_p)c(e_n)c(e_\gamma)c(e_\beta)
	\nonumber\\
	&\quad
	-2h'(0)\mathbb C_4
	\sum_{p,\gamma,\beta=1}^{n}
	\sum_{l,q,\alpha,k=1}^{n-1}
	\xi_l\xi_q\xi_\alpha\xi_kX_nY_l
	a_p^q a_\gamma^\alpha a_\beta^k
	c(e_p)c(e_n)c(e_\gamma)c(e_\beta)
	\nonumber\\
	&\quad
	-2h'(0)\mathbb C_4
	\sum_{p,\gamma,\beta=1}^{n}
	\sum_{j,q,\alpha,k=1}^{n-1}
	\xi_j\xi_q\xi_\alpha\xi_kX_jY_n
	a_p^q a_\gamma^\alpha a_\beta^k
	c(e_p)c(e_n)c(e_\gamma)c(e_\beta)
	\nonumber\\
	&\quad
	-2h'(0)\mathbb C_4
	\sum_{p,\gamma,\beta=1}^{n}
	\sum_{j,l,\alpha,k=1}^{n-1}
	\xi_j\xi_l\xi_\alpha\xi_kX_jY_l
	a_p^n a_\gamma^\alpha a_\beta^k
	c(e_p)c(e_n)c(e_\gamma)c(e_\beta)
	\nonumber\\
	&\quad
	-2h'(0)\mathbb C_5
	\sum_{p,\gamma,\beta=1}^{n}
	\sum_{j,l,q,\alpha=1}^{n-1}
	\xi_j\xi_l\xi_q\xi_\alpha X_jY_l
	a_p^q a_\gamma^\alpha a_\beta^n
	c(e_p)c(e_n)c(e_\gamma)c(e_\beta).
\end{align}
\begin{align}
	&\operatorname{Tr}\left[
	\pi_{\xi_n}^{+}\mathcal A_4\,
	\partial_{\xi_n}\sigma_{-3}(D_J^{-3})
	\right](x_0)\big|_{|\xi'|=1}
	\nonumber\\
	&\sim
	8h'(0)(\mathbb C_3-\mathbb C_2)X_nY_n
	+
	8h'(0)(\mathbb C_5-\mathbb C_4)
	\sum_{j,l=1}^{n-1}X_jY_l\xi_j\xi_l
	\nonumber\\
	&=
	\frac{h'(0)(-3\xi_n^3+8i\xi_n^2+5\xi_n)}
	{2(\xi_n-i)^5(\xi_n+i)^3}X_nY_n
	+
	\frac{h'(0)(-3\xi_n^3+6i\xi_n^2+5\xi_n-2i)}
	{2(\xi_n-i)^5(\xi_n+i)^3}
	\sum_{j,l=1}^{n-1}X_jY_l\xi_j\xi_l.
\end{align}
Taking the trace and integrating over $\xi_n$ and $|\xi'|=1$, we obtain
\begin{align}
	&-i\int_{|\xi'|=1}\int^{+\infty}_{-\infty}
	\mathrm{Tr} \bigg(\pi^+_{\xi_n}\mathcal{A}_4 \times
	\partial_{\xi_n}\sigma_{-3}(D_J^{-3})\bigg)(x_0)d\xi_n\sigma(\xi')dx'\nonumber\\
	=&\;\;\;\frac{1}{16} X_nY_n h'(0)\pi\Omega_3 dx'+\frac{5}{48}\sum_{j=1}^{n-1} X_jY_jh'(0)\pi\Omega_3 dx'.
\end{align}
Combining (4.42), (4.48), (4.54), and (4.58), we obtain
\begin{align}
	&-i\int_{|\xi'|=1}\int^{+\infty}_{-\infty}
	\mathrm{Tr} \bigg(\pi^+_{\xi_n}(\mathcal{A}_1 -h'(0) \mathcal{A}_2 + \mathcal{A}_3  + \mathcal{A}_4) \times
	\partial_{\xi_n}\sigma_{-3}(D_J^{-3})\bigg)(x_0)d\xi_n\sigma(\xi')dx'\nonumber\\
	=&-\frac{1}{8} \sum_{v=1}^{n-1} a_v^v X_nY_nh'(0)\pi\Omega_3 dx'-\frac{1 }{12}\sum_{j,v=1}^{n-1} a_v^v X_jY_jh'(0)\pi\Omega_3 dx'\nonumber\\
	&-\frac{5}{32} X_nY_n  h'(0) \pi\Omega_3 dx'-\frac{5}{32}\sum_{j=1}^{n-1} X_j Y_j  h'(0) \pi\Omega_3 dx'\nonumber\\
	&+\frac{1}{8} \sum_{y,\gamma=1}^{n} \sum_{j=1}^{n-1} \bigg( X_j Y_n + X_n Y_j \bigg) a_y^{\gamma} \partial_{x_{\gamma}}(a_y^j)\pi\Omega_3 dx'.
\end{align}

(2) For the second term in \eqref{b41}, we have
\begin{align}
	&\sigma_{1}(\widetilde{\nabla}^J_{X}\widetilde{\nabla}^J_{Y} )\sigma_{-1}(D_{J}^{-1})(x_0)|_{|\xi'|=1}\nonumber\\
	=&-\bigg[\sum_{j,l=1}^{n}X_j\frac{\partial{Y_l}}{\partial{x_j}}\xi_l
	+\sum_{l=1}^{n}\bigg(\big(M(Y)+N(Y)\big)X_l\xi_l +\big(M(X)+N(X)\big)Y_l\xi_l\bigg)\bigg]\frac{c[J(\xi)]}{1+\xi_n^{2}}\nonumber\\
	=&-\frac{1}{1+\xi_n^{2}}\sum_{h=1}^{n}\sum_{l,k=1}^{n-1}\bigg(\sum_{j=1}^{n}X_j\frac{\partial{Y_l}}{\partial{x_j}}
	+\big(M(Y)+N(Y)\big)X_l +\big(M(X)+N(X)\big)Y_l\bigg)a_h^k c(e_h) \xi_l \xi_k\nonumber\\
	&-\frac{\xi_n}{1+\xi_n^{2}}\sum_{l=1}^{n-1}\bigg(\sum_{j=1}^{n}X_j\frac{\partial{Y_l}}{\partial{x_j}}
	+\big(M(Y)+N(Y)\big)X_l +\big(M(X)+N(X)\big)Y_l\bigg) c(e_n) \xi_l \nonumber\\
	&-\frac{\xi_n}{1+\xi_n^{2}}\sum_{h=1}^{n}\sum_{k=1}^{n-1}\bigg(\sum_{j=1}^{n}X_j\frac{\partial_{Y_n}}{\partial_{x_j}}
	+\big(M(Y)+N(Y)\big)X_n +\big(M(X)+N(X)\big)Y_n\bigg)a_h^k c(e_h) \xi_k\nonumber\\
	&-\frac{\xi_n^2}{1+\xi_n^{2}}\bigg(\sum_{j=1}^{n}X_j\frac{\partial_{Y_n}}{\partial_{x_j}}
	+\big(M(Y)+N(Y)\big)X_n +\big(M(X)+N(X)\big)Y_n\bigg) c(e_n),
\end{align}
It follows that
\begin{align}
	&\left(\pi^+_{\xi_n}\Big(\sigma_{1}( \widetilde{\nabla}^J_{X}\widetilde{\nabla}^J_{Y} )\sigma_{-1}(D_{J}^{-1})\Big)\times
	\partial_{\xi_n}\sigma_{-3}(D_J^{-3})\right)(x_0)|_{|\xi'|=1}\nonumber\\
	\sim &\;\;\mathbb{H}_1 \sum_{h=1}^{n}\sum_{l,k=1}^{n-1}\bigg(\sum_{j=1}^{n}X_j\frac{\partial{Y_l}}{\partial{x_j}}
	+\big(M(Y)+N(Y)\big)X_l +\big(M(X)+N(X)\big)Y_l\bigg)a_h^k c(e_h) c(e_n) \xi_l \xi_k\nonumber\\
	&+\mathbb{H}_2 \sum_{\beta=1}^{n}\sum_{l,\alpha=1}^{n-1}\bigg(\sum_{j=1}^{n}X_j\frac{\partial{Y_l}}{\partial{x_j}}
	+\big(M(Y)+N(Y)\big)X_l +\big(M(X)+N(X)\big)Y_l\bigg) a_{\beta}^{\alpha} c(e_n) c(e_\beta) \xi_l \xi_{\alpha} \nonumber\\
	&+\mathbb{H}_2 \sum_{h,\beta=1}^{n}\sum_{k,\alpha=1}^{n-1}\bigg(\sum_{j=1}^{n}X_j\frac{\partial_{Y_n}}{\partial_{x_j}}
	+\big(M(Y)+N(Y)\big)X_n +\big(M(X)+N(X)\big)Y_n\bigg)a_h^k a_{\beta}^{\alpha} c(e_h) c(e_\beta) \xi_k \xi_{\alpha} \nonumber\\
	&+\mathbb{H}_1 \bigg(\sum_{j=1}^{n}X_j\frac{\partial_{Y_n}}{\partial_{x_j}}
	+\big(M(Y)+N(Y)\big)X_n +\big(M(X)+N(X)\big)Y_n\bigg),
\end{align}
where
\begin{align}
	\mathbb{H}_1:=\frac{3\xi_n^2-1}{2(\xi_n-i)^4(\xi_n+i)^3},\;\;
	\mathbb{H}_2:=\frac{2i\xi_n}{(\xi_n-i)^4(\xi_n+i)^3}.\nonumber
\end{align}

Since 
\begin{align}
	\int_{|\xi'|=1}\xi_{\beta_{1}}\xi_{\beta_{2}}\cdots\xi_{\beta_{2d+1}}\sigma(\xi')=0,~\beta_i<n,\nonumber
\end{align}
the terms containing odd monomials in $\xi'$ make no contribution after integration over $|\xi'|=1$. We obtain
\begin{align}
&\mathrm{Tr} \left(\pi^+_{\xi_n}\Big(\sigma_{1}(\widetilde{\nabla}^J_{X}\widetilde{\nabla}^J_{Y} )\sigma_{-1}(D_{J}^{-1})\Big)\times
\partial_{\xi_n}\sigma_{-3}(D_J^{-3})\right)(x_0)|_{|\xi'|=1}\nonumber\\
\sim &\;\;4(\mathbb{H}_1-\mathbb{H}_2)\sum_{j=1}^{n} X_j\frac{\partial_{Y_n}}{\partial_{x_j}}+(\mathbb{H}_2-\mathbb{H}_1) h'(0)  \sum_{\alpha,l,k=1}^{n-1} (X_l Y_{\alpha}+X_{\alpha} Y_l) a_{\alpha}^k \xi_l \xi_k \nonumber\\
&-2(\mathbb{H}_1-\mathbb{H}_2)\sum_{h=1}^{n} \sum_{l,k=1}^{n-1} a_h^k g\big((\nabla^L_{e_{n}}J)Y,e_h\big)X_l \xi_l \xi_k+2(\mathbb{H}_1-\mathbb{H}_2)\sum_{l,k=1}^{n-1}  g\big((\nabla^L_{e_{k}}J)Y,e_n\big)X_l \xi_l \xi_k\nonumber\\
&-2(\mathbb{H}_1-\mathbb{H}_2)\sum_{h=1}^{n} \sum_{l,k=1}^{n-1} a_h^k g\big((\nabla^L_{e_{n}}J)X,e_h\big)Y_l \xi_l \xi_k+2(\mathbb{H}_1-\mathbb{H}_2)\sum_{l,k=1}^{n-1}  g\big((\nabla^L_{e_{k}}J)X,e_n\big)Y_l \xi_l \xi_k\nonumber\\
&+2(\mathbb{H}_1-\mathbb{H}_2)\sum_{\alpha=1}^{n} g\big((\nabla^L_{e_{\alpha}}J)Y,J(e_\alpha)\big)X_n+2(\mathbb{H}_1-\mathbb{H}_2)\sum_{\alpha=1}^{n} g\big((\nabla^L_{e_{\alpha}}J)X,J(e_\alpha)\big)Y_n.
\end{align}
Integrating over $\xi_n$ and $|\xi'|=1$, we obtain
\begin{align}
	&-i\int_{|\xi'|=1}\int^{+\infty}_{-\infty}
	\mathrm{Tr} [\pi^+_{\xi_n}\sigma_{1}(\widetilde{\nabla}^J_{X}\widetilde{\nabla}^J_{Y} )\sigma_{-1}(D_{J}^{-1})\times
	\partial_{\xi_n}\sigma_{-3}(D_{J}^{-3})](x_0)d\xi_n\sigma(\xi')dx'\nonumber\\
	=&-\frac{3}{4}X(Y_n)\pi\Omega_3dx'+ \frac{1}{8} h'(0) g\big(JX,Y\big)\pi\Omega_3dx' - \frac{1}{8}X_n Y_n h'(0) \pi\Omega_3dx'\nonumber\\
	&+ \frac{1}{8}g\big((\nabla^L_{e_{n}}J)Y,J(X^T)\big)\pi\Omega_3dx'+ \frac{1}{8}g\big((\nabla^L_{e_{n}}J)X,J(Y^T)\big)\pi\Omega_3dx'\nonumber\\
	&- \frac{1}{8}g\big((\nabla^L_{X^T}J)Y,e_n\big)\pi\Omega_3dx' - \frac{1}{8}g\big((\nabla^L_{Y^T}J)X,e_n\big)\pi\Omega_3dx'\nonumber\\
	&- \frac{3}{8}\sum_{\alpha=1}^{n} g\big((\nabla^L_{e_{\alpha}}J)Y,J(e_\alpha)\big)X_n \pi\Omega_3dx' - \frac{3}{8}\sum_{\alpha=1}^{n} g\big((\nabla^L_{e_{\alpha}}J)X,J(e_\alpha)\big)Y_n \pi\Omega_3dx'.
\end{align}
	where $X(Y_n)$ is $\sum_{j=1}^nX_j\frac{\partial_{Y_n}}{\partial_{x_j}}.$
	
(3) For the third term in \eqref{b41}, we have
\begin{align}
	&\sum_{j=1}^{n}\partial_{\xi_j}\big[\sigma_{2}(\widetilde{\nabla}^J_{X}\widetilde{\nabla}^J_{Y})\big]
	D_{x_j}\big[\sigma_{-1}(D_{J}^{-1})\big](x_0)|_{|\xi'|=1}\nonumber\\
	=&\;\;\sum_{j=1}^{n}\sum_{l=1}^{n}\sqrt{-1}(X_{j}Y_l+X_{l}Y_j)\xi_{l}\partial_{x_{j}}(\frac{\sqrt{-1}c[J(\xi)]}{|\xi|^{2}})\nonumber\\
	=&\;\;\frac{\xi_n}{(1+\xi_n^2)^2} \sum_{l=1}^{n-1} (X_{n} Y_l + X_{l} Y_n ) h'(0) c(e_n) \xi_{l} - \frac{\xi_n}{1+\xi_n^2} \sum_{j,h=1}^{n} \sum_{k=1}^{n-1} (X_{j} Y_n + X_{n} Y_j ) \partial_{x_j}(a_h^k) c(e_h) \xi_{k}\nonumber\\
	&-\frac{2\xi_n}{1+\xi_n^2}\sum_{h=1}^{n}\sum_{k=1}^{n-1}X_n Y_n a_h^k \left.\partial_{x_n}c(dx_h)\right|_{x_0}\xi_{k} + \frac{2\xi_n}{(1+\xi_n^2)^2}\sum_{h=1}^{n}\sum_{k=1}^{n-1}X_n Y_n h'(0) a_h^k c(e_h) \xi_{k}\nonumber\\
	&+\frac{2\xi_n^2}{(1+\xi_n^2)^2}X_n Y_n h'(0)c(e_n)
	-\frac{1}{1+\xi_n^2}
	\sum_{j,h=1}^{n}\sum_{l,k=1}^{n-1}
	(X_jY_l+X_lY_j)
	\partial_{x_j}(a_h^k)c(e_h)\xi_l\xi_k
	\nonumber\\
	&-\frac{1}{1+\xi_n^2}
	\sum_{h=1}^{n}\sum_{l,k=1}^{n-1}
	(X_nY_l+X_lY_n)
	a_h^k\left.\partial_{x_n}c(dx_h)\right|_{x_0}
	\xi_l\xi_k
	\nonumber\\
	&+\frac{1}{(1+\xi_n^2)^2}
	\sum_{h=1}^{n}\sum_{l,k=1}^{n-1}
	(X_nY_l+X_lY_n)
	h'(0)a_h^k c(e_h)\xi_l\xi_k.
\end{align}
Consequently,
\begin{align}
	&\left(
	\pi_{\xi_n}^{+}\Biggl(
	\sum_{j=1}^{n}
	\partial_{\xi_j}
	\bigl[\sigma_2(\widetilde{\nabla}^J_{X}\widetilde{\nabla}^J_{Y})\bigr]
	D_{x_j}\bigl[\sigma_{-1}(D_J^{-1})\bigr]
	\Biggr)
	\times\partial_{\xi_n}\sigma_{-3}(D_J^{-3})
	\right)(x_0)\big|_{|\xi'|=1}
	\nonumber\\
	&\sim
	-\frac{\xi_n}{(\xi_n-i)^5(\xi_n+i)^3}
	\sum_{\beta=1}^{n}\sum_{l,\alpha=1}^{n-1}
	(X_nY_l+X_lY_n)h'(0)a_\beta^\alpha
	c(e_n)c(e_\beta)\xi_l\xi_\alpha
	\nonumber\\
	&\quad
	+\frac{2i\xi_n}{(\xi_n-i)^4(\xi_n+i)^3}
	\sum_{j,h,\beta=1}^{n}\sum_{k,\alpha=1}^{n-1}
	(X_jY_n+X_nY_j)
	a_\beta^\alpha\partial_{x_j}(a_h^k)
	c(e_h)c(e_\beta)\xi_k\xi_\alpha
	\nonumber\\
	&\quad
	+\frac{4i\xi_n}{(\xi_n-i)^4(\xi_n+i)^3}
	\sum_{h,\beta=1}^{n}\sum_{k,\alpha=1}^{n-1}
	X_nY_n a_h^k a_\beta^\alpha
	\left.\partial_{x_n}c(dx_h)\right|_{x_0}
	c(e_\beta)\xi_k\xi_\alpha
	\nonumber\\
	&\quad
	-\frac{2\xi_n}{(\xi_n-i)^5(\xi_n+i)^3}
	\sum_{h,\beta=1}^{n}\sum_{k,\alpha=1}^{n-1}
	X_nY_nh'(0)a_h^k a_\beta^\alpha
	c(e_h)c(e_\beta)\xi_k\xi_\alpha
	\nonumber\\
	&\quad
	-\frac{\xi_n(1-3\xi_n^2)}
	{2(\xi_n-i)^5(\xi_n+i)^3}
	X_nY_nh'(0)
	\nonumber\\
	&\quad
	-\frac{1-3\xi_n^2}
	{2(\xi_n-i)^4(\xi_n+i)^3}
	\sum_{j,h=1}^{n}\sum_{l,k=1}^{n-1}
	(X_jY_l+X_lY_j)
	\partial_{x_j}(a_h^k)
	c(e_h)c(e_n)\xi_l\xi_k
	\nonumber\\
	&\quad
	-\frac{1-3\xi_n^2}
	{2(\xi_n-i)^4(\xi_n+i)^3}
	\sum_{h=1}^{n}\sum_{l,k=1}^{n-1}
	(X_nY_l+X_lY_n)a_h^k
	\left.\partial_{x_n}c(dx_h)\right|_{x_0}
	c(e_n)\xi_l\xi_k
	\nonumber\\
	&\quad
	+\frac{(\xi_n-2i)(1-3\xi_n^2)}
	{4(\xi_n-i)^5(\xi_n+i)^3}
	\sum_{h=1}^{n}\sum_{l,k=1}^{n-1}
	(X_nY_l+X_lY_n)h'(0)a_h^k
	c(e_h)c(e_n)\xi_l\xi_k.
\end{align}
Again, since
$\int_{|\xi'|=1}
\xi_{\beta_1}\xi_{\beta_2}\cdots\xi_{\beta_{2d+1}}
\,\sigma(\xi')=0$,
the terms containing odd monomials in $\xi'$
make no contribution after integration over $|\xi'|=1$.
Taking the trace, we therefore obtain
\begin{align}\label{39}
	&\mathrm{Tr} \left(\pi^+_{\xi_n}\Big(\sum_{j=1}^{n}\partial_{\xi_j}\big[\sigma_{2}(\widetilde{\nabla}^J_{X}\widetilde{\nabla}^J_{Y} )\big]
	D_{x_j}\big[\sigma_{-1}(D_{J}^{-1})\big]\Big)\times
	\partial_{\xi_n}\sigma_{-3}(D_{J}^{-3})\right)(x_0)|_{|\xi'|=1}\nonumber\\
	&\sim 2X_nY_nh'(0)\frac{3\xi_n^2-i\xi_n}{(\xi_n-i)^4(\xi_n+i)^3}.
\end{align}
Integrating over \(\xi_n\) and \(|\xi'|=1\), we obtain
\begin{align}\label{39i}
&-i\int_{|\xi'|=1}\int^{+\infty}_{-\infty}
\mathrm{Tr}\Big[\pi^+_{\xi_n}\Big(\sum_{j=1}^{n}\partial_{\xi_j}
\big[\sigma_{2}(\widetilde{\nabla}^J_{X}\widetilde{\nabla}^J_{Y})\big]
D_{x_j}\big[\sigma_{-1}(D_{J}^{-1})\big]\Big)
\times
\partial_{\xi_n}\sigma_{-3}(D_{J}^{-3})\Big](x_0)d\xi_n\sigma(\xi')dx'\nonumber\\
&=\frac{1}{4}X_{n}Y_n\pi h'(0)\Omega_3dx'.
\end{align}
Combining the contributions from (1), (2), and (3), we obtain
\begin{align}\label{41ii}
\widetilde{\Phi}_4
=&-\frac{1}{8} \sum_{v=1}^{n-1} a_v^v X_nY_nh'(0)\pi\Omega_3 dx'-\frac{1 }{12}\sum_{j,v=1}^{n-1} a_v^v X_jY_jh'(0)\pi\Omega_3 dx'\nonumber\\
&-\frac{1}{32} X_nY_n  h'(0) \pi\Omega_3 dx'-\frac{5}{32}\sum_{j=1}^{n-1} X_j Y_j  h'(0) \pi\Omega_3 dx'\nonumber\\
&+\frac{1}{8} \sum_{y,\gamma=1}^{n} \sum_{j=1}^{n-1} \bigg( X_j Y_n + X_n Y_j \bigg) a_y^{\gamma} \partial_{x_{\gamma}}(a_y^j)\pi\Omega_3 dx'\nonumber\\
&-\frac{3}{4}X(Y_n)\pi\Omega_3dx'+ \frac{1}{8} h'(0) g\big(JX,Y\big)\pi\Omega_3dx'\nonumber\\
&+ \frac{1}{8}g\big((\nabla^L_{e_{n}}J)Y,J(X^T)\big)\pi\Omega_3dx'+ \frac{1}{8}g\big((\nabla^L_{e_{n}}J)X,J(Y^T)\big)\pi\Omega_3dx'\nonumber\\
&- \frac{1}{8}g\big((\nabla^L_{X^T}J)Y,e_n\big)\pi\Omega_3dx' - \frac{1}{8}g\big((\nabla^L_{Y^T}J)X,e_n\big)\pi\Omega_3dx'\nonumber\\
&- \frac{3}{8}\sum_{\alpha=1}^{n} g\big((\nabla^L_{e_{\alpha}}J)Y,J(e_\alpha)\big)X_n \pi\Omega_3dx' - \frac{3}{8}\sum_{\alpha=1}^{n} g\big((\nabla^L_{e_{\alpha}}J)X,J(e_\alpha)\big)Y_n \pi\Omega_3dx'.
\end{align}

 {\bf  case c)}~$r=1,~\ell=-4,~k=j=|\alpha|=0$.

By (4.2), we get
\begin{align}\label{61i}
\widetilde{\Phi}_5&=-i\int_{|\xi'|=1}\int^{+\infty}_{-\infty}\mathrm{Tr} [\pi^+_{\xi_n}
\sigma_{1}(\widetilde{\nabla}^J_{X}\widetilde{\nabla}^J_{Y} D_{J}^{-1})\times
\partial_{\xi_n}\sigma_{-4}(D_{J}^{-3})](x_0)d\xi_n\sigma(\xi')dx'\nonumber\\
&=i\int_{|\xi'|=1}\int^{+\infty}_{-\infty}\mathrm{Tr}
[\partial_{\xi_n}\pi^+_{\xi_n}\sigma_{1}( \widetilde{\nabla}^J_{X}\widetilde{\nabla}^J_{Y} D_{J}^{-1})\times
\sigma_{-4}(D_{J}^{-3})](x_0)d\xi_n\sigma(\xi')dx'.
\end{align}
By Lemma 4.3, we have
\begin{align}
	&\sigma_{-4}({D}_{J}^{-3})(x_{0})|_{|\xi'|=1}\nonumber\\
	&=-\frac{1}{4(1+\xi_{n}^{2})^{2}}h'(0)\sum_{\mu=1}^{n}\sum_{\nu=1}^{n-1}a_{\nu}^{\mu}c(e_{\mu})c(e_{n})c(e_{\nu})\nonumber\\
	&-\frac{1}{(1+\xi_{n}^{2})^{3}}h'(0)\sum_{\omega=1}^{n}\sum_{\lambda,k=1}^{n-1}\xi_{k}\xi_{\lambda}a_{\omega}^{\lambda}c(e_{k})c(e_{n})c(e_{\omega})\nonumber\\
	&+\frac{\xi_{n}}{(1+\xi_{n}^{2})^{3}}h'(0)\sum_{k=1}^{n-1}\xi_{k}c(e_{k})+\frac{\xi_{n}^{3}+7\xi_{n}}{(1+\xi_{n}^{2})^{4}}h'(0)\sum_{\omega=1}^{n}\sum_{\lambda=1}^{n-1}\xi_{\lambda}a_{\omega}^{\lambda}c(e_{\omega})\nonumber\\
	&+\frac{3\xi_{n}^{4}+9\xi_{n}^{2}}{(1+\xi_{n}^{2})^{4}}h'(0)c(e_n)+\frac{2}{(1+\xi_{n}^{2})^{3}}\sum_{\alpha,\beta,\omega=1}^{n}\sum_{\lambda=1}^{n-1}\xi_{\lambda}a_{\alpha}^{\beta}a_{\omega}^{\lambda}c(e_{\beta})c[(\nabla^{L}_{e_{\alpha}}J)(\xi^{*})]c(e_{\omega})\nonumber\\
	&+\frac{2\xi_{n}}{(1+\xi_{n}^{2})^{3}}\sum_{\alpha,\beta=1}^{n}a_{\alpha}^{\beta}c(dx_{\beta})c[(\nabla^{L}_{e_{\alpha}}J)(\xi^{*})]c(e_n)-\frac{4}{(1+\xi_{n}^{2})^{3}}\sum_{h=1}^{n}\sum_{j,p=1}^{n-1}\xi_{j}\xi_{p}\partial_{x_{j}}(a_{h}^{p})c(e_{h})\nonumber\\
	&-\frac{4\xi_{n}}{(1+\xi_{n}^{2})^{3}}\sum_{h=1}^{n}\sum_{p=1}^{n-1}\xi_{p}\partial_{x_{n}}(a_{h}^{p})c(e_{h}).
\end{align}

and
\begin{align}
	&\partial_{\xi_n}\pi_{\xi_n}^{+}
	\sigma_1\bigl(
	\widetilde{\nabla}^J_{X}\widetilde{\nabla}^J_{Y} D_J^{-1}
	\bigr)(x_0)\big|_{|\xi'|=1}
	\nonumber\\
	&=-\frac{i}{2(\xi_n-i)^2}
	X_nY_n c(e_n)
	-\frac{1}{2(\xi_n-i)^2}
	\sum_{h=1}^{n}\sum_{p=1}^{n-1}
	X_nY_n a_h^p\xi_p c(e_h)
	\nonumber\\
	&\quad
	-\frac{1}{2(\xi_n-i)^2}
	\sum_{k=1}^{n-1}
	X_kY_n\xi_k c(e_n)
	-\frac{1}{2(\xi_n-i)^2}
	\sum_{l=1}^{n-1}
	X_nY_l\xi_l c(e_n)
	\nonumber\\
	&\quad
	+\frac{i}{2(\xi_n-i)^2}
	\sum_{h=1}^{n}\sum_{p,k=1}^{n-1}
	X_kY_n a_h^p\xi_p\xi_k c(e_h)
	+\frac{i}{2(\xi_n-i)^2}
	\sum_{h=1}^{n}\sum_{p,l=1}^{n-1}
	X_nY_l a_h^p\xi_p\xi_l c(e_h)
	\nonumber\\
	&\quad
	+\frac{i}{2(\xi_n-i)^2}
	\sum_{k,l=1}^{n-1}
	X_kY_l\xi_k\xi_l c(e_n)
	+\frac{1}{2(\xi_n-i)^2}
	\sum_{h=1}^{n}\sum_{p,k,l=1}^{n-1}
	X_kY_l a_h^p\xi_p\xi_k\xi_l c(e_h).
\end{align}
Since the integrals of odd monomials in 
$\xi'$, over $|\xi'|=1$ vanish, the terms with odd powers of $\xi'$ make no contribution. Therefore, a direct computation gives 
\begin{align}
	&\mathrm{Tr}\Bigl[
	\partial_{\xi_n}\pi_{\xi_n}^{+}
	\sigma_1\bigl(
	\widetilde{\nabla}^J_{X}\widetilde{\nabla}^J_{Y} D_J^{-1}
	\bigr)
	\times\sigma_{-4}(D_J^{-3})
	\Bigr](x_0)\big|_{|\xi'|=1}
	\nonumber\\
	&\sim
	h'(0)
	\left(
	\sum_{k,l=1}^{n-1}X_kY_l\xi_k\xi_l-X_nY_n
	\right)
	\Biggl[
	\frac{2(\xi_n^2-4i\xi_n-1)}
	{(\xi_n-i)^5(\xi_n+i)^4}
	-\frac{3i\xi_n+1}
	{2(\xi_n-i)^4(\xi_n+i)^3}
	\sum_{v=1}^{n-1}a_v^v
	\Biggr]
	\nonumber\\
	&\quad
	-\frac{4}{(\xi_n-i)^4(\xi_n+i)^3}
	\sum_{y,\gamma=1}^{n}\sum_{k,p=1}^{n-1}
	(X_kY_n+X_nY_k)
	a_y^\gamma\partial_{x_\gamma}(a_y^p)
	\xi_k\xi_p.
\end{align}
Using the adapted orthonormal frame and spinor
trivialization fixed in Section~3, we express the
derivatives of the coefficients of $J$ in terms of
its covariant derivative.
All the following identities are evaluated at $x_0$.

The Levi-Civita connection satisfies
\begin{align*}
	\nabla^L_{e_a}e_b
	&=\frac{h'(0)}2\delta_{ab}e_n,
	&
	\nabla^L_{e_a}e_n
	&=-\frac{h'(0)}2e_a,
	\\
	\nabla^L_{e_n}e_a&=0,
	&
	\nabla^L_{e_n}e_n&=0,
	\qquad a,b<n.
\end{align*}
Since $J(e_n)=e_n$ throughout the collar and $J$ is
orthogonal, we have
$a_n^k=a_k^n=0$ for $k<n$ and $a_n^n=1$.
Consequently, for $\gamma,k<n$,
\begin{align*}
	(\nabla^L_{e_\gamma}J)e_k
	&=
	\sum_{h=1}^{n-1}
	\partial_{x_\gamma}(a_h^k)e_h
	+\frac{h'(0)}2
	\bigl(a_\gamma^k-\delta_{\gamma k}\bigr)e_n,
	\\
	(\nabla^L_{e_n}J)e_k
	&=
	\sum_{h=1}^{n-1}
	\partial_{x_n}(a_h^k)e_h.
\end{align*}
Here $J(e_\gamma)$ is tangential for $\gamma<n$,
whereas $J(e_n)=e_n$. Taking the corresponding inner
products therefore gives
\begin{equation}\label{eq:partial-a-covariant-J}
	\sum_{\gamma,h=1}^{n}
	a_h^\gamma\partial_{x_\gamma}(a_h^k)
	=
	\sum_{\gamma=1}^{n}
	g\bigl(
	J(e_\gamma),(\nabla^L_{e_\gamma}J)e_k
	\bigr),
	\qquad k<n.
\end{equation}
Using \eqref{eq:partial-a-covariant-J} and performing
the integrations over $\xi_n$ and $|\xi'|=1$, we obtain
\begin{align}
	\widetilde{\Phi}_5
	={}&-\frac14
	\sum_{v=1}^{n-1}
	a_v^v h'(0)X_nY_n\pi\Omega_3\,dx'
	-\frac{25}{32}
	h'(0)X_nY_n\pi\Omega_3\,dx'
	\nonumber\\
	&+\frac1{12}
	\sum_{k,v=1}^{n-1}
	a_v^vX_kY_kh'(0)\pi\Omega_3\,dx'
	+\frac{25}{96}
	\sum_{k=1}^{n-1}
	X_kY_kh'(0)\pi\Omega_3\,dx'
	\nonumber\\
	&+\frac5{12}
	\sum_{\alpha=1}^{n}\sum_{k=1}^{n-1}
	(X_kY_n+X_nY_k)
	g\!\left(J(e_\alpha),(\nabla^L_{e_\alpha}J)e_k\right)
	\pi\Omega_3\,dx'.
\end{align}

Write
$X=X^T+X_n\partial_{x_n}$ and
$Y=Y^T+Y_n\partial_{x_n}$,
where $X^T$ and $Y^T$ are the tangential components.
Set
\[
\tau=\sum_{a=1}^{n-1}a_a^a.
\]
All the following identities are evaluated at $x_0$.
The above connection formulas imply
\begin{align*}
	(\nabla^L_{e_a}J)e_n
	&=\frac{h'(0)}2\bigl(J(e_a)-e_a\bigr),
	\qquad a<n,
	\\
	(\nabla^L_{e_n}J)e_n&=0.
\end{align*}
Hence
\begin{equation}\label{eq:J-normal-contraction}
	\sum_{\alpha=1}^{n}
	g\bigl(
	J(e_\alpha),(\nabla^L_{e_\alpha}J)e_n
	\bigr)
	=
	\frac{h'(0)}2(n-1-\tau).
\end{equation}

Since $J$ is self-adjoint, we also have
\begin{align*}
	&g\bigl((\nabla^L_{X^T}J)Y,e_n\bigr)
	+g\bigl((\nabla^L_{Y^T}J)X,e_n\bigr)
	\\
	&\qquad
	=h'(0)\bigl(
	g(JX^T,Y^T)-g(X^T,Y^T)
	\bigr).
\end{align*}
Moreover, differentiating
$g(JU,JV)=g(U,V)$ covariantly in the normal direction
and using $(\nabla^L_{e_n}J)e_n=0$ gives
\begin{align*}
	&g\bigl((\nabla^L_{e_n}J)Y,JX^T\bigr)
	+g\bigl((\nabla^L_{e_n}J)X,JY^T\bigr)
	=0.
\end{align*}

Substituting these identities and
\eqref{eq:partial-a-covariant-J} into the boundary
contributions, the terms involving $g(JX^T,Y^T)$ cancel.

The coefficient of
\[
\sum_{k=1}^{n-1}(X_kY_n+X_nY_k)
\sum_{\alpha=1}^{n}
g\bigl(
J(e_\alpha),(\nabla^L_{e_\alpha}J)e_k
\bigr)
\]
vanishes, as do the coefficients of
$\tau h'(0)X_nY_n$ and
$\tau h'(0)g(X^T,Y^T)$.
Therefore, for $n=4$, summing the five cases yields
\begin{align*}
	\widetilde{\Phi}
	={}&
	\Bigg[
	h'(0)\left(
	\frac18\sum_{j=1}^{n-1}X_jY_j-X_nY_n
	\right)
	-\frac34X(Y_n)
	\\
	&\qquad
	-\frac16\sum_{j=1}^{n-1}
	\partial_{x_n}(X_jY_j)
	+\frac12\partial_{x_n}(X_nY_n)
	\Bigg]\pi\Omega_3\,dx'.
\end{align*}

Using the relations established in Section~3 between
the coordinate components and $g(X^T,Y^T)$, together
with their normal derivatives, we obtain

\begin{align}
	\widetilde{\Phi}
	=\sum_{\beta=1}^{5}\widetilde{\Phi}_{\beta}
	={}&
	\Biggl[
	-\left(
	\frac1{24}g(X^T,Y^T)+X_nY_n
	\right)h'(0)
	-\frac34X(Y_n)
	\nonumber\\
	&\qquad
	-\frac16\partial_{x_n}
	\bigl(g(X^T,Y^T)\bigr)
	+\frac12\partial_{x_n}(X_nY_n)
	\Biggr]\pi\Omega_3\,dx'.
	\nonumber
\end{align}

So, we are reduced to prove the following.
\begin{thm}\label{thmb2}
Let $(M,g,J)$ be a compact oriented four-dimensional
almost product Riemannian spin manifold with smooth boundary.
Assume that a collar neighborhood of $\partial M$
is identified with $\partial M\times[0,\varepsilon)$,
on which the metric has the form
\[
g=h(x_n)^{-1}g^{\partial M}+dx_n^2,
\qquad h(0)=1,
\]
where $h$ is a smooth positive function and $x_n$
increases in the inward normal direction.
Assume further that
\[
J(\partial_{x_n})=\partial_{x_n}
\]
throughout this collar neighborhood.
Let $\widetilde{\nabla}^{J}$ be the connection defined
in \eqref{eq:connection-J}.
Then, for vector fields $X,Y$ smooth up to the boundary,
the following identity holds:
	\begin{align}\label{b2773}
	&\widetilde{\mathrm{Wres}}\left[
	\pi^+\bigl(
	\widetilde{\nabla}^J_{X}\widetilde{\nabla}^J_{Y} D_J^{-1}
	\bigr)
	\circ\pi^+(D_J^{-3})
	\right]
	\nonumber\\
	&=\frac{4\pi^2}{3}
	\int_M
	\left[
	\operatorname{Ric}(X,Y)-\frac{s}{2}g(X,Y)
	\right]\,d\mathrm{vol}_g
	\nonumber\\
	&\quad
	+2\pi^2\int_M\sum_{l=1}^{4}
	\Bigl\{
	X\!\left[
	g\bigl(Je_l,(\nabla^L_{e_l}J)Y\bigr)
	\right]
	-Y\!\left[
	g\bigl(Je_l,(\nabla^L_{e_l}J)X\bigr)
	\right]
	-g\bigl(Je_l,(\nabla^L_{e_l}J)[X,Y]\bigr)
	\Bigr\}\,d\mathrm{vol}_g
	\nonumber\\
	&\quad
	+\pi^2\int_M g(X,Y)
	\Biggl[
	-\sum_{i,j=1}^{4}R(Je_i,Je_j,e_j,e_i)
	+2\sum_{\nu,j=1}^{4}
	g\bigl(
	(\nabla^L_{e_j}J)e_\nu,
	(\nabla^L_{e_\nu}J)e_j
	\bigr)
	\nonumber\\
	&\hspace{23mm}
	+2\sum_{\nu,j=1}^{4}
	g\Bigl(
	Je_\nu,
	\bigl[
	\nabla^L_{e_j}(\nabla^L_{e_\nu}J)
	-\nabla^L_{\nabla^L_{e_j}e_\nu}J
	\bigr]e_j
	\Bigr)
	+\sum_{\alpha,\nu,j=1}^{4}
	g\bigl(Je_\alpha,(\nabla^L_{e_\alpha}J)e_j\bigr)
	g\bigl(Je_\nu,(\nabla^L_{e_\nu}J)e_j\bigr)
	\nonumber\\
	&\hspace{23mm}
	-\sum_{\nu,j=1}^{4}
	g\bigl(
	(\nabla^L_{e_\nu}J)e_j,
	(\nabla^L_{e_\nu}J)e_j
	\bigr)
	+\sum_{\alpha,\nu,j=1}^{4}
	g\bigl(Je_\alpha,(\nabla^L_{e_\nu}J)e_j\bigr)
	g\bigl((\nabla^L_{e_\alpha}J)e_j,Je_\nu\bigr)
	\Biggr]\,d\mathrm{vol}_g
	\nonumber\\
	&\quad
	+\pi\Omega_3\int_{\partial M}
	\Biggl[
	-\left(
	\frac1{24}g(X^T,Y^T)+X_nY_n
	\right)h'(0)
	-\frac34X(Y_n)
	-\frac16\partial_{x_n}\bigl(g(X^T,Y^T)\bigr)
	\nonumber\\
	&\hspace{40mm}+\frac12\partial_{x_n}(X_nY_n)
	\Biggr]\,dx'.
\end{align}
where $s$ is the scalar curvature.
\end{thm}

\section{Applications in physics and noncommutative geometry}
\label{sec:applications}
The preceding formulas supply local geometric data that may be used in spectral models of gravity and in noncommutative geometry. We give concrete consequences and examples, and distinguish them from proposed extensions. Throughout this section, physical interpretations are Euclidean; a Lorentzian theory or field equations do not follow from the residue identities alone.

\subsection{Tensor responses associated with an almost product splitting}
Set $r=2^m$ and $v_{n-1}=2\pi^m/\Gamma(m)$, where $n=2m$. The orthogonal projections
\begin{equation}
 P_\pm=\frac12(\mathrm{Id}\pm J),\qquad TM=T_+M\oplus T_-M
\end{equation}
define two distinguished distributions. Their failure to be parallel is measured by $\nabla^LJ$; integrability is a separate condition. Since $c(J\xi^\sharp)^2=-|\xi|_g^2\mathrm{Id}$, the principal symbol of $D_J^2$ is the original metric symbol. Thus the twist introduces lower-order geometric couplings without replacing the Riemannian principal metric.

Write $\mathcal G_J(V,W)=\operatorname{Wres}(\widetilde\nabla^J_V\widetilde\nabla^J_W D_J^{-2m})$ on a closed manifold and define
\begin{equation}
 \alpha_J(Z)=\sum_l g(Je_l,(\nabla^L_{e_l}J)Z),\qquad
 q_J=r^{-1}\operatorname{Tr}_{S(TM)}E_J.
\end{equation}
The traced curvature calculation in Section~2 gives
$\operatorname{Tr}F^J(V,W)=\frac r2 d\alpha_J(V,W)$. Consequently Lemma~2.2 yields the useful decomposition
\begin{align}
 \mathcal G_J(V,W)
 ={}&\frac{rv_{n-1}}6\int_M G(V,W)\,d\mathrm{vol}_g
 +\frac{rv_{n-1}}4\int_M d\alpha_J(V,W)\,d\mathrm{vol}_g\nonumber\\
 &+\frac{rv_{n-1}}2\int_M q_Jg(V,W)\,d\mathrm{vol}_g.
 \label{eq:app-decomp}
\end{align}
Here $G=\operatorname{Ric}-\frac s2g$, and $q_J$ is explicitly determined by the trace formula of Section~2. If $\mathcal G_J^\pm(V,W)=\frac12(\mathcal G_J(V,W)\pm\mathcal G_J(W,V))$, then
\begin{align}
 \mathcal G_J^+(V,W)&=rv_{n-1}\int_M\left[\frac16G(V,W)+\frac12q_Jg(V,W)\right]d\mathrm{vol}_g,\\
 \mathcal G_J^-(V,W)&=\frac{rv_{n-1}}4\int_M d\alpha_J(V,W)\,d\mathrm{vol}_g.
\end{align}
This separates a symmetric gravitational response from an antisymmetric response of the splitting. The latter is an exact two-form and therefore carries no nonzero de Rham cohomology class. Exactness does not force its evaluation and integration against arbitrary $V,W$ to vanish. These formulas provide a direct diagnostic in examples: insert tangent vectors from $T_+M$ and $T_-M$ to compare diagonal and mixed responses. Neither response is asserted to be a conserved stress tensor.

A reference case is a Riemannian product $M=M_+\times M_-$, equipped with $J=+\mathrm{Id}$ on the first factor and $J=-\mathrm{Id}$ on the second. Then $\nabla^LJ=0$ and $\alpha_J=0$. In the product spin description, $D_J$ changes the sign of one factor Dirac operator; the two factor summands anticommute, so $D_J^2=D^2$. Hence $E_J=s\mathrm{Id}/4$ and
\begin{equation}
 \mathcal G_J(V,W)=rv_{n-1}\int_M
 \left[\frac16\operatorname{Ric}(V,W)+\frac1{24}s g(V,W)\right]d\mathrm{vol}_g.
\end{equation}
In particular, mixed factor vectors give zero. A flat product torus gives zero for the entire functional. These are concrete consistency tests for numerical or symbolic calculations with nonparallel splittings.

For perturbation theory with fixed $g$, let $A$ be a skew-adjoint bundle endomorphism and set
\begin{equation}
 J_t=e^{tA}Je^{-tA},\qquad \dot J_0=[A,J].
 \label{eq:app-variation}
\end{equation}
This preserves orthogonality and $J_t^2=\mathrm{Id}$ exactly. On a flat product torus with constant $J$, a smooth periodic $A(x)$ produces $\partial_kJ_t=t[\partial_kA,J]+O(t^2)$. Substitution into the explicit $q_J$ and $d\alpha_J$ formulas gives a controlled expansion in derivatives of $A$. Constant $A$ leaves a parallel splitting and produces no such derivative terms. Comparing variable and constant $A$ therefore isolates the part of the residue sensitive to spatial variation of the splitting.

\subsection{Low-energy spectral actions and constrained variations}
The spectral action principle \cite{CCSA} suggests using heat coefficients to organize possible effective couplings. For a positive self-adjoint Laplace-type operator $P=\Delta^\nabla+E$ on a closed four-manifold, write
\begin{equation}
 \operatorname{Tr}(e^{-tP})\sim t^{-2}a_0(P)+t^{-1}a_2(P)+a_4(P)+\cdots.
\end{equation}
With the sign convention used in this paper, the standard heat formula \cite{PBG} gives
\begin{align}
 a_0(P)&=(4\pi)^{-2}\int_M\operatorname{Tr}(\mathrm{Id})\,d\mathrm{vol}_g,\\
 a_2(P)&=(4\pi)^{-2}\int_M\operatorname{Tr}\left(\frac s6\mathrm{Id}-E\right)d\mathrm{vol}_g.
 \label{eq:app-a2}
\end{align}
For a sufficiently regular rapidly decreasing cutoff $f$, this produces
\begin{equation}
 \operatorname{Tr} f(P/\Lambda^2)
 \sim\Lambda^4F_4a_0(P)+\Lambda^2F_2a_2(P)+f(0)a_4(P)+\cdots,
\end{equation}
where $F_4=\int_0^\infty u f(u)du$ and $F_2=\int_0^\infty f(u)du$. For $P=D_J^2$ under the self-adjointness condition specified below, the four-dimensional spinor rank is four and
\begin{equation}
 a_2(D_J^2)=(4\pi)^{-2}\int_M\left(\frac23s-4q_J\right)d\mathrm{vol}_g.
 \label{eq:app-action}
\end{equation}
Thus the previously computed $q_J$ gives an explicit candidate coupling between curvature and the almost product structure. For $J=\mathrm{Id}$ this density becomes $-s/3$, consistent with the adopted heat convention. Its sign alone is not a stability criterion for a proposed gravitational theory.

There is an essential analytic distinction. On a closed manifold, integration by parts gives
\begin{equation}
 D_J^*=D_J+c(B_J),\qquad B_J=\sum_j(\nabla^L_{e_j}J)e_j.
\end{equation}
Therefore $D_J$ is formally self-adjoint when $B_J=0$; this is not automatic for an arbitrary almost product structure. If $B_J\ne0$, a natural symmetric elliptic replacement is
\begin{equation}
 \mathcal D_J=\frac12(D_J+D_J^*)=D_J+\frac12c(B_J).
\end{equation}
It has a self-adjoint realization on a closed manifold. Its square has the same principal symbol as $D_J^2$ but generally different connection and endomorphism. Equations \eqref{eq:app-a2} and the spectral expansion then apply to $\mathcal D_J^2$ after recomputing those lower-order data; one must not simply reuse $E_J$. The algebraic residue calculation for $D_J^2$ itself does not require this self-adjointness restriction.

A concrete next step is to study a truncated action built from $a_0$ and $a_2$, adding $a_4$ when curvature-squared effects are required. At fixed metric, use \eqref{eq:app-variation} to vary the integrated scalar density. If an unconstrained variational derivative is written $\delta S=\int_M\operatorname{tr}(T\delta J)d\mathrm{vol}_g$, then
\begin{equation}
 \delta S=\int_M\operatorname{tr}([J,T]A)\,d\mathrm{vol}_g.
\end{equation}
Stationarity means that this pairing vanishes for every skew-adjoint $A$, with boundary terms handled separately when present. This is a practical constrained Euler--Lagrange prescription. A simultaneous metric variation must also account for the varying spin connection and the orthogonality constraint. The bilinear residue \eqref{eq:app-decomp} is not, without that calculation, the metric derivative of the candidate action.

\subsection{Unchanged spectral distance and additional curvature data}
For $a\in C^\infty(M)$, multiplication by $a$ satisfies
\begin{equation}
 [D_J,a]=c(J\operatorname{grad}a),\qquad
 [\mathcal D_J,a]=c(J\operatorname{grad}a),\qquad
 \|[\mathcal D_J,a]\|=\sup_M|da|_g.
 \label{eq:app-distance}
\end{equation}
On a closed connected spin manifold, the self-adjoint elliptic operator $\mathcal D_J$ has compact resolvent, and these bounded commutators define a spectral triple on $L^2(S(TM))$. Equation \eqref{eq:app-distance} yields
\begin{equation}
 \sup\{|a(x)-a(y)|:a\in C^\infty(M),\ \|[\mathcal D_J,a]\|\leq1\}=d_g(x,y).
\end{equation}
Thus this construction preserves the Riemannian distance, whereas its lower-order spectral curvature data may depend on $J$. It provides a specific setting in which distance and additional geometric responses can be studied separately. If the residue is computed using $\mathcal D_J$ rather than $D_J$, the correction explained above remains necessary. Real structures, orientability and additional axioms must be verified independently when required.

The spectral Einstein functionals of \cite{DL} provide a route toward genuinely noncommutative analogues. A workable program would first choose a regular spectral triple admitting a residue calculus, then specify an involution on an appropriate differential module, and finally construct a compatible self-adjoint Dirac-type operator. The metric and Einstein residues could then be tested against the commutative product and perturbation examples above. A tangent-bundle tensor $J$ alone does not supply these algebraic data, so no automatic extension to an arbitrary noncommutative algebra is claimed.

There is also a concrete almost commutative starting point. Let $(\mathcal A_F,H_F,D_F)$ be finite spectral data and let $\Gamma$ be spinor chirality. For the even-dimensional operator $\mathcal D_J$, which anticommutes with $\Gamma$, consider
\begin{equation}
 \mathcal D=\mathcal D_J\otimes1+\Gamma\otimes D_F,
 \qquad
 \mathcal D^2=\mathcal D_J^2\otimes1+1\otimes D_F^2.
\end{equation}
The unfluctuated heat trace factorizes. If $N_F=\dim H_F$, its first two coefficients are
\begin{equation}
 a_0(\mathcal D^2)=N_Fa_0(\mathcal D_J^2),\qquad
 a_2(\mathcal D^2)=N_Fa_2(\mathcal D_J^2)
 -\operatorname{Tr}_{H_F}(D_F^2)a_0(\mathcal D_J^2).
\end{equation}
This exhibits explicitly how finite internal mass data modify the volume coefficient while multiplying the geometric contribution by $N_F$. Inner fluctuations and gauge or scalar couplings would require a new heat calculation; the factorized example does not establish a particular particle model.

\subsection{Boundary probes and comparison with gravitational boundary actions}
Let $\mathcal B_{22}(X,Y)$ and $\mathcal B_{13}(X,Y)$ denote the two boundary residue functionals of Sections~3 and~4. Under their common hypotheses, their interior parts agree, hence
\begin{equation}
 \mathcal B_{22}(X,Y)-\mathcal B_{13}(X,Y)
 =\int_{\partial M}(\Phi_{22}(X,Y)-\Phi_{13}(X,Y)).
\end{equation}
This gives an observable of the boundary factorization. Theorem~\ref{thmb2} is particularly useful because its explicit $J$-dependent boundary terms cancel when $J\nu=\nu$ throughout the collar.

For the inward unit normal $\nu=\partial_{x_n}$, fix
$K(U,V)=g(\nabla^L_U\nu,V)$ and $H=\operatorname{tr}K$. The collar metric gives
\begin{equation}
 K=-\frac12h'(0)g^{\partial M},\qquad H=-\frac32h'(0)
\end{equation}
in dimension four. Writing the boundary density in Theorem~\ref{thmb2} as $\pi\Omega_3 b_{13}(X,Y)dx'$, one obtains
\begin{align}
 b_{13}(X,Y)={}&\left(\frac1{36}g(X^T,Y^T)+\frac23X_nY_n\right)H
 -\frac34X(Y_n)\nonumber\\
 &-\frac16\partial_{x_n}g(X^T,Y^T)
 +\frac12\partial_{x_n}(X_nY_n).
 \label{eq:app-boundary}
\end{align}
Two choices of probes make its meaning more explicit. First, let $X,Y$ be tangential throughout the collar and extend them by Levi-Civita parallel transport along $\nu$. Then $X_n=Y_n=0$ and $\partial_{x_n}g(X,Y)=0$, so
\begin{equation}
 \mathcal B_{13}^{\partial}(X,Y)=\frac{\pi\Omega_3}{36}
 \int_{\partial M}H g(X,Y)\,dx'.
\end{equation}
Second, set $X=Y=\nu$ throughout the collar. All derivative terms vanish and
\begin{equation}
 \mathcal B_{13}^{\partial}(\nu,\nu)=\frac{2\pi\Omega_3}{3}
 \int_{\partial M}H\,dx'.
\end{equation}
These formulas isolate mean-curvature responses in tangential and normal channels with different coefficients. They supply direct tests on a warped cylinder; for a product collar $h'(0)=0$, both responses vanish for these extensions. For general extensions the derivative terms in \eqref{eq:app-boundary} survive even when $H=0$, so the boundary value of a probe alone is insufficient to determine the functional.

The appearance of $\int_{\partial M}H$ invites comparison with the gravitational boundary action \cite{GHBoundary}. However, the general functional depends on inserted vector fields, their normal derivatives and the chosen factorization. It is not itself the Gibbons--Hawking--York action or a boundary stress tensor. To obtain such an interpretation one must first choose a scalar action and a variational problem, specify the fixed boundary data, and verify cancellation of unwanted normal variations. Moreover, the residue in Boutet de Monvel's calculus is not a choice of self-adjoint elliptic boundary condition. A comparison with a heat trace or effective action on a bounded region requires a compatible operator realization. These steps identify precise further calculations rather than physical conclusions already established by the present theorems.

\section*{ Acknowledgements}
This work is sponsored by Natural Science Foundation of Xinjiang Uygur Autonomous Region 2024D01C341 and supported by the National Natural Science Foundation of China 11771070.
 The authors thank the referee for his (or her) careful reading and helpful comments.

\section*{Declarations}

\begin{itemize}
	\item  {\bf Funding} Sponsored by Natural Science Foundation of Xinjiang Uygur Autonomous Region 2024D01C341 and the National Natural Science Foundation of China 11771070.
\item \textbf{Conflict of interest}
The authors declare that they have no conflict of interest.
	\item \textbf{Ethics approval and consent to participate}
	Not applicable.
	\item {\bf Consent for publication} All the authors agreed to publish this research.
	\item {\bf Author contribution} All authors contributed to the study conception and design. Material preparation, data collection and analysis were performed by JH and YW. The first draft of the manuscript was written by JH and all authors commented on previous versions of the manuscript. All authors read and approved the final manuscript.
\end{itemize}

\end{document}